\documentclass[11pt]{amsart}

\usepackage{graphicx}

\newtheorem{theorem}{Theorem}[section]
\newtheorem{corollary}[theorem]{Corollary}
\newtheorem{proposition}[theorem]{Proposition}
\newtheorem{lemma}[theorem]{Lemma}

\newtheoremstyle{defn}{12pt}{12pt}{}{}{\bfseries}{.}{ }{}
\theoremstyle{defn}
\newtheorem{definition}[theorem]{Definition}

\theoremstyle{remark}
\newtheorem*{notation}{Notation}
\newtheorem{remark}[theorem]{Remark}
\newtheorem*{example}{Example}
\newtheorem*{claim}{Claim}

\numberwithin{figure}{section}

\input xy
\xyoption{all}

\def\Z{\mathbb{Z}}
\def\R{\mathbb{R}}
\def\RR{\mathcal{R}}
\def\C{\mathbb{C}}
\def\Q{\mathbb{Q}}
\def\alg{\hat{\mathcal{A}}}
\def\d{d}
\def\dd{d}
\def\dsft{\d_{\operatorname{SFT}}}
\def\dstr{\d_{\operatorname{str}}}
\def\I{\mathcal{I}}
\def\J{\mathcal{J}}
\def\CC{\mathcal{C}}
\def\A{\mathcal{A}}
\def\F{\mathcal{F}}
\def\comm{\hat{\mathcal{A}}^{\operatorname{comm}}}
\def\cyc{\hat{\mathcal{A}}^{\operatorname{cyc}}}
\def\tocyc{^{\operatorname{cyc}}}
\def\tocomm{^{\operatorname{comm}}}
\def\co{\colon\thinspace}
\def\sftb#1{\{#1\}}
\def\Id{\operatorname{Id}}

\def\secspace{\vspace{33pt}}
\begin{document}

\title{Rational Symplectic Field Theory for Legendrian knots}
\author{Lenhard Ng}
\address{Mathematics Department, Duke
University, Durham, NC 27708}
\email{ng@math.duke.edu}

\begin{abstract}
We construct a combinatorial invariant of Legendrian knots in standard
contact three-space. This invariant, which encodes rational relative
Symplectic Field Theory and extends contact homology, counts
holomorphic disks with an arbitrary number of positive punctures. The
construction uses ideas from string topology.
\end{abstract}

\maketitle

\section{Introduction}

The theory of Legendrian knots plays a key role in contact and symplectic
topology and has recently shown surprising connections to low
dimensional topology; see \cite{EtLegsurvey} for a survey of the subject.
A key breakthrough in the study of Legendrian
knots, and symplectic topology generally, was the
introduction of Gromov-type holomorphic-curve techniques in the
1990s. This led in particular to the development of Legendrian contact
homology, outlined by Eliashberg and Hofer \cite{Eli98} and fleshed
out famously by Chekanov \cite{Ch02} for standard contact $\R^3$
and later by others in more general setups (e.g.,
\cite{EES05a,EES07,NT,Sab03}). Besides applications to contact
topology, Legendrian contact homology has
been closely linked to standard knot theory (e.g., \cite{Ngframed}).

Contact homology is part of a much larger construction,
Symplectic Field
Theory (SFT), which was introduced by Eliashberg, Givental, and Hofer
about a decade
ago \cite{EGH}. The relevant portion of the SFT package for our
purposes is a filtered theory for contact manifolds whose first order
comprises contact
homology.  Somewhat more precisely, while contact homology counts holomorphic
disks in the symplectization of a contact manifold with exactly one
positive boundary puncture, SFT counts holomorphic curves with arbitrarily many
positive punctures.

In the ``closed'' case (in the absence of a Legendrian or
Lagrangian boundary condition), SFT is now fairly well understood,
both algebraically and analytically, and has produced a number of
spectacular applications in symplectic topology; see, e.g.,
\cite{EliSFT} and references therein. However, in the ``relative'' case
that is the focus of this paper, much less is currently understood. In
particular, technical problems involving bubbling of holomorphic
curves have thus far prevented a formulation of SFT with Legendrian
boundary condition even for the basic case of standard contact
$\R^3$. The development of contact homology for Legendrian knots
involves two steps, a fairly easy proof that $\d^2=0$ and a more
difficult invariance proof; it has proven surprisingly difficult to
extend this to a reasonable algebraic setup for Legendrian SFT that
even satisfies $d^2=0$, not to mention invariance.




In this paper, we will give an algebraic formulation, \`a la Chekanov
\cite{Ch02},
of Legendrian SFT for standard contact $\R^3$; this allows
us to skirt the analytical issues that would usually beset the proofs
of $d^2=0$ and invariance. We note that we present not the full
Legendrian SFT, which would
consider holomorphic curves of arbitrary genus and possibly marked
points and gravitational descendants, but ``rational'' SFT, which
counts only holomorphic disks.\footnote{This is a slight misuse of the
  term ``rational'' since we do
  not count genus-$0$ surfaces with more than one boundary component.}

The technique that we use to overcome the bubbling
problems comes from string topology \cite{CS}.
Cieliebak and Latschev \cite{CiL}, motivated by work of Fukaya, have
developed a program for using string topology to deal with
compactification issues in Legendrian SFT; see especially the
appendix to \cite{CiL} jointly written with Mohnke. The program currently has significant unresolved technical issues, but one can avoid these issues in the case of $\R^3$ by using the combinatorial approach we employ here.
On a related note, we remark without proof that a separate approach to Legendrian SFT, based on the cluster homology of
  Cornea and Lalonde \cite{CL}, seems in the $\R^3$ case to yield the same theory as ours, or at least the commutative quotient that we call $(\comm,\d)$.

We now outline the mathematical content of this paper. In
Section~\ref{sec:lsft}, we associate to any Legendrian knot in
standard contact $\R^3$ a filtered version of a familiar structure
from algebra, a curved dg-algebra, which itself is a type of a
curved $A_\infty$ algebra. Our particular filtered curved
dg-algebra, which we call the \textit{LSFT algebra} $(\alg,\d)$ of the
Legendrian knot, takes the following form: $\alg$ is the tensor
algebra over $\Z$ freely generated by two
generators $p,q$ for each Reeb chord, along with one more generator
$t$ and its inverse $t^{-1}$ essentially encoding the homology of the knot.
The map $\d$ on $\alg$ is a derivation
\[
\d = \dsft + \dstr,
\]
where $\dsft$ is an ``SFT differential'' obtained by counting rational
holomorphic curves in the symplectization $\R\times\R^3$ with boundary on
the Lagrangian cylinder over the knot and boundary punctures
approaching Reeb chords at $\pm\infty$ in the distinguished $\R$
direction, and $\dstr$ is a ``string differential'' encoding a string
cobracket operation that glues trivial holomorphic strips to broken
closed strings on the knot. The Hamiltonian that produces the SFT
differential lives naturally in the quotient $\cyc$ of $\alg$ by
cyclic permutations but acts on $\alg$ as well.

The string differential is a necessary
correction that accounts for the aforementioned bubbling and ensures
a result analogous to $\d^2=0$. More precisely, we have the following
two main results.


\begin{theorem}[see Theorem~\ref{thm:d2}]
The algebra $(\alg,\d)$ associated to a Legendrian knot is a curved
dg-algebra; that is, there is an element
$F_\d$ of $\alg$ such that $\d^2(x) = F_\d x - x
F_\d$ for all $x\in\alg$.
\label{thm:main1}
\end{theorem}


\begin{theorem}[see Theorem~\ref{thm:invariance}]
$(\alg,\d)$ is invariant under restricted Legendrian isotopies.
\label{thm:main2}
\end{theorem}

\noindent
Here ``restricted'' is a minor technical condition (see
Definition~\ref{def:restricted}) that we conjecture
can be removed, but that in any case can still be used to produce an
invariant of Legendrian knots under arbitrary Legendrian isotopies;
see Corollary~\ref{cor:invariance}. It is possible that we can remove
the ``restricted'' condition if we allow arbitrary equivalences of
curved $A_\infty$ algebras rather than the specific equivalences of
LSFT algebras defined in Section~\ref{ssec:alg}, but we do not pursue
this point in this paper.

The LSFT algebra has a filtration whose associated graded object, in
the bottom filtration level, is Legendrian contact homology (cf.\
Remark~\ref{rem:LCH}).
Theorems~\ref{thm:main1} and~\ref{thm:main2} contain Chekanov's
$d^2=0$ and invariance results for contact homology
(Corollary~\ref{cor:LCH}).

One possible and desirable application of Legendrian SFT would be the
construction of invariants of Legendrian knots that do not vanish for
stabilized knots, which in some sense comprise ``most'' Legendrian
knots. This could produce invariants of topological knots (which can
be viewed as Legendrian knots modulo stabilization) and transverse
knots (Legendrian knots modulo one particular stabilization), among
other things. Contact homology famously vanishes under
stabilization \cite{Ch02}, but it was hoped for some time that
Legendrian SFT would not. Unfortunately, rational Legendrian SFT, as
constructed in this paper, also loses all interesting information
under stabilization; see Appendix~B. There is some hope that one could
apply rational SFT to the double of a Legendrian knot \cite{NT} to obtain an
interesting invariant, but this is unclear as yet. We note that the
contact homology of the double of a stabilized knot contains no
information \cite{Ngthesis}, but rational SFT may encode significantly
more information.

We remark that we develop the theory over $\Z$, and a fair amount of
work throughout the paper is devoted to keeping track of
signs. In particular, we include an appendix that computes all
possible sign rules, in some suitable sense, and shows that they are
all equivalent.
However, the entire theory works over $\Z/2$ as well as $\Z$,
with the notable exception of invariance for cyclic and commutative
complexes (Proposition~\ref{prop:quasi}), and the reader may find it
easier to ignore all signs.

In this paper, we omit discussion of the relation between our algebraic
version of rational Legendrian SFT and the
more general, more geometric string-topology version, though we may
return to this topic in the future. We also postpone concrete applications of
the Legendrian SFT formalism presented here, such as the construction
of an $L_\infty$
structure on cyclic Legendrian contact homology, to a future paper.

The main results of this paper are contained in
Section~\ref{sec:lsft}. Their proofs, some of which involve a
discussion of a rudimentary version of string topology, occupy
Sections~\ref{sec:string} (for Theorem~\ref{thm:main1}) and
\ref{sec:invariance} (for Theorem~\ref{thm:main2}). Appendices A and B
deal with sign choices and triviality for stabilized knots, respectively.




\subsection*{Acknowledgments}

I would like to give significant thanks to Yasha Eliashberg, whose
many conversations with me
about candidates for Legendrian SFT played a key role in the
present work. In addition, the crucial
catalyst for this paper was the September 2007 workshop ``Towards
Relative Symplectic Field Theory'' sponsored by the American Institute
of Mathematics, the NSF, the CUNY Graduate Center, and the Stanford
Mathematical Research Center. I am deeply indebted to all of the
workshop's participants, particularly Mohammed Abouzaid, Fr\'ed\'eric Bourgeois, Kai Cieliebak, Tobias Ekholm,
John Etnyre, Eleny Ionel, Janko Latschev, and Josh Sabloff, for their ideas and
suggestions, and to Mikhail Khovanov for a separate illuminating
conversation. The combinatorial version of Legendrian SFT presented
here was largely formulated in discussions at the AIM workshop. I also thank the referee for helpful corrections and suggestions.
This work is partially supported by the following NSF grants:
DMS-0706777, FRG-0244663, and CAREER grant DMS-0846346.

\secspace

\section{The SFT Invariant}
\label{sec:lsft}

In this section, we describe the algebraic object to be associated to
a Legendrian knot, the LSFT algebra, and state
the main ``$d^2=0$'' and invariance results, though their proofs are
deferred to Sections~\ref{sec:string} and~\ref{sec:invariance}.
The LSFT algebra is a special case of a familiar construction from
homological algebra, the curved dg-algebra, whose salient features we
review in Section~\ref{ssec:curved}. We then present the
definition of an LSFT algebra in Section~\ref{ssec:alg}, followed by
a combinatorial definition for the LSFT algebra associated to the
$xy$ projection of a Legendrian knot in Section~\ref{ssec:comb}. In
Section~\ref{ssec:cyccomm}, we discuss two quotient invariants derived from the
LSFT algebra, the cyclic and commutative complexes.

\subsection{Algebraic setup: curved dg-algebras}
\label{ssec:curved}

Throughout this section and the paper, we use the convention that the
commutator on a graded associative algebra is
$[x,y] = xy-(-1)^{|x||y|} yx.$

\begin{definition}
A \textit{curved dg-algebra} consists of a triple
\label{def:cdg}
$(\A,d,F)$, where:
\begin{itemize}
\item $\A$ is a graded associative algebra over $\Z$;
\item $d :\thinspace \A\to\A$ is a derivation, i.e.,
$d(xy) = (dx)y+(-1)^{|x|}x(dy)$, and $d$ lowers degree by $1$;
\item $F$ is a degree $-2$ element of $\A$ (the \textit{curvature}) for
  which $dF = 0$;
\item for all $x\in\A$, $d^2(x) = [F,x]$.
\end{itemize}
A \textit{filtered curved dg-algebra} is a curved dg-algebra with a
descending filtration of subalgebras
\[
\A = \F^0\!\A \supset \F^1\!\A \supset \F^2\!\A \supset \cdots
\]
with respect to which $d$ is a filtered derivation and $F \in \F^1\!\A$.
\end{definition}

\begin{remark}
Curved dg-algebras have been studied extensively in the literature,
though sometimes under other names, e.g., CDG-algebra
\cite{Positselski} and $Q$-algebra \cite{Schwarz} (note however that
the standard definition involves an algebra over a field rather than
over $\Z$).
In particular, a curved dg-algebra is essentially a special case of a curved (or
``weak'') $A_\infty$ algebra. A curved $A_\infty$ algebra is
a graded vector space $V$ with multilinear maps $m_n
:\thinspace V^{\otimes n} \to V$ of degree $n-2$
 for all $n \geq 0$, satisfying the curved $A_\infty$ relations
\[
\sum_{i+j+k=n} (-1)^{i+j(1+i+|a_1|+\cdots+|a_i|)}
m_{i+k+1}(a_1,\dots,a_i,m_j(a_{i+1},\dots,a_{i+j}),a_{i+j+1},\dots,a_{i+j+k})
= 0
\]
for $n \geq 0$. Except for the aforementioned discrepancy in base
ring, a curved dg-algebra is a curved
$A_\infty$ algebra where $m_n = 0$ for $n \geq 3$; we then have $m_0
= F$, $m_1(a_1) = da_1$, and $m_2(a_1,a_2) = a_1a_2$, and the curved
$A_\infty$ relations become the relations in
Definition~\ref{def:cdg}, along with multiplicative associativity.
For comparison, note that
a usual $A_\infty$ algebra is a curved $A_\infty$ algebra with $m_0 =
0$, while a usual dg-algebra has $m_n=0$ for all $n \neq 1,2$.
\end{remark}

A special case of morphisms of curved $A_\infty$ algebras is the
following.

\begin{definition}
A \textit{morphism} of curved dg-algebras is a map $(\varphi,\alpha)
:\thinspace (\A,d,F) \to (\A',d',F')$, where:
\begin{itemize}
\item
$\varphi :\thinspace \A \to \A'$ is a graded algebra map;
\item
$\alpha$ is a degree $-1$ element of $\A'$;
\item
$d'\varphi(\cdot) = \varphi d(\cdot) + [\alpha,\varphi(\cdot)]$;
\item
$F' = \varphi(F) +d\alpha + \alpha^2$.
\end{itemize}
\label{def:fcdgmorph}
A \textit{filtered morphism} of filtered curved dg-algebras is a
morphism for which $\varphi$ respects the filtration and $\alpha \in
\F^1\!\A'$.
\end{definition}

It is easy to check that a composition of morphisms is a
morphism, where we define
$(\varphi',\alpha') \circ (\varphi,\alpha) =
(\varphi'\circ\varphi,\alpha'+\varphi'(\alpha))$. There is an identity
morphism $(\Id,0)$, and if
$(\varphi,\alpha)$ is a morphism for which $\varphi$ is an
isomorphism, then $(\varphi^{-1},-\varphi^{-1}\alpha)$ provides an
inverse to $(\varphi,\alpha)$.

We can now define chain homotopy and homotopy equivalence in the usual
way. We state the definitions for filtered curved dg-algebras; there
is an obvious analogue in the unfiltered case.

\begin{definition}
Two filtered morphisms of filtered curved dg-algebras
$(\varphi,\alpha),(\varphi',\alpha') :\thinspace (\A,d,F) \to
(\A',d',F')$ are \textit{chain homotopic} if $\alpha = \alpha'$ and
there exists a filtered $\Z$-module map $H :\thinspace
\A'\to\A'$ of degree $1$ such that
\[
\varphi-\varphi' = Hd'+d'H.
\]
\label{def:htpyequiv}
A filtered morphism $(\varphi,\alpha) :\thinspace (\A,d,F) \to
(\A',d',F')$ is a \textit{homotopy equivalence} if there exists a
filtered morphism $(\varphi',\alpha') :\thinspace (\A',d',F') \to
(\A,d,F)$ such that $(\varphi',\alpha') \circ (\varphi,\alpha)$ and
$(\varphi,\alpha) \circ (\varphi',\alpha')$ are each chain homotopic
to the identity morphism.
\end{definition}

We can now state a preliminary version of the main result of this
paper; see Theorems~\ref{thm:d2} and~\ref{thm:invariance} for the
precise statements.

\begin{theorem}
Rational SFT gives a map from Legendrian knots in $\R^3$ modulo
Legendrian isotopy to filtered curved dg-algebras modulo homotopy
equivalence. \label{thm:fcdgalg}
\end{theorem}

Because of the curvature term $F$, a curved dg-algebra does not typically
comprise a complex. One can produce a complex and thus homology from a
filtered curved dg-algebra in several ways. See
Remark~\ref{rem:assocgraded} for discussion of the associated graded
complex, and Section~\ref{ssec:cyccomm} for the cyclic and commutative complexes.

\subsection{Algebraic setup: LSFT algebras}
\label{ssec:alg}

The invariant we associate to a Legendrian knot is a particular type
of filtered curved dg-algebra that we term an LSFT algebra.
Besides being a specialization of the construction in the previous
section, our definition of LSFT algebra generalizes (and
contains) Chekanov's DGAs and stable tame
isomorphisms from Legendrian contact homology.

Underlying an LSFT algebra is a (based) tensor algebra $\A$ over
$\Z$ generated by
$q_1,\dots,q_n,p_1,\dots,p_n,t,t^{-1}$; this is noncommutative and has
sole relations $t\cdot t^{-1} = t^{-1} \cdot t = 1$. We consider
$q_1,\dots,q_n,p_1,\dots,p_n$ to be distinguished generators that are
included in the data of the LSFT algebra, where we view $q_i$ and
$p_i$ as being paired
together for $i=1,\ldots,n$, but the indices $1,\ldots,n$ can be
permuted without changing $\A$. Each generator of $\A$ is
$\Z$-graded with $|q_i|+|p_i|=-1$ for all $i$, and $|t| = -|t^{-1}| =
-2r$ for some $r\in\Z$; this grading induces a grading on $\A$.

There is a filtration
\[
\A = \F^0\!\A \supset \F^1\!\A \supset \F^2\!\A \supset \cdots,
\]
where $\F^k\!\A$ is generated by words containing at least
$k$ $p$'s. (Note that $\F^k\!\A = (\F^1\!\A)^k$.) We will sometimes write
$O(p^k)$ to denote an element of $\F^k\!\A$ (or $\F^k\!\alg$, defined
below), and $x \equiv y \pmod{p^k}$ for $x = y + O(p^k)$.

Let $\alg$ be the ``$p$-adic completion'' of $\mathcal{A}$ consisting
of possibly infinite sums $\sum_{k=0}^\infty z_k$ with $z_k \in
\F^k\!\mathcal{A}$ for all $k$. That is, $\alg$ includes infinite sums
in $\mathcal{A}$ as long as for each $k$, all but finitely many terms
in the sum do not lie in $\F^k\!\mathcal{A}$. Then $\alg$ inherits from
$\mathcal{A}$ the structure of a graded algebra with filtration $\alg
= \F^0\!\alg \supset \F^1\!\alg \supset \cdots$.


\begin{definition}
An \textit{LSFT algebra} is a filtered graded tensor algebra $\alg
=\Z\langle q_1,\dots,q_n,p_1,\dots,p_n,t,t^{-1}\rangle$, as above, with a
derivation\footnote{As in the previous section, a derivation is a
  $\Z$-linear map
$d\co\alg\to\alg$ such that $d(xy) =
(dx)y+(-1)^{|x|} x(dy)$ for all $x,y\in\alg$
for which $x$ is of pure degree. Note that, for an LSFT algebra, $d$
necessarily satisfies $d(1)=0$ and
$d(t^{-1}) = -t^{-1} \cdot d(t) \cdot t^{-1}$.}
$\d\co \alg\to\alg$ satisfying the following conditions:
\label{def:lsftalg}
\begin{enumerate}
\item
$\d$ has degree $-1$ and preserves the filtration;

\item
$\d(t) \in \F^1\!\alg$;

\item \label{cond3}
there is an element $F_{\d} \in \F^1\!\alg$, the
\textit{curvature} of $\d$, such that
$\d^2 x = [F_{\d},x]$ for all $x\in\alg$.
\end{enumerate}
We denote an LSFT algebra by $(\A,\d)$, omitting the curvature
$F_{\d}$, which is uniquely determined by $\d$.
\end{definition}

Condition (\ref{cond3}) ensures that $dF_{\d} = 0$, since
$[F_{\d},\d x] = \d^2(\d x) = \d(\d^2x) = [\d F_{\d},x]+[F_{\d},\d x]$ for
all $x\in\alg$;
thus an LSFT algebra is a filtered curved dg-algebra in the sense of
Section~\ref{ssec:curved}.

\begin{remark}[The Chekanov--Eliashberg DGA]
Given a curved dg-algebra $(\alg,d,F)$, one can consider
the complex given by the
associated graded object $\oplus_{i=0}^\infty \F^i\!\alg/\F^{i+1}\!\alg$
with the induced differential. In the case when $(\alg,d)$ is an
LSFT algebra generated by $q_1,\dots,q_n,p_1,\dots,p_n,t,t^{-1}$, the
$i=0$ summand $(\F^0\!\alg/\F^1\!\alg,\d)$ of the associated graded
complex is generated by $q_1,\dots,q_n,t,t^{-1}$, with $\d(t)=\d(t^{-1})=0$.

This quotient
$(\F^0\!\alg/\F^1\!\alg,\d)$ is essentially Chekanov's
differential graded
algebra (usually abbreviated DGA), also formulated by Eliashberg, that encodes Legendrian contact homology.
Indeed, it will be clear
from the definition of $\d$ in Section~\ref{ssec:comb} that the
differential on $\F^0\!\alg/\F^1\!\alg$, and in fact the entire
associated graded object, counts precisely the same holomorphic disks
as contact homology, namely disks with exactly one positive puncture. It should be noted, however, that
$(\F^0\!\alg/\F^1\!\alg,\d)$ is not precisely the same as the Chekanov
DGA; see
Remark~\ref{rem:LCH} below.

%


\label{rem:assocgraded}
\end{remark}

\begin{notation}
We will sometimes want to treat the $q$'s and $p$'s
together, and will use $s$ to denote any $q_j$ or $p_j$ (or sometimes
$t^{\pm 1}$ as well). The $q$'s and
$p$'s are paired together, and we use $^*$ to denote the pairing; that
is, write $p_j^* = q_j$, $q_j^* = p_j$. We reserve $w$ to mean a word
in the $q$'s, $p$'s, and $t^{\pm 1}$.

If $s$ is a $q_j$ or $p_j$, then define $\sftb{s,s^*}$ to be $+1$ if
$s$ is a $p$ and $-1$ if $s$ is a $q$; this is a special case of the
SFT bracket to be defined in Section~\ref{ssec:sftb}.
\end{notation}

We next define a notion of equivalence between LSFT algebras, which is
a special case of homotopy equivalence between filtered curved
dg-algebras (see Proposition~\ref{prop:htpyequiv} below). To do this,
we introduce a specific family of curved dg-morphisms between LSFT algebras.

\begin{definition}
Let $\alg$ be an LSFT algebra (without its differential).
\label{def:elem}
An \textit{elementary
  automorphism} of $\alg$ is a grading-preserving algebra automorphism
$\phi$ of $\alg$ of one of the following forms:

\begin{enumerate}
\item
$\phi(q_i) = \pm t^{\alpha_i} q_i t^{\beta_i}$, $\phi(p_i) = \pm
t^{\gamma_i} p_i t^{\delta i}$,
\label{def:elem0}
$\phi(t) = t$, $\phi(t^{-1}) = t^{-1}$ for some integers
$\alpha_i,\beta_i,\gamma_i,\delta_i$;

\item
for some $j$, $\phi(q_i) = q_i$ for all $i \neq j$, $\phi(p_i) = p_i$
for all $i$, $\phi(t^{\pm 1}) = t^{\pm 1}$,
\label{def:elem1}
and
\[
\phi(q_j) = q_j + v + u
\]
where $v \in \alg$ does not involve $q_j$ and $u \in \F^1\!\alg$;

\item
for some $j$, $\phi(p_i) = p_i$ for all $i \neq j$, $\phi(q_i) = q_i$
for all $i$, $\phi(t^{\pm 1}) = t^{\pm 1}$,
\label{def:elem2} and
\[
\phi(p_j) = p_j + v + u
\]
where $v \in \alg$ does not involve $p_j$ and $u \in \F^2\!\alg$;

\item
$\phi(q_i)=q_i$ and $\phi(p_i)=p_i$ for all $j$, \label{def:elem3} and
$\phi(t) = t + v$ for some $v \in \F^1\!\alg$.
\end{enumerate}
In the last three cases, we say that the elementary automorphism is
\textit{supported at} the generator of $\alg$ on which it is
nontrivial: $q_j$ for (\ref{def:elem1}), $p_j$ for (\ref{def:elem2}),
$t$ for (\ref{def:elem3}).
\end{definition}

Implicit in the above definition is the following fact.

\begin{lemma}
Each of the maps in Definition~\ref{def:elem} is invertible.
\end{lemma}

\begin{proof}
Maps of type (\ref{def:elem0}) in the statement of
Definition~\ref{def:elem} are obviously invertible.
Next consider a map $\phi$ of type (\ref{def:elem1}).
It suffices to show that $\phi$ is invertible if either $v=0$ or
$u=0$, since in the general case, $\phi = \phi_1 \circ \phi_2$, where
$\phi_1$, $\phi_2$ are supported on the same $q_j$ and $\phi_1(q_j) =
q_j+v$, $\phi_2(q_j) = q_j + \phi_1^{-1}(u)$ (note $\phi_1^{-1}(u) \in
\F^1\!\alg$). Now if $u=0$, then
$\phi$ is clearly invertible: $\phi^{-1}(q_j) = q_j - v$. If $v=0$,
define $\psi \co \alg \to \alg$ by $\psi(x) = \phi(x)-x$ for all
$x\in\alg$; then $\psi$ increases filtration level by $1$, and
\[
\phi^{-1}(q_j) = q_j - u + \psi(u) - \psi \circ \psi(u) +
\psi \circ \psi \circ \psi(u) - \cdots
\]
gives the inverse for $\phi$.

The same proof works for a map of type (\ref{def:elem2}).

Finally, suppose that $\phi$ is of type (\ref{def:elem3}). We can define
\[
\phi(t^{-1}) = t^{-1} - t^{-1} \cdot v \cdot t^{-1} + t^{-1} \cdot v
\cdot t^{-1} \cdot v \cdot t^{-1} -
\cdots
\]
and then $\phi$ is an algebra map on $\alg$; also, $\phi$ is
invertible for the same reason as in case (\ref{def:elem1}).
\end{proof}

\begin{definition}
We say that LSFT algebras $(\alg,\d)$ and $(\alg,\d')$ are related by a
\textit{basis change} if there is a sequence of elementary
automorphisms of $\alg$ sending $\d$ to $\d'$.
\label{def:basischange}
\end{definition}

\noindent
We remark that the composition $\phi$ of elementary automorphisms in
Definition~\ref{def:basischange} yields a curved dg-morphism
$(\phi,0)$ between the curved dg-algebras given by $(\alg,\d)$ and
$(\alg',\d')$.

If $(\alg,\d)$ and $(\alg,\d')$ are related by a basis change, then the
quotient differential graded algebras $(\F^0\!\alg/\F^1\!\alg,\d)$ and
$(\F^0\!\alg/\F^1\!\alg,\d')$ are
related by a tame isomorphism in the sense of Chekanov (see \cite{Ch02,ENS}
for the precise definition). Note that on
the quotient level, any basis change fixes $t$ and $t^{-1}$.

We need two more operations on LSFT algebras, gauge change and
stabilization.

\begin{definition}
We say that LSFT algebras $(\alg,\d)$ and $(\alg,\d')$ are related by a
\textit{gauge change} if there exists $\alpha\in\F^1\!\alg$ with
$|\alpha|=-1$ such that
\begin{equation}
\d'(x) = \d(x) + [\alpha,x]
\label{eq:gauge}
\end{equation}
for all $x\in\alg$.
\label{def:gauge}
\end{definition}

It is easy to check that if $(\alg,\d)$ is an LSFT algebra, then
(\ref{eq:gauge}) defines an LSFT algebra $(\alg,\d')$ with
$F_{\d'} = F_{\d} + \d(z) + z^2$.
Note that a gauge change is nothing more than a curved dg-morphism of the form $(\Id,\alpha)$.

\begin{remark}
Our notion of a gauge change coincides with the standard algebraic
notion of changing by an inner derivation. One can view the derivation
$\d$ on $\alg$ as an element of the Hochschild cohomology
$\textit{HH}^1(\alg)$; two derivations on $\alg$ related by gauge change
represent the same element of $\textit{HH}^1(\alg)$.
\end{remark}

Finally, we define stabilization.

\begin{definition}
\label{def:stab}
Let $(\alg,\d)$ be an LSFT algebra
generated by $q_1,\dots,q_n,p_1,\dots,p_n,t,t^{-1}$. The degree-$i$
\textit{(algebraic) stabilization} of $(\alg,\d)$ is the LSFT algebra
$(S_i\alg,\d)$ generated by
$q_1,\dots,q_n,p_1,\dots,p_n,t,t^{-1},q_a,q_b,p_a,p_b$, where $q_a,q_b,p_a,p_b$
are four new generators with $|q_a|=|q_b|+1=-1-|p_a|=-|p_b|=i$, and $\d$ is
defined on $S_i\alg$ by extending the existing derivation by
\[
\d(q_a) = q_b,
\qquad
\d(q_b) = [F_\d,q_a],
\qquad
\d(p_b) = p_a,
\qquad
\d(p_a) = [F_\d,p_b].
\]
If $(S_i\alg,\d)$ is a stabilization of $(\alg,\d)$, then we say that
$(\alg,\d)$ is a destabilization of $(S_i\alg,\d)$.
\end{definition}

On $\F^0\!\alg/\F^1\!\alg$, this definition reduces to Chekanov's notion of
stabilization for DGAs. The following
definition then generalizes
Chekanov's stable tame isomorphism.

\begin{definition}
Two LSFT algebras are \textit{equivalent} if they are related by some
finite sequence of basis changes, gauge changes, stabilizations, and
destabilizations.
\end{definition}

\begin{proposition}
An equivalence of LSFT algebras is a homotopy equivalence of filtered
curved dg-algebras.
\label{prop:htpyequiv}
\end{proposition}

\begin{proof}
Since basis changes and gauge changes are isomorphisms of the
underlying algebra, it is easy to check that they are homotopy
equivalences in the
sense of Definition~\ref{def:htpyequiv}.
It thus suffices to show that
stabilization is a homotopy equivalence as well.

Let $(\alg,\d)$ be an LSFT algebra with stabilization $(S_i\alg,\d)$. Let
$\iota :\thinspace \alg \to S_i\alg$ and $\pi :\thinspace S_i\alg
\to \alg$ denote the usual inclusion and projection maps, where $\pi$
projects away any word involving the four additional generators
$q_a,q_b,p_a,p_b$. Then $(\iota,0) :\thinspace (\alg,\d,F_\d) \to
(S_i\alg,\d,F_\d)$ and $(\pi,0) :\thinspace (S_i\alg,\d,F_\d) \to
(\alg,\d,F_\d)$ are morphisms of filtered curved dg-algebras, and it
is clear that $(\pi,0) \circ (\iota,0) = (\Id,0)$.

As for $(\iota,0)
\circ (\pi,0)$, define a $\Z$-linear map $H :\thinspace
S_i\alg \to S_i\alg$ by its action on words $w$:
\[
H(w) = \begin{cases}
(-1)^{|w_1|} w_1 q_a w_2 & \text{if $w = w_1 q_b w_2$ for words
  $w_1,w_2$ with $w_1\in\A$} \\
(-1)^{|w_1|} w_1 p_b w_2 & \text{if $w = w_1 p_a w_2$ for words $w_1,w_2$
  with $w_1 \in \A$} \\
0 & \text{if $w\in\A$ or $w = w_1 q_a w_2$ or $w = w_1 p_b w_2$ for} \\
& \qquad \text{words
  $w_1,w_2$ with $w_1\in\A$}.
\end{cases}
\]
The proof is complete once we check that $H$ is a homotopy between
the identity and $\iota \circ \pi$, a fact that we defer to the ensuing lemma.
\end{proof}

\begin{lemma}
On $S_i\alg$, we have $\Id_{S_i\alg} - \iota \circ \pi = H \circ
\d + \d \circ H$.
\label{lem:htpy}
\end{lemma}

\begin{proof}
It suffices to check
\begin{equation}
w - \iota \circ \pi(w) = (H \circ \d)(w) + (\d \circ H)(w)
\label{eq:htpy}
\end{equation}
for all words $w$ in $S_i\A$. If $w\in\A$, both sides of
(\ref{eq:htpy}) are $0$. Otherwise, the left hand side of
(\ref{eq:htpy}) is $w$.
If $w = w_1 q_a w_2$ for $w_1 \in \A$, then
\[
(H \circ \d)(w) + (\d \circ H)(w) = (H \circ \d)(w_1 q_a w_2) =
w_1 q_b w_2 = w;
\]
if $w = w_1 q_b w_2$ for $w_1 \in \A$, then
\begin{eqnarray*}
(H \circ \d)(w) + (\d \circ H)(w) &=& H\left(
(\d w_1) q_b w_2 + (-1)^{|w_1|} w_1[F_{\d},q_a] w_2 \right. \\
&& \left. + (-1)^{|w_1|+i+1}
w_1 q_b (\d w_2) \right) + (-1)^{|w_1|} \d(w_1 q_a w_2) \\
&=& w.
\end{eqnarray*}
The cases $w = w_1 p_a w_2$ and $w = w_1 p_b w_2$ for $w_1 \in\A$ are
similar.
\end{proof}


\subsection{Combinatorial description of the invariant}
\label{ssec:comb}

\begin{figure}
\centerline{
\includegraphics[height=0.75in]{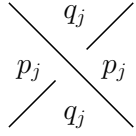}
}
\caption{
Labeling the quadrants at a crossing $R_j$ in the $xy$ projection
$\pi_{xy}(\Lambda)$ by $p_j$ and $q_j$.
}
\label{fig:pq}
\end{figure}

Let $\Lambda$ be a Legendrian knot in $\R^3$ with the standard contact
structure $\ker(dz-y\,dx)$, that is, a knot everywhere tangent to the
contact structure. In this
section, we associate an LSFT algebra to $\Lambda$. A generic knot
$\Lambda$ has finitely many Reeb chords $R_1,\dots,R_n$. To each Reeb
chord $R_j$, we assign two indeterminates $q_j,p_j$.
Let $\pi_{xy}(\Lambda)$ be the knot diagram given by projecting $\Lambda$ to the
$xy$ plane; then the crossings of $\pi_{xy}(\Lambda)$ are the Reeb chords of
$\Lambda$, and the four quadrants at each crossing can be labeled with a $q$
or a $p$ as shown in Figure~\ref{fig:pq}. We also fix two points
$\ast,\bullet$ on $\Lambda$, neither of which lies at an endpoint of a
Reeb chord; the LSFT algebra will depend on the choices of
$\ast,\bullet$, though the equivalence class of the LSFT algebra will not.

Recall that $\Lambda$ has two classical invariants $tb$ and $r$. The
Thurston--Bennequin number $tb(\Lambda)$ is the writhe of the knot
diagram $\pi_{xy}(\Lambda)$. The rotation number $r(\Lambda)$ is the
Whitney index of $\pi_{xy}(\Lambda)$. More precisely, if $\gamma\co
[a,b]\to\R^2$ is any immersed path, then define $r(\gamma) \in\R$ to
be the number of counterclockwise revolutions made by the unit tangent
vector $\gamma'(t)/|\gamma'(t)|$ around $S^1$ as $t$ goes from $a$ to
$b$; $\pi_{xy}(\Lambda)$ is a closed immersed path and $r(\Lambda) =
r(\pi_{xy}(\Lambda)) \in \Z$.

We now construct the LSFT algebra $\alg$ associated to
$(\Lambda,\ast,\bullet)$. This is generated by
$q_1,\dots,q_n,p_1,\dots,p_n,t,t^{-1}$, with grading as follows. For
each $j=1,\dots,n$, there is a unique path $\gamma_j$ along
$\pi_{xy}(\Lambda)$ beginning at the overcrossing of crossing $R_j$,
ending at the undercrossing of $R_j$, and not passing through
$\ast$. If we assume the crossings of $\pi_{xy}(\Lambda)$ are
transverse, then $r(\gamma_j)$ is neither an integer nor a
half-integer. Define
\begin{eqnarray*}
|q_j| &=& \lfloor 2r(\gamma_j) \rfloor \\
|p_j| &=& \lfloor -2r(\gamma_j) \rfloor = -1-|q_j| \\
|t| &=& -2r(\Lambda) \\
|t^{-1}| &=& 2r(\Lambda).
\end{eqnarray*}
We note that the gradings for the $q$'s and $t$ are the same as in
Legendrian contact homology.

When considering signs in the theory, we will often draw an arrow
alongside a section of $\Lambda$; such an arrow is understood to
correspond to a sign $\pm 1$, namely $+1$ if the arrow agrees with the
given orientation of $\Lambda$, and $-1$ if it disagrees. In this
vein, we have the following easy result.

\begin{figure}
\centerline{
\includegraphics[height=0.75in]{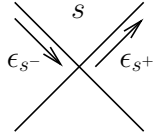}
}
\caption{
Signs at a quadrant $s$.
}
\label{fig:degsign}
\end{figure}

\begin{lemma}
Let $s$ be a $p$ or $q$, corresponding to a corner at a crossing of
$\pi_{xy}(\Lambda)$. Define the signs $\epsilon_{s^-},\epsilon_{s^+}
\in\{\pm 1\}$ to be the orientations along the sides of $s$ relative
to the orientation of $\Lambda$, as shown in
Figure~\ref{fig:degsign}.
\label{lem:degsign}
Then $(-1)^{|s|} = \epsilon_{s^-} \epsilon_{s^+}$.
\end{lemma}

\noindent In the language of \cite{ENS}, $s$ is ``coherent''
($\epsilon_{s^-}\epsilon_{s^+}=1$) if and only
if $|s|$ is even.

\begin{figure}
\centerline{
\includegraphics[height=1.75in]{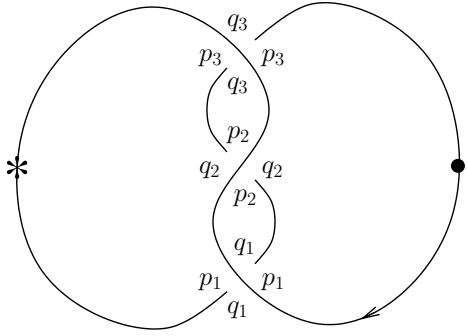}
}
\caption{
The $xy$ projection of a Legendrian unknot $\Lambda_0$.
Note that this is the usual Legendrian unknot after a Reidemeister II move.
}
\label{fig:unknot2}
\end{figure}

\begin{example}
Through this section and Section~\ref{sec:string},
we will use the Legendrian knot $\Lambda_0$ depicted in
Figure~\ref{fig:unknot2} as a running example. Here the gradings are
given by
\[
|q_2|=2, \quad |q_1|=|q_3|=1, \quad
|t|=0, \quad |p_1|=|p_3|=-2, \quad |p_2|=-3.
\]
This agrees with the fact that $p_1,q_2,p_3$ are coherent while
$q_1,p_2,q_3$ are not.
\end{example}

We define the derivation $\d$ on $\alg$ as the sum of two derivations
$\dsft+\dstr$, where $\dsft$ is the ``SFT differential'' and $\dstr$
is the ``string differential''.

For any $k \geq 1$, let $D^2_k$ denote the unit disk $\{|z|\leq 1\}
\subset\C$ minus $k$ fixed points $*_1,\dots,*_k$ on the boundary,
ordered sequentially in counterclockwise order. The punctures divide
the boundary $\partial D^2_k$ into $k$ arcs denoted by $(\partial
D^2_k)_1,\dots,(\partial D^2_k)_k$, where $(\partial D^2_k)_i$ is the
portion of $\partial D^2_k$ between $*_i$ and $*_{i+1}$ (or between
$*_k$ and $*_1$ if $i=k$).

\begin{definition}
For any $s_1,\dots,s_k$ where $k \geq 1$ and each
$s$ is a $q$ or $p$, let $\Delta(s_1,\dots,s_k)$ denote the set of all
orientation-preserving immersions $f\co (D^2_k,\partial D^2_k) \to
(\R^2,\Lambda)$, up to domain reparametrization, such that $f(\partial D^2_k)
\subset \Lambda$ and $f$ sends neighborhoods of the boundary punctures
$\ast_1,\dots,\ast_k$ to quadrants labeled $s_1, \dots, s_k$ in
succession.
\label{def:disks}
\end{definition}

We will call the quadrants described in Definition~\ref{def:disks},
labeled by $s_1,\dots,s_k$, the \textit{corners} of $f$.
Note that $\Delta(s_1,\dots,s_k)$ is unchanged by cyclic permutation
of the $s$'s. We also have the following ``index formula''.

\begin{lemma}
Suppose that $f \in \Delta(s_1,\dots,s_k)$, and let $\alpha$ be the
number of times $f(\partial D^2_k)$ passes through $\ast$, counted
according to the orientation of $\Lambda$. Then
\label{lem:index}
\[
|s_1| + \dots + |s_k| - 2\alpha \, r(\Lambda) = -2.
\]
\end{lemma}

\begin{proof}
For each $s_j$, define $\gamma_{s_j}$ to be the path in
$\pi_{xy}(\Lambda)$ given by $\gamma_{k_j}$ if $s_j =
q_{k_j}$ and $-\gamma_{k_j}$ (i.e., $\gamma_{k_j}$
with the opposite
orientation) if $s_j = p_{k_j}$. Also define $\gamma_{f,j}$ to be the
image in $\pi_{xy}(\Lambda)$ of $f|_{(\partial D^2_k)_j}$. Then
\[
\gamma_{s_1} \cup \gamma_{f,1} \cup \gamma_{s_2} \cup
\gamma_{f,2} \cup \cdots \cup \gamma_{s_k} \cup \gamma_{f,k}
\]
represents a closed loop in $\pi_{xy}(\Lambda)$ (more precisely, the
projection of a closed loop in $\Lambda$) wrapping around
$\pi_{xy}(\Lambda)$ $\alpha$ times. It follows that
\[
\sum_{j=1}^k r(\gamma_{s_j}) + \sum_{j=1}^k r(\gamma_{f,j}) = \alpha
\, r(\Lambda).
\]
Now if $\theta_j$ is the angle (between $0$ and $\pi$) determined by
the image of $f$ at $*_j$, then $r(\gamma_{s_j}) = n/2-\theta_j/(2\pi)$ for
some integer $n$, and thus $|s_j| = 2r(\gamma_{s_j})-\theta_j/\pi$. On
the other hand, since $f$ is an immersed disk,
$\sum_{j=1}^k (r(\gamma_{f,j}) + \theta_j/(2\pi)) = 1$. It
follows that
\[
\sum |s_j| = \sum \left(2r(\gamma_{s_j})-\theta_j/\pi\right) = 2\alpha
\, r(\Lambda) - \sum \left( 2r(\gamma_{f,j})-\theta_j/\pi \right) =
2\alpha \,r(\Lambda) - 2,
\]
as desired.
\end{proof}

\begin{figure}
\centerline{
\includegraphics[height=0.75in]{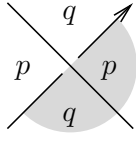}
}
\caption{
Orientation signs for corners. The two unshaded corners are given the
sign $+1$, the shaded corners $-1$. The arrow indicates the
orientation of the knot and ensures that each crossing can be uniquely
viewed as this local picture.
}
\label{fig:pqsigns}
\end{figure}

\begin{figure}
\centerline{
\includegraphics[height=0.75in]{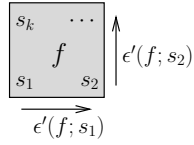}
}
\caption{
The immersed disk $f \in\Delta(s_1,\dots,s_k)$ contributes the term
$\epsilon s_2\dots s_k$ to $\dsft(s_1^*)$, where $\epsilon$ is the
product of: $\epsilon'(f;s_1)$ shown here; the orientation signs for
the $k$ corners $s_1,\ldots,s_k$; and $\sftb{s_1,s_1^*}$.
The $\epsilon'(f;s_2)$ sign will be used in the
proof of Lemma~\ref{lem:hamsign}.
}
\label{fig:hamsign}
\end{figure}

For each map $f \in \Delta(s_1,\dots,s_k)$, we can
define a word $w(f;s_1) \in \A$ by
\[ w(f;s_1) = t^{\alpha_1} s_2 t^{\alpha_2} s_3 \dots t^{\alpha_{k-1}} s_k
t^{\alpha_k},
\]
where $\alpha_i$ is the number of times $f|(\partial D^2_k)_i$
passes through $*$, counted according to the orientation of
$\Lambda$. We also associate a sign $\epsilon(f;s_1) \in \{\pm 1\}$ to $f$
as follows. Each quadrant of a crossing of $\pi_{xy}(\Lambda)$ can be
given an orientation sign according to Figure~\ref{fig:pqsigns}. For
each of the $k$ corners of $f$, we thus obtain a sign
$\epsilon_i(f)$. Further define a sign $\epsilon'(f;s_1)$ to be $+1$
if the image of $f|_{[*_1,*_2]} \subset \pi_{xy}(\Lambda)$, oriented
from $f(*_1)$ to $f(*_2)$, has the same orientation as
$\pi_{xy}(\Lambda)$, and $-1$ if it has the opposite
orientation. Finally, define
\[
\epsilon(f;s_1) = \epsilon'(f;s_1) \epsilon_1(f) \cdots \epsilon_k(f).
\]
See Figure~\ref{fig:hamsign}.

\begin{example}
Consider the bigon $f$ in Figure~\ref{fig:unknot2} with
corners at $p_2$ and $q_3$, which can be considered as an element of
$\Delta(p_2,q_3)$ and of $\Delta(q_3,p_2)$. The orientation signs of
both corners are $-1$. If we consider $f \in \Delta(p_2,q_3)$, then
$\epsilon(f;p_2) = \epsilon'(f;p_2) = 1$; if we consider $f\in
\Delta(q_3,p_2)$, then $\epsilon(f;q_3) = \epsilon'(f;q_3) = -1$.
\end{example}

The following observation will be useful in Section~\ref{ssec:ham}.

\begin{lemma}
Any two diagonally-opposite corners at a crossing have opposite
orientation signs. Also,
if $s,s^*$ denote consecutive corners at a crossing, and $s$ lies
counterclockwise from $s^*$, then the product of the orientation signs
of $s$ and $s^*$ is $\sftb{s,s^*}$ (recall that this is $+1$ if $s$ is
a $p$, $-1$ if $s$ is a $q$).
\label{lem:orsigns}
\end{lemma}

\begin{definition}
Define the SFT differential on $\alg$ by
\begin{align*}
\dsft(q_i) &= \sum_{f\in\Delta(p_i)} \epsilon(f;p_i) w(f) \\
\dsft(p_i) &= - \sum_{f\in\Delta(q_i)} \epsilon(f;q_i) w(f) \\
\dsft(t) &= \dsft(t^{-1}) = 0,
\end{align*}
where $\Delta(s) = \coprod
\Delta(s,s_2\dots,s_k)$ is the set of all immersed disks with
a corner at $s$ (i.e., over all possible $k$ and $s_2,\dots,s_k$).
An immersed disk with multiple corners at $s$ contributes multiple
times to $\dsft(s)$.
Extend $\dsft$ to all of $\alg$ as a derivation.
\label{def:dsft}
\end{definition}

\noindent
It is possible that $\dsft(q_i)$ or $\dsft(p_i)$ may be an infinite
sum, but it will always be a sum in the $p$-adic completion
$\alg$; see the discussion of $h$ in Section~\ref{sec:string}.

We note that $\dsft$ preserves the $p$ filtration on $\alg$. This is a
consequence of a basic area estimate originally due to Chekanov.
Define a height function on the
$p$'s and $q$'s as follows: let $h(p_j)$ be the length of the Reeb
chord $R_j$ (i.e., the difference in the $z$ coordinates of its
endpoints), and let $h(q_j) = -h(p_j)$.

\begin{lemma}
If $\Delta(s_1,\dots,s_k)$ is nonempty, then $\sum_{j=1}^k h(s_j) > 0$.
\label{lem:stokes}
\end{lemma}

\begin{proof}
Since $dz = y\,dx$ along $\Lambda$, it is easy to show from Stokes'
Theorem that $\sum h(s_j)$ is the area of an immersed-disk element of
$\Delta(s_1,\dots,s_k)$. See \cite{Ch02,ENS}.
\end{proof}

\begin{lemma}
$\dsft$ has degree $-1$ and preserves the $p$ filtration on $\alg$.
\label{lem:hdegree}
\end{lemma}

\begin{proof}
The fact that $\dsft$ has degree $-1$ follows from Lemma~\ref{lem:index}.
Since $h(q_j) < 0$ and $h(p_j) > 0$ for all $j$,
Lemma~\ref{lem:stokes} implies that any term in $\dsft(p_j)$ must
contain a $p$, and hence that $\dsft$ preserves the $p$
filtration.
\end{proof}

\begin{example}
For $\Lambda_0$, we have
\begin{align*}
\dsft(p_1) &= -p_2 \\
\dsft(q_1) &= 1 + t + p_3 q_2 + q_2 p_3 t \\
\dsft(p_2) &= -p_1p_3 - p_3tp_1 \\
\dsft(q_2) &= -q_1 + q_3 \\
\dsft(p_3) &= p_2 \\
\dsft(q_3) &= 1+t+q_2p_1 + tp_1q_2 \\
\dsft(t) &= \dsft(t^{-1}) = 0.
\end{align*}
Note that $\dsft^2 \neq 0$, a fact that remains true even if we
quotient by cyclic permutations of words. This is an example of the bubbling problem mentioned in the Introduction.
\end{example}

We next define the string differential $\dstr$. For each Reeb chord
$R_j$ of $\Lambda$, write $R_j^+,R_j^-$ for the endpoints of $R_j$,
with the Reeb vector field flowing from $R_j^-$ to $R_j^+$ (i.e.,
$R_j^+$ has the greater $z$ coordinate). View $q_j$ and $p_j$ as
the line segment $R_j$, oriented from $R_j^-$ to $R_j^+$ for $q_j$ and
from $R_j^+$ to $R_j^-$ for $p_j$. Let $\RR$ denote the set
of Reeb-chord endpoints $\{R_1^+,R_1^-,\dots,R_n^+,R_n^-\}$.



Let $\Gamma$ be the set of embedded paths $\gamma \co [0,1] \to
\Lambda$ such that
$\gamma^{-1}(\RR)$ is finite and $\gamma'(\tau) \neq 0$ whenever
$\gamma(\tau) \in \RR$. If $\gamma\in\Gamma$ and $\gamma(\tau)\in\RR$,
then we can define signs
$\epsilon_1(\gamma;\tau),\epsilon_2(\gamma;\tau),\epsilon(\gamma;\tau)$
as follows: $\epsilon_1(\gamma;\tau)$ is $+1$ if $\gamma(\tau) =
R_i^-$ and $-1$ if $\gamma(\tau) = R_i^+$;
$\epsilon_2(\gamma;\tau)$ is the sign of the orientation of $\gamma$
near $\tau$, relative to the orientation of $\Lambda$ there;
and $\epsilon(\gamma;\tau) =
\epsilon_1(\gamma;\tau) \epsilon_2(\gamma;\tau)$.
Define a map $\delta \co \Gamma \to \A$ as follows:
for each $\tau$ such that $\gamma(\tau) \in \RR$,
define
\[
\tilde{\delta}(\gamma;\tau) = \begin{cases}
q_i p_i & \text{if $\gamma(\tau) = R_i^+$} \\
p_i q_i & \text{if $\gamma(\tau) = R_i^-$},
\end{cases}
\]
and then set
\[
\delta(\gamma) = \sum_{\tau\in
  \gamma^{-1}(\RR),~\tau\neq 0,1} \epsilon(\gamma;\tau)
\tilde{\delta}(\gamma;\tau).
\]

We can now define the string portion $\dstr$ of the differential. Let
$s$ denote one of the $q_i$ or $p_i$. Define $s^+,s^-$ as follows: if
$s = p_i$, then $s^\pm = R_i^\pm$; if $s = q_i$, then $s^\pm =
R_i^\mp$. Recall that we are given two distinct points
$\ast,\bullet\in\Lambda$. There are uniquely defined (up to
reparametrization) injective paths
$\gamma_s^\pm$ in $\Lambda$ that begin at $\bullet$, end at $s^\pm$,
and do not pass through $\ast$. Note that $\gamma_s^+ =
\gamma_{s^*}^-$ and $\gamma_s^- = \gamma_{s^*}^+$.

\begin{figure}
\centerline{
\includegraphics[height=1in]{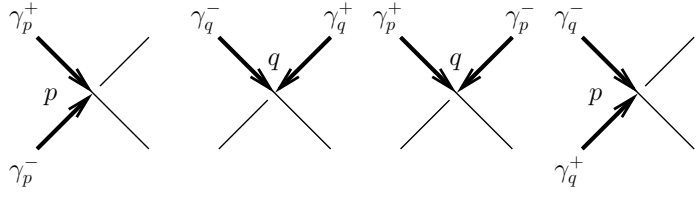}
}
\caption{
Four cases, left to right: $p$ has holomorphic capping paths; $q$
holomorphic; $p$ antiholomorphic; $q$ antiholomorphic. Note that the
middle two diagrams are identical, as are the outer two.
}
\label{fig:holcorner}
\end{figure}

We distinguish two cases: if the
quadrant at $R_i$ in $\pi_{xy}(\Lambda)$ determined by the ends of
$\gamma_s^{\pm}$ is labeled by $s$, we say $s$ \textit{has holomorphic
  capping paths}; if it is labeled by $s^*$, we say $s$
\textit{has antiholomorphic capping paths}. Equivalently, $s$ has
holomorphic capping paths if and only if $\gamma_s^-$ approaches the
crossing in $\pi_{xy}(\Lambda)$ to the right of $\gamma_s^+$. See
Figure~\ref{fig:holcorner}. Note that $s$ has holomorphic capping
paths if and only if $s^*$ has antiholomorphic capping paths.

\begin{definition}
Define the string differential on $\alg$ as follows. If $s$ is a $p$
or $q$ with holomorphic capping paths,
\[
\dstr(s) = \delta(\gamma_s^-) \cdot s + (-1)^{|s|} s \cdot
\delta(\gamma_s^+);
\]
if $s$ is a $p$ or $q$ with antiholomorphic capping paths,
\[
\dstr(s) = \delta(\gamma_s^-) \cdot s + (-1)^{|s|} s \cdot
\delta(\gamma_s^+) + \eta(s) s \cdot s^* \cdot s
\]
where
\[
\eta(s) = \begin{cases}
+1 & \text{if $s=p$ and $\gamma_p^-$ is oriented like $\Lambda$} \\
-1 & \text{if $s=p$ and $\gamma_p^-$ is oriented unlike $\Lambda$} \\
-1 & \text{if $s=q$ and $\gamma_q^-$ is oriented like $\Lambda$} \\
+1 & \text{if $s=q$ and $\gamma_q^-$ is oriented unlike $\Lambda$}.
\end{cases}
\]
Furthermore, $\Lambda$ itself can be viewed as a union of two
injective paths $\gamma_{\Lambda,1},\gamma_{\Lambda,2}$ where
$\gamma_{\Lambda,1}$ begins at $\bullet$ and ends at $\ast$,
$\gamma_{\Lambda,2}$ begins at $\ast$ and ends at $\bullet$, and each
path follows the orientation of $\Lambda$; then set
\[
\dstr(t) = \delta(\gamma_{\Lambda,1}) \cdot t + t \cdot \delta(\gamma_{\Lambda,2})
\]
and $\dstr(t^{-1}) = - t^{-1} \cdot \dstr(t) \cdot t^{-1}$.
Extend $\dstr$ to all of $\alg$ as a derivation.
\label{def:dstr}
\end{definition}

Note that $\dstr$ is well defined since $\dstr(t\cdot
t^{-1}) = \dstr(t^{-1}\cdot t) = 0$.

\begin{figure}
\centerline{
\includegraphics[height=1.5in]{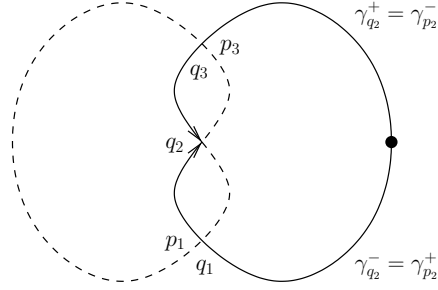}
}
\caption{
Capping paths $\gamma_{q_2}^{\pm},\gamma_{p_2}^{\pm}$ of $q_2,p_2$ for
the Legendrian knot from Figure~\ref{fig:unknot2}.
}
\label{fig:unknot3}
\end{figure}

\begin{example}
For $\Lambda_0$, the capping paths for
$p_2,q_2$ are depicted in Figure~\ref{fig:unknot3}, leading to
$\dstr(q_2) = -q_1p_1q_2 + q_2p_3q_3$ and $\dstr(p_2) = -p_3q_3p_2 -
p_2q_2p_2 - p_2q_1p_1$. The full string differential is given below.
\end{example}

\begin{theorem}
$(\alg,\d = \dsft+\dstr)$ is an LSFT algebra. \label{thm:d2}
\end{theorem}

The fact that
$\d$ preserves the filtration on $\alg$ follows from the facts that
$\dsft$ and $\dstr$ also preserve the filtration; this property for
$\dsft$ has already been established, while for $\dstr$ this is clear
by construction.

Theorem~\ref{thm:d2} is the LSFT analogue of the $\dd^2=0$ result
in Legendrian contact homology, and indeed implies it.
It will be
proven in Section~\ref{sec:string}; see Proposition~\ref{prop:d2}.

\begin{example}
For $\Lambda_0$, the full derivation $\d$ is
given by
\begin{align*}
\d(p_1) &= -p_2 + (-p_3 q_3 p_1 - p_2 q_2 p_1) \\
\d(q_1) &= 1 + t + p_3 q_2 + q_2 p_3 t + (-q_1 p_1 q_1 - q_1 p_2 q_2 -
q_1 p_3 q_3) \\
\d(p_2) &= -p_1p_3 - p_3tp_1 + (-p_3q_3p_2 - p_2q_2p_2 - p_2q_1p_1) \\
\d(q_2) &= -q_1 + q_3 + (-q_1p_1q_2 + q_2p_3q_3) \\
\d(p_3) &= p_2 + (p_3q_2p_2 + p_3q_1p_1) \\
\d(q_3) &= 1+t+q_2p_1 + tp_1q_2 + (-q_1p_1q_3 - q_2p_2q_3 - q_3p_3q_3)
\\
\d(t) &= (-q_1p_1t-q_2p_2t-q_3p_3t + tp_1q_1+tp_2q_2+tp_3q_3),
\end{align*}
where the $\dstr$ contributions are enclosed in parentheses.
The curvature for this differential is $F_{\d} =
-p_1-p_3-p_3q_2p_1$, and indeed it is straightforward to check that
$\d^2 s = [-p_1-p_3-p_3q_2p_1,s]$ for all generators $s$ of the
algebra, whence $\d^2 x = [-p_1-p_3-p_3q_2p_1,x]$ for all $x\in\alg$.
\end{example}

\begin{figure}[b]
\centerline{
\includegraphics[height=0.5in]{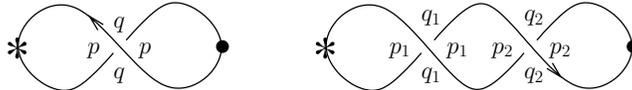}
}
\caption{
Two Legendrian unknots, the standard one with $tb=-1$ and $r=0$ (left)
and a once-stabilized one with $tb=-2$ and $r=1$ (right).
}
\label{fig:unknots}
\end{figure}

\begin{example}
For reference and comparison, we give here the derivations for the
standard Legendrian unknot and a once-stabilized Legendrian unknot
with $r=1$, as shown in Figure~\ref{fig:unknots}. The former has
\begin{align*}
\d(p) &= 0 \\
\d(q) &= 1 + t - qpq \\
\d(t) &= -qpt + tpq
\end{align*}
and $|p|=-2$, $|q|=1$, $|t|=0$, $F_{\d} = -p$;
the latter has
\begin{align*}
\d(p_1) &= q_2p_2p_1 - p_1p_2q_2 \\
\d(q_1) &= t - p_2 + p_2q_2q_1 - q_1p_1q_1 + q_1q_2p_2 \\
\d(p_2) &= 0 \\
\d(q_2) &= 1 - p_1 + q_2p_2q_2 \\
\d(t) &= p_2q_2t - q_1p_1t + tp_1q_1 - tq_2p_2
\end{align*}
and $|p_1| = 0$, $|q_1| = -1$, $|p_2| = -2$, $|q_2| = 1$, $|t|=-2$,
$F_{\d} = p_2$.
\end{example}

We next state the invariance result for LSFT algebras. Our invariance
proof requires us to restrict to a special class of Legendrian
isotopies, though we will see that this restriction covers all
Legendrian isotopies if we instead restrict to particular types of
$xy$ projections.

\begin{figure}
\centerline{
\includegraphics[height=3in]{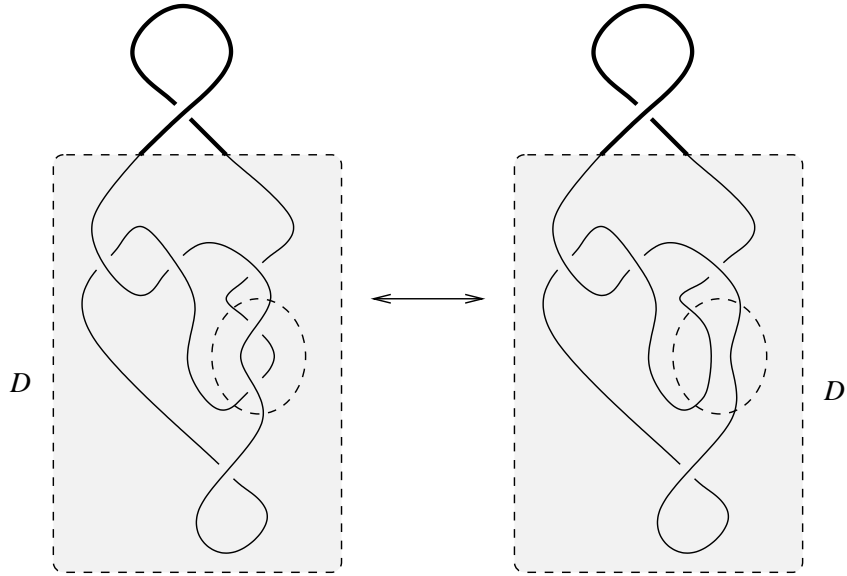}
}
\caption{
Restricted Reidemeister II move.
}
\label{fig:ReidIIrestricted}
\end{figure}

\begin{definition}
Two $xy$ projections $\Lambda_1,\Lambda_2$ of Legendrian knots are
related by a \textit{restricted Reidemeister II move} if there is an
embedded disk
$D \subset \R^3$ such that $\Lambda_1,\Lambda_2$ are identical
outside $D$, each with exactly one crossing outside $D$,
$\Lambda_1\cap \partial D = \Lambda_2\cap\partial D$ consists of two
points, and $\Lambda_1 \cap D, \Lambda_2 \cap D$ are related by a
Reidemeister II move inside $D$. See Figure~\ref{fig:ReidIIrestricted}.

\label{def:restricted}

Two $xy$ projections are related by \textit{restricted Reidemeister
  moves} if they are related by a sequence of Reidemeister III moves
and restricted Reidemeister II moves; a \textit{restricted Legendrian isotopy} is a Legendrian isotopy given in the $xy$ projection by restricted Reidemeister moves.
\end{definition}

Note that the Legendrian knot in Figure~\ref{fig:unknot2} is
related to the standard one-crossing unknot by a Reidemeister II move
eliminating crossings $2$ and $3$, but this move is not a restricted
Reidemeister II move. It is clear that the knot from
Figure~\ref{fig:unknot2} is not related to the standard unknot by
restricted Reidemeister moves, though one can show that its LSFT
algebra is equivalent to that of the standard unknot.

Recall that there is a standard procedure, called ``morsification''
\cite{Fer} or ``resolution'' \cite{NgCLI}, to obtain an $xy$ projection
from a front ($xz$) projection of a Legendrian knot, by smoothing out
left cusps and replacing right cusps by loops.

\begin{proposition}
The resolutions of the fronts of two Legendrian isotopic knots can be
related by restricted Reidemeister moves.
\label{prop:front}
\end{proposition}

\begin{proof}
Examine the resolutions of Legendrian Reidemeister moves for
fronts: front Reidemeister III resolves to a usual Reidemeister III
move; front Reidemeister I and II both resolve to Reidemeister II
moves that are restricted since they do not involve the
rightmost cusp of the front.
\end{proof}

The next result
is the LSFT version of invariance, and again implies the analogous
result in contact homology.

\begin{theorem}
  If $\Lambda$ and $\Lambda'$ are related by restricted Legendrian isotopy, then the LSFT algebras for $\Lambda$ and $\Lambda'$ are equivalent.
\label{thm:invariance}
\end{theorem}

\noindent
Theorem~\ref{thm:invariance} will be proved in Section~\ref{sec:invariance}.

\begin{corollary}
The LSFT algebra associated to the resolution of a Legendrian front is
an invariant of the corresponding Legendrian knot.
\label{cor:invariance}
\end{corollary}

As mentioned in the Introduction, it is not unreasonable to guess that one can
extend Theorem~\ref{thm:invariance} to cover all Legendrian isotopies
and not just restricted ones, but one might need to broaden the notion
of equivalence to allow arbitrary curved $A_\infty$ morphisms.


\begin{corollary}
The stable tame isomorphism type of the contact homology DGA $(\F^0\!\alg/\F^1\!\alg,\dd)$ is invariant under restricted Legendrian isotopy.
\label{cor:LCH}
\end{corollary}

%
%
%

\noindent
In fact, an examination of the proof of Theorem~\ref{thm:invariance}
shows that the contact homology DGA is
invariant under all Legendrian isotopies, not just restricted ones;
this recovers the original result of \cite{Ch02}.


\begin{remark}[\textit{The LSFT algebra and Legendrian contact
  homology}]
\label{rem:LCH}
In Remark~\ref{rem:assocgraded}, we identified $\F^0\!\alg/\F^1\!\alg$
with the Chekanov--Eliashberg differential graded algebra \cite{Ch02,Eli98}
calculating Legendrian contact homology. This holds not
only in Chekanov's original formulation over $\Z/2$, but also
over $\Z$ in the formulation of \cite{EES05a,EES07,ENS}.
There are, however, two caveats to this identification. First, the
signs used here do not coincide precisely with the signs from
\cite{ENS}, though they do agree with another sign assignment for
Legendrian contact homology given in \cite{EES05b}. However, up to a
basis change, all possible sign assignments are equivalent. The
precise statement is given and proven in Appendix~A.

Second, there is a base ring issue. In the standard formulation of
Legendrian contact homology, the differential graded algebra is
generated by Reeb chords (the $q_j$'s) over the group ring
$\Z[H_1(\Lambda)]$, which for knots is $\Z[t,t^{-1}]$. In
particular, $t^{\pm 1}$ commutes with all of the $q_j$'s. By contrast,
$\F^0\!\alg/\F^1\!\alg$ is generated by the $q_j$'s and also $t^{\pm 1}$,
with $\d(t^{\pm 1}) = 0$, and $t^{\pm 1}$ does not commute with the
$q_j$'s. We can think of the contact homology differential graded
algebra as a quotient of
$\F^0\!\alg/\F^1\!\alg$ by commutators involving $t^{\pm 1}$.

On the other hand, there is no obvious reason why, in formulating
Legendrian contact homology, we should impose the
relation that $t$ commutes with the $q_j$'s. One could reasonably
define Legendrian contact homology (even in situations more general
than knots in $\R^3$) without this relation. In our case, we would
precisely recover $\F^0\!\alg/\F^1\!\alg$.
\end{remark}


\subsection{The cyclic and commutative complexes}
\label{ssec:cyccomm}

We now discuss two quotient complexes
that can be derived from the LSFT algebra or any curved
dg-algebra. The cyclic complex has close relations to string topology
and the geometric motivation for the LSFT algebra; see
Section~\ref{sec:string}.\footnote{Cyclic constructions are
common in Symplectic Field Theory and related topics. See for instance \cite{BEE}.}
The commutative complex may be useful from a computational standpoint,
especially since it has a particularly simple formulation in the case
of the LSFT algebra, as we discuss at the end of this section.

\begin{definition}
Let $(\alg,d,F)$ be a curved dg-algebra.
\begin{enumerate}
\item
Let $\I$ be the submodule of $\alg$ generated (over $\Z$, not
over $\alg$) by commutators\footnote{The restriction that one of $x,y \in
\F^1\!\alg$ is unnecessary for most purposes, but is needed for
the theory to include full Legendrian contact homology, rather than a
cyclic version, as a quotient.}
\[
\{[x,y] \, | \, x,y\in\alg \text{ and at least one of $x,y \in
  \F^1\!\alg$}\}.
\]
The \textit{cyclic complex} associated to
$(\alg,d,F)$ is
$(\alg\tocyc = \alg/\I, \d)$, where $\d$ is the induced differential on $\alg\tocyc$.
\item
Let $\J$ be the subalgebra of $\alg$ generated (over
$\alg$) by commutators $[x,y]$ for all $x,y\in\alg$. The
\textit{commutative complex} associated to $(\alg,d,F)$ is $(\alg\tocomm =
\alg/\J,\d)$, where $\d$ is the induced differential on $\alg\tocomm$.
\end{enumerate}
\end{definition}

\noindent
The key point here is that $\d^2=0$ on $\alg\tocyc$ and $\alg\tocomm$, by the
definition of curved dg-algebra.

When $\alg$ is a tensor
algebra (as for the LSFT algebra), $\alg\tocyc$ is generated by ``cyclic
words'', or words modulo cyclic permutations of the letters (for words
in $\F^1\!\alg$). Note
that $\alg\tocyc$ is a $\Z$-module and not an algebra; it however still
inherits the grading and filtration from $\alg$.
By contrast, $\alg\tocomm$ is a (sign-)commutative algebra over $\Z$, the
polynomial algebra generated by the generators of $\alg$.
There are obvious quotient maps
\[
\xymatrix{
(\alg,\d,F) \ar@{->>}[r] & (\alg\tocyc,\d) \ar@{->>}[r] & (\alg\tocomm,\d),
}
\]
and the latter induces a map on homology.

We next show that the cyclic and commutative complexes associated to
the LSFT algebra of a Legendrian knot are invariant. This is a direct
consequence of the following result.

\begin{proposition}
If $(\alg,\d)$ and $(\alg',\d')$ are equivalent LSFT algebras, then their
rational cyclic and commutative quotient complexes are filtered chain
homotopy equivalent. In particular, they are quasi-isomorphic:
\[
H_*(\cyc\otimes\Q,\dd) \cong
H_*((\alg')\tocyc\otimes\Q,\dd')
\]
as filtered
graded $\Q$-modules, and
\[
H_*(\comm\otimes\Q,\dd) \cong
H_*((\alg')^{\operatorname{comm}}\otimes\Q,\dd')
\]
as filtered graded $\Q$-algebras.
\label{prop:quasi}
\end{proposition}

The substance of the proof of Proposition~\ref{prop:quasi}, which we
give below, is invariance under stabilization. Recall from the proof of
Proposition~\ref{prop:htpyequiv} the maps
$\iota,\pi$ between an LSFT algebra $(\alg,\d)$ and its stabilization
$(S_i\alg,\d)$. Though the homotopy operator $H$ from that proof does not
descend from $S_i\alg$ to $(S_i\alg)^{\text{cyc}}$, we can define a
slight variant that
serves as the corresponding homotopy operator for
$(S_i\alg)\tocyc$. Let $h\co S_i\alg \to S_i\alg$ be the derivation
defined by $h(q_b)=q_a$, $h(p_a)=p_b$, $h(q_a)=h(p_b)=0$,
$h(q_j)=h(p_j)=h(t^{\pm 1})=0$ for all $j$ (besides $a,b$).
For a word $w \in \alg$,
let $\sigma(w)$ be the total number of occurrences of $q_a,q_b,p_a,p_b$ in
$x$. Now define $H\tocyc\co (S_i\alg) \otimes\Q \to (S_i\alg) \otimes\Q$ by
\[
H\tocyc(w) = \begin{cases}
0 & \text{if $\sigma(w) = 0$} \\
\frac{1}{\sigma(w)} h(w) & \text{if $\sigma(w) > 0$}.
\end{cases}
\]

\begin{lemma}
On $(S_i\alg)\otimes\Q$, we have
\label{lem:cychtpy}
$\Id_{S_i\alg} - \iota \circ \pi = H\tocyc \circ
\d + \d \circ H\tocyc$.
\end{lemma}

\begin{proof}
  It suffices to show that $(h\circ\d + d\circ h)(w) = \sigma(w)
  w$ for all words $w$, since $\d$ preserves the number of occurrences
  of $q_a,q_b,p_a,p_b$. But both $h\circ\d + d\circ h$ and the map
  (generated on words by) $w \mapsto \sigma(w) w$ are derivations, and
  they agree on the generators $q_a,q_b,p_a,p_b,q_j,p_j,t$.
\end{proof}

\begin{proof}[Proof of Proposition~\ref{prop:quasi}]
  Let $(\alg,\d)$ and $(\alg',\d')$ be equivalent LSFT algebras. We
  show that $(\cyc\otimes\Q,\dd)$ and
  $((\alg')\tocyc\otimes\Q,\dd')$ are filtered chain homotopy
  equivalent; the proof of the corresponding result for the
  commutative complexes is nearly identical.
  The result clearly holds if $(\alg,\d)$ and $(\alg',\d')$
  are related by a basis change or gauge change. Thus we may assume
  that $(\alg',\d')=(S_i\alg,\d)$ is a stabilization of $(\alg,\d)$.
  In this case, the inclusion and projection maps $\iota,\pi$ between $\alg$ and
  $S_i\alg$ induce chain maps $\iota\co \cyc \to (S_i\alg)\tocyc$ and
  $\pi\co (S_i\alg)\tocyc \to \cyc$. Furthermore, $\pi\circ\iota =
  \Id_{\cyc}$, while $\iota\circ\pi$ is chain homotopic over
  $\Q$ to $\Id_{(S_i\alg)\tocyc}$ by Lemma~\ref{lem:cychtpy}.
\end{proof}

\begin{corollary}
The cyclic and commutative complexes associated to the LSFT algebra of
a Legendrian knot are invariant, up to filtered chain homotopy equivalence, under restricted Legendrian isotopy.
\end{corollary}

We remark that it can be checked that the powers of
$F$, $F^n \in \alg$ for $n\geq 0$, descend to invariant classes in $H(\alg\tocyc\otimes\Q,d)$. See
\cite{Positselski} for a fuller discussion, where these invariant
classes are called ``Chern classes''.


To end this section, we observe that the commutative complex for an
LSFT algebra has a rather simpler formulation than the full LSFT algebra.
More precisely, on $\comm$ we can still define
$\d = \dsft + \dstr$, with $\dsft$ defined as for $\alg$, but now
$\dstr$ can be given as follows.


\begin{figure}[b]
\centerline{
\includegraphics[height=0.5in]{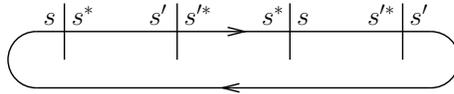}
}
\caption{
A schematic diagram traversing $\pi_{xy}(\Lambda)$ along its
orientation. In this picture, the crossings corresponding to $s$ and
$s'$ are interlaced, and we have
$s \dasharrow s' \dasharrow s^* \dasharrow s'^* \dasharrow s$.
}
\label{fig:interlaced}
\end{figure}

\begin{definition}
Two crossings in $\pi_{xy}(\Lambda)$ are \textit{interlaced}
if, in traversing the knot one full time, we encounter one crossing,
then the other, then the first again, then the second again. If
crossings corresponding to $s$ and $s'$ are interlaced, we say that
$s$ is \textit{interlaced toward} $s'$ and write $s \dasharrow s'$
if, when we traverse $\Lambda$
along its orientation starting from $s^-$, we encounter $(s')^-$
before $(s')^+$.
\end{definition}

\noindent See Figure~\ref{fig:interlaced} for an illustration. If the
crossings corresponding to $s$ and $s'$ are interlaced, then there are
two possibilities, $s \dasharrow s' \dasharrow s^* \dasharrow s'^*
\dasharrow s$ or $s \dasharrow s'^* \dasharrow s^* \dasharrow s'
\dasharrow s$, depending on the orientation of $\Lambda$.

\begin{proposition}
In $\alg\tocomm$, we have $\dstr(t^{\pm 1}) = 0$, while if $s$ is a
$p$ or $q$, then
\label{prop:commdiff}
\[
\dstr(s) = \sum_{s \dasharrow s'} \sftb{s'^*,s'} \, s'^*s's
\]
where the sum is over all $s'$ such that $s$ is
interlaced toward $s'$, and $\sftb{s'^*,s'}$ is $1$ if $s'$
is a $q$, $-1$ if $s'$ is a $p$.
\end{proposition}

The proof of Proposition~\ref{prop:commdiff} is an exercise in chasing
signs, and we leave it to the interested reader. Note that $\dstr$
does not depend on the choice of base points $\bullet,\ast$ on
$\Lambda$. Thus the differential $\d$ on $\alg\tocomm$ is independent
of $\bullet$, and only depends on $\ast$ insofar as $\ast$ keeps track
of powers of $t$ in the SFT differential, cf.\ group-ring coefficients
in Legendrian contact homology \cite{ENS}.

\secspace


\section{String Interpretation of the LSFT Algebra}
\label{sec:string}

It will be useful to have another description of the LSFT algebra,
closer to the standard SFT formalism and string topology. This allows
us to prove the ``$\d^2=0$ result'', Theorem~\ref{thm:d2}.

\subsection{Broken closed strings and the SFT bracket}
\label{ssec:sftb}

Generators of the LSFT algebra of a Legendrian knot $\Lambda$
are more conveniently seen as strings
on $\Lambda$. From the holomorphic perspective, these are the
boundaries of holomorphic disks with boundary on $\Lambda\times\R$. It
will be fruitful, however, to consider all possible strings, not just
those that arise as the boundary of a disk.

\begin{definition}
Let $\Lambda\in\R^3$ be a Legendrian knot with Reeb chords
$R_1,\dots,R_n$, and let $R_i^\pm\in\Lambda$ denote the endpoints of Reeb chord
$R_i$. For fixed $k \geq 0$, choose $k$ distinct points
$\tau_1,\dots,\tau_k$ on an oriented
circle $S^1$ so that they appear sequentially in order; we refer to
these points as \textit{punctures} of $S^1$, and the punctures
divide $S^1$ into $k$ intervals, which we denote
$[\tau_1,\tau_2],[\tau_2,\tau_3],\dots,[\tau_k,\tau_1]$.
A \textit{broken closed string} of length $k$ is a piecewise continuous
map $\gamma \co S^1 \to \Lambda$ such that:
\begin{enumerate}
\item
$\gamma|_{[\tau_j,\tau_{j+1}]}$ is continuous for each $j$;
\item
for each $j=1,\dots,k$, either $\lim_{\tau\to \tau_j^{\pm}}
\gamma(\tau) = R_{i_j}^{\pm}$ or
$\lim_{\tau\to \tau_j^{\pm}} \gamma(\tau) = R_{i_j}^{\mp}$ for some $i_j$.
\end{enumerate}
We consider broken closed strings up to orientation-preserving
reparametrization of the domain $S^1$.

If we are given a point $\bullet\in\Lambda$ distinct from any of the
$R_i^\pm$, then choose a point $\tau_0\in(\tau_k,\tau_1)\subset S^1$; a
\textit{based broken closed string} of length $k$ is a broken closed
string $\gamma$ of length $k$ such that $\gamma(\tau_0)=\bullet$.
\end{definition}

Given distinct points $\ast,\bullet\in\Lambda$, we obtain a map $w$
between based broken closed strings and words in $\alg$. If $\gamma$
is a based broken closed string of length $k$, then define the word
associated to $\gamma$ to be
\[
w(\gamma) = t^{a_0} s_1 t^{a_1} s_2 \cdots s_k t^{a_k}
\]
where $t^{a_j}$ is the number of times $\gamma|_{[\tau_j,\tau_{j+1}]}$
passes through $\ast$, counted with sign according to the orientation
of $\Lambda$, and
\[
s_j = \begin{cases}
p_{i_j} & \text{if $\lim_{\tau\to \tau_j^{\pm}} \gamma(\tau) = R_{i_j}^{\pm}$}
\\
q_{i_j} & \text{if $\lim_{\tau\to \tau_j^{\pm}} \gamma(\tau) = R_{i_j}^{\mp}$}.
\end{cases}
\]
Note that the correspondence between based broken closed strings and
words in $\alg$ is bijective if we mod out the strings by homotopy.

Similarly, we can define a map, which we also denote by $w$, between
broken closed strings and words in $\cyc$. Note that this map does not
depend on the choice of $\bullet$, as changing $\bullet$ corresponds
to conjugation by some power of $t$.

\begin{figure}
\centerline{
\includegraphics[height=2in]{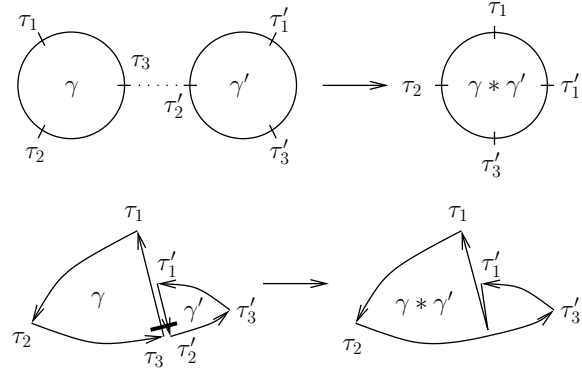}
}
\caption{
Gluing broken closed strings $\gamma,\gamma'$ at punctures $\tau_3,\tau_2'$
to produce a new broken closed string $\gamma * \gamma'$. The top
illustrates the picture in the domain; the bottom illustrates the
image of $\gamma,\gamma',\gamma*\gamma'$ in the $xy$ projection
$\pi_{xy}(\Lambda)$. The heavy bar indicates the glued corners.
}
\label{fig:SFTglue}
\end{figure}

If $\gamma,\gamma'$ are broken closed strings of length $k,k'$
respectively, and a puncture from each is mapped to (the endpoints of)
the same Reeb chord but in opposite directions, then we can glue
$\gamma,\gamma'$ at this puncture to obtain another broken closed
string of length $k+k'-2$. More precisely, suppose that the $S^1$
domain of $\gamma$ has punctures $\tau_1,\dots,\tau_k$, the $S^1$ domain
of $\gamma'$ has punctures $\tau_1',\dots,\tau_{k'}'$, and there are $j$ and
$j'$ such that
\[
\lim_{\tau\to \tau_j^{\pm}} \gamma(\tau) = \lim_{\tau\to
  \tau_{j'}^{\mp}} \gamma'(\tau);
\]
then we can define a broken closed string $\gamma * \gamma'$ on $S^1$
with sequential punctures
$\tau_{j+1},\dots,\tau_k,\tau_1,\dots,\tau_{j-1},\tau_{j'+1}',\dots,
\tau_{k'}',\tau_1',\dots,\tau_{j'-1}'$
by
\begin{align*}
(\gamma * \gamma')|_{[\tau_i,\tau_{i+1}]} &= \gamma|_{[\tau_i,\tau_{i+1}]}\\
\intertext{for all $i=1,\dots,k$ with $i \neq j$,}
(\gamma * \gamma')|_{[\tau_i',\tau_{i+1}']} &= \gamma'|_{[\tau_i',\tau_{i+1}']} \\
\intertext{for all $i=1,\dots,k'$ with $i\neq j'$,}
(\gamma * \gamma')|_{[\tau_{j-1},\tau_{j'+1}']} &= \gamma|_{[\tau_{j-1},\tau_j]} \cup
\gamma'|_{[\tau_{j'}',\tau_{j'+1}']}, \\
(\gamma * \gamma')|_{[\tau_{j'-1}',\tau_{j+1}]} &=
\gamma'|_{[\tau_{j'-1}',\tau_{j'}']} \cup \gamma|_{[\tau_{j},\tau_{j+1}]}.
\end{align*}
See Figure~\ref{fig:SFTglue} for an illustration.

\begin{figure}
\centerline{
\includegraphics[height=2in]{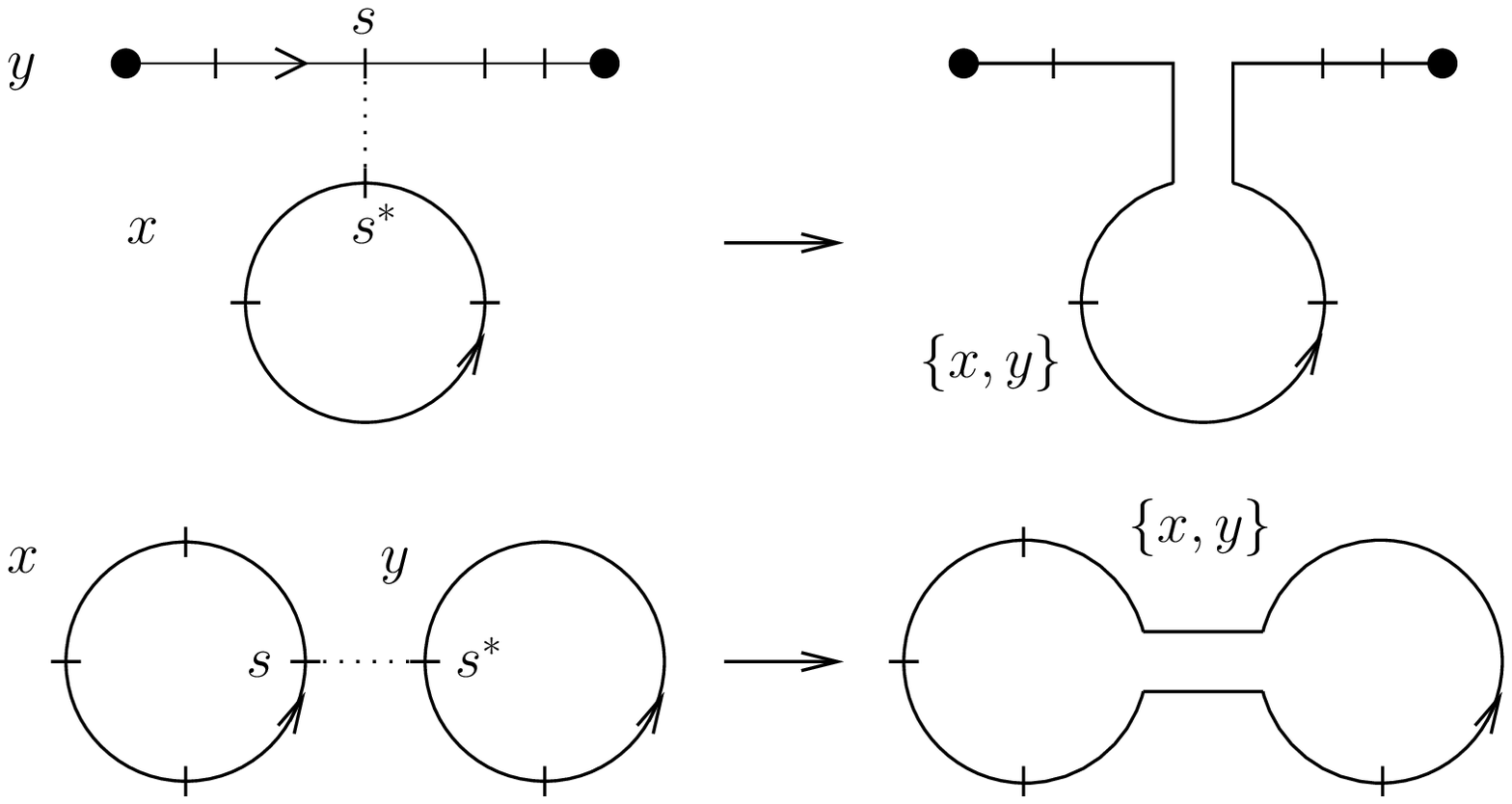}
}
\caption{
Gluing $x\in\cyc$ to $y\in\alg$ to get $\sftb{x,y}\in\alg$ (top);
gluing $x\in\cyc$ to $y\in\cyc$ to get $\sftb{x,y}\in\cyc$
(bottom). The notches represent corners (punctures), and the dots
represent the base point $\bullet$.
}
\label{fig:bracket}
\end{figure}

Using the gluing operation, we can define the SFT bracket of two
broken closed strings to be the sum of all possible gluings of the
broken closed strings. This gives an operation
$\sftb{\cdot,\cdot}\co\cyc\otimes\cyc\to\cyc$ (at least mod $2$). In
the same way, we can define the SFT bracket of a broken closed string
and a based broken closed string to be the sum of all based broken
closed strings obtained by gluing; this gives a mod $2$ map
$\sftb{\cdot,\cdot}\co\cyc\otimes\alg\to\alg$. We refer to either
operation as the \textit{SFT bracket}. See Figure~\ref{fig:bracket}.

We can define the SFT bracket in a more precise algebraic manner, with
the added benefit of lifting to $\Z$, as follows.
First, given a word $w \in \A$ ending in a $p$ or $q$, define a contraction map
$\iota_w\co\A\to\A$ as follows. Write $w = w'
p_i$ or $w = w' q_i$ for some word $w'$, and set
\[
\iota_{w'p_i}(s) = \begin{cases}
\delta_{ij} w'
& \text{if $x =
q_j$} \\ 0 & \text{if $x = p_j$ or $x = t^{\pm 1}$}
\end{cases}
\]
and
\[ \iota_{w'q_i}(x) = \begin{cases}
- \delta_{ij} w' & \text{if $x = p_j$} \\
0 & \text{if $x = q_j$ or $x = t^{\pm 1}$},
\end{cases}
\]
where $\delta_{ij}$ is the Kronecker delta function;
extend $\iota_w$ to a map $\A\to\A$ by linearity and the following
modified Leibniz rule:
\[
\iota_w(xy) = (\iota_w x) y + (-1)^{(|w|+1)|x|} x (\iota_w y).
\]
The unusual sign ensures that $\iota_w$ descends to a map
$\A\tocyc\to\A\tocyc$.



Before extending contraction to cyclic words and defining the SFT
bracket, we
introduce notation for the set of all words that project to a
particular cyclic word.

\begin{definition}
  Let $w$ be a word in $\A$. The \textit{length} $l(w)$ of
  $w$ is the number of $q$'s and $p$'s in $w$. The
  \textit{cyclic word set} $\CC(w)$ of $w$ is the $l(w)$-element
  multiset in $\alg$ of words equal to $w$ in $\cyc$. More precisely,
  if $w = t^{a_0} s_1 t^{a_1} s_2 \cdots s_{l(w)} t^{a_{l(w)}}$, where
  each $s_i$ is a $q$ or $p$, then
\begin{align*}
\CC(w) &= \left\{ t^{a_0+a_{l(w)}} s_1 t^{a_1} s_2 \cdots s_{l(w)}, \right.\\
&
(-1)^{|s_1|(|w|-|s_1|)} t^{a_1} s_2 \cdots s_{l(w)} t^{a_0+a_{l(w)}} s_1,\dots,\\
& \left. (-1)^{|s_{l(w)}|(|w|-|s_{l(w)}|)}
t^{a_{l(w)-1}} s_{l(w)} t^{a_0+a_{l(w)}} s_1 \cdots s_{l(w)-1}\right\}.
\end{align*}
Note that if $w,w'$ represent the same element in $\cyc$, then
$\CC(w)=\CC(w')$.
\end{definition}

Now if $w$ is a word in $\A$ and $[w]$ is the image of $w$
in $\A\tocyc$, then we define a contraction map
$\iota_{[w]}\co\alg\to\alg$ by
\[
\iota_{[w]}(x) = \sum_{y\in\CC(w)} \iota_y(x).
\]
Here we use the convention that $\iota_{-w}(x) = -\iota_w(x)$ if $w$
is a word in $\A$.
The contraction map extends by linearity to a map
$\iota_{\cdot}(\cdot)\co\cyc\otimes\alg\to\alg$.

\begin{definition}
The \textit{SFT bracket} $\sftb{\cdot,\cdot}\co\cyc\otimes\alg\to\alg$
is defined by $\sftb{x,y} = \iota_x(y)$.
This descends to a map $\cyc\otimes\cyc\to\cyc$, which
we also denote by $\sftb{\cdot,\cdot}$.
\label{def:sftb}
\end{definition}

For reference, the full sign rule, which can be deduced from the
definition of $\iota$, is as follows: the $s$ and $s^*$ entries in two
words $w_1sw_2$, $w_3s^*w_4$ pair together to give
\begin{equation}
\sftb{w_1sw_2,w_3s^*w_4} = (-1)^{|w_2||w_1s|+(|w_1sw_2|+1)|w_3|}
\sftb{s,s^*} w_3w_2w_1w_4 + \cdots.
\label{eq:sftb}
\end{equation}

\begin{proposition}[Properties of the SFT bracket]
\begin{enumerate}
\item
Let $x,y\in\cyc$.  \label{it:sftb1} Then
\[
\sftb{y,x} = (-1)^{|x||y|+|x|+|y|} \sftb{x,y}.
\]
\item
Let $x\in\cyc$ and $y,z\in\alg$. \label{it:sftb2} Then
\[
\sftb{x,yz} = \sftb{x,y} z + (-1)^{(|x|+1)|y|} y \sftb{x,z}.
\]
\item
Let $x,y\in\cyc$ and $z\in\alg$.  \label{it:sftb3} Then we have the
following version of the Jacobi identity:
\[
\sftb{x,\sftb{y,z}} + (-1)^{|x||y|+|x|+|y|} \sftb{y,\sftb{x,z}} =
\sftb{\sftb{x,y},z}.
\]
\end{enumerate}
\label{prop:sftb}
\end{proposition}

\begin{proof}
We first establish the proposition mod $2$. Note that (\ref{it:sftb1}) is
clear, while (\ref{it:sftb2}) can be pictorially represented:
\[
\sftb{x,yz} =
\raisebox{-11pt}{\includegraphics[height=33pt]{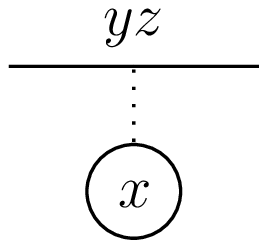}} =
\raisebox{-11pt}{\includegraphics[height=33pt]{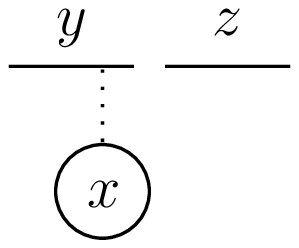}} +
\raisebox{-11pt}{\includegraphics[height=33pt]{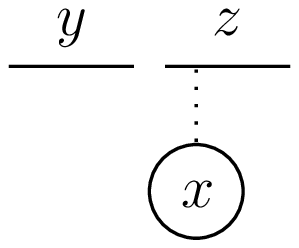}} =
\sftb{x,y} z + y
\sftb{x,z}.
\]
For (\ref{it:sftb3}), we have
\begin{align*}
\sftb{x,\sftb{y,z}} &=
\raisebox{-11pt}{\includegraphics[height=33pt]{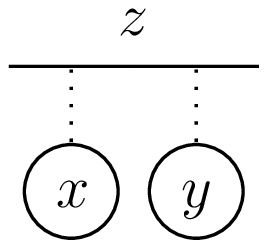}} +
\raisebox{-11pt}{\includegraphics[height=33pt]{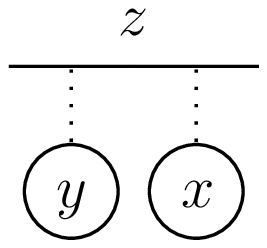}} +
\raisebox{-11pt}{\includegraphics[height=33pt]{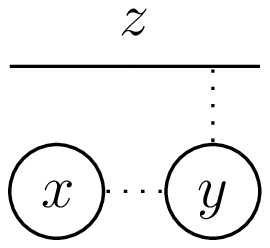}} \\
\sftb{y,\sftb{x,z}} &=
\raisebox{-11pt}{\includegraphics[height=33pt]{figures/sftb4.eps}} +
\raisebox{-11pt}{\includegraphics[height=33pt]{figures/sftb5.eps}} +
\raisebox{-11pt}{\includegraphics[height=33pt]{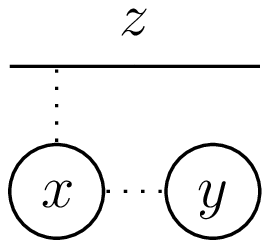}} \\
\sftb{\sftb{x,y},z} &=
\raisebox{-11pt}{\includegraphics[height=33pt]{figures/sftb7.eps}} +
\raisebox{-11pt}{\includegraphics[height=33pt]{figures/sftb6.eps}}.
\end{align*}
Checking the signs is now a routine exercise using equation (\ref{eq:sftb}).
\end{proof}

%

\subsection{The $\delta$ map}
\label{ssec:delta}

Having defined the SFT bracket, we now define another operation on
strings, the $\delta$ map. This is essentially a string cobracket
operation in the language of string topology. First we need to take a
slightly closer look at broken closed strings.

\begin{definition}
A \textit{generic broken closed string} is a broken closed string
$\gamma\co (S^1;\tau_1,\dots,\tau_k) \to \Lambda$ such that whenever
$\gamma'(\tau) =
0$, $\gamma(\tau) \not\in \RR$, where we recall that $\RR$ is the set
of Reeb chord endpoints; in particular, $\gamma'({\tau_i}^\pm) \neq 0$, where
$\gamma'(\tau_i^{\pm}) = \lim_{\tau\to \tau_i^{\pm}} \gamma'(\tau)$.

A generic broken closed string \textit{has holomorphic corners} if for
each $i$,
$(\gamma'(\tau_i^-),\gamma(\tau_i^+)-\gamma(\tau_i^-),\gamma'(\tau_i^+))$ is a
positively oriented frame in $\R^3$. This condition is most easily
interpreted in the $xy$ projection: the image of $\gamma$ in
$\pi_{xy}(\Lambda)$ near each $\tau_i$ makes a left turn at the corner.
\end{definition}

\begin{figure}
\centerline{
\includegraphics[height=0.5in]{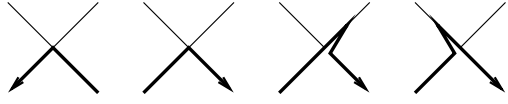}
}
\caption{
From left to right: a holomorphic corner of a broken closed string; a
non-holomorphic corner; and two broken closed strings with holomorphic
corners homotopic to the one with the non-holomorphic corner.
}
\label{fig:holcorner2}
\end{figure}

\begin{figure}
\centerline{
\includegraphics[height=1.5in]{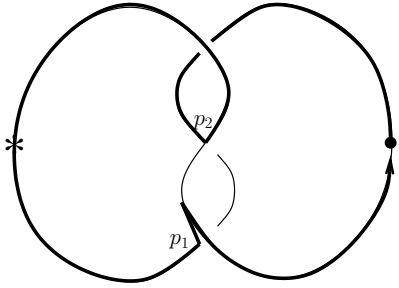}
}
\caption{
A generic broken closed string $\gamma$
with holomorphic corners. The word associated to $\gamma$ is $p_2tp_1$.
}
\label{fig:bcs}
\end{figure}

See Figure~\ref{fig:holcorner2} for an illustration of holomorphic
corners and Figure~\ref{fig:bcs} for an example of a generic broken
closed string with holomorphic corners.
It is easy to see that any broken closed string is homotopic to a
generic broken closed string with holomorphic corners.

Now suppose that $\gamma\co(S^1;\tau_1,\dots,\tau_k)\to\Lambda$ is a
generic broken
closed string of length $k$, and suppose $\tau\in (\tau_i,\tau_{i+1})$
satisfies $\gamma(\tau) = R_j^{\pm}$ for some $j$ and some choice of
$\pm$; in this case, we say
that $\tau$ is \textit{interior Reeb} for $\gamma$. We can then define a
broken closed string $\delta(\gamma;\tau)$ of length $k+2$ to have
sequential punctures
$\tau_1,\dots,\tau_i,\tau^{(1)},\tau^{(2)},\tau_{i+1},\dots,\tau_k$, and
\begin{align*}
\delta(\gamma;\tau)|_{[\tau_i,\tau^{(1)}]} &= \gamma|_{[\tau_i,t]} \\
\delta(\gamma;\tau)|_{[\tau^{(1)},\tau^{(2)}]} &= \text{constant path at
  $R_j^{\mp}$} \\
\delta(\gamma;\tau)|_{[\tau^{(2)},\tau_{i+1}]} &= \gamma|_{[\tau,\tau_{i+1}]} \\
\delta(\gamma;\tau)|_{[\tau_j,\tau_{j+1}]} &= \gamma|_{[\tau_j,\tau_{j+1}]} \quad
\text{for $j \neq i$.}
\end{align*}
Note that $\delta(\gamma;\tau)$ is not generic, but it can be perturbed
to become generic; furthermore, if we stipulate that the perturbed broken
closed string has holomorphic corners on the domain interval
$(\tau_i,\tau_{i+1})$, then the perturbation is unique up to homotopy
through generic broken closed strings with holomorphic corners.

If $\gamma$ is a generic based broken closed curve and $\tau$ is
interior Reeb for
$\gamma$,  we can similarly define a
based broken closed curve $\delta(\gamma;\tau)$.

\begin{figure}
\centerline{
\includegraphics[height=1.5in]{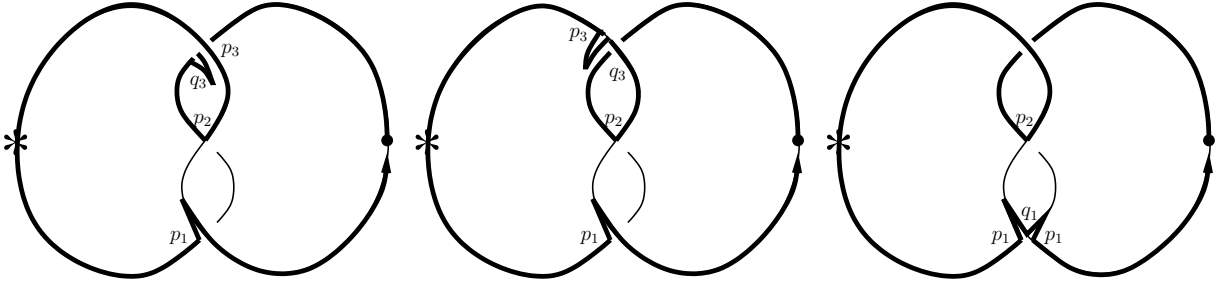}
}
\caption{
The three terms in $\delta(p_2tp_1)$ for the broken closed string
$\gamma=\gamma(p_2tp_1)$ from Figure~\ref{fig:bcs}. Each new broken
closed string has been perturbed to have holomorphic corners.
}
\label{fig:bcsdelta}
\end{figure}

We can now define a map $\delta\co\A\to\A$.
In fact, $\delta$ is just a string reformulation of $\dstr$ from
Section~\ref{ssec:comb}; see Proposition~\ref{prop:dstrdelta} below.

\begin{definition}
Let
 $w$ be a word in $\A$; we have $w = w(\gamma)$ for
some generic based broken closed string with holomorphic corners.
Define $\delta(w) \in \A$ by
\[
\delta(w) = \sum_{\text{$\tau$ interior Reeb for $\gamma$}}
\epsilon(\gamma;\tau) w(\delta(\gamma;\tau)),
\]
with $\epsilon(\gamma;\tau)$ a sign to be defined in the next paragraph.
Extend $\delta$ to $\A$ by linearity.

\label{def:deltabcs}
We define $\epsilon(\gamma;\tau) \in\{\pm
1\}$ as follows. Suppose that $w(\gamma)$, with powers of $t$
omitted, is of the form $s_1\cdots s_k$, where $s_j$ corresponds to
$\gamma(\tau_j)$, and that $\tau_i < \tau < \tau_{i+1}$. Let
$\epsilon_1 = \pm 1$ according to whether $\gamma(\tau) = R_j^{\pm
  1}$, and let $\epsilon_2 = \pm 1$ according to whether the
orientation of $\gamma$ in a neighborhood of $\tau$ agrees or
disagrees with the orientation of $\Lambda$. Finally, define
\[
\epsilon(\gamma;\tau) = (-1)^{|s_1\cdots s_i|} \epsilon_1 \epsilon_2.
\]
\end{definition}

See Figure~\ref{fig:bcsdelta} for an illustration of
Definition~\ref{def:deltabcs}.
In short, if $w = s_1 \cdots s_k$, then $\delta(w)$ is a sum of
terms of the form
\[
\epsilon_2 (-1)^{|s_1\cdots s_i|} s_1 \cdots s_i (p q) s_{i+1} \cdots s_k
\quad \text{and} \quad
- \epsilon_2 (-1)^{|s_1\cdots s_i|} s_1 \cdots s_i (q p) s_{i+1} \cdots
s_k,
\]
where $\epsilon_2$ measures the orientation of the broken closed
string for $w$ at the point where $pq$ or $qp$ is attached.

\begin{example}
Consider the broken closed string from Figure~\ref{fig:bcs},
corresponding to the word $p_2tp_1$. From Figure~\ref{fig:bcsdelta},
we see that $\delta(p_2tp_1)$ has three terms corresponding to
$(p_3q_3)p_2tp_1$, $p_2(q_3p_3)tp_1$, and $p_2tp_1(q_1p_1)$. In fact,
we have
\[
\delta(p_2tp_1) = -p_3q_3p_2tp_1 + p_2q_3p_3tp_1 - p_2tp_1q_1p_1.
\]
We verify the sign of the third term as an example. There are three
$-1$'s contributing to this sign, one from $(-1)^{|p_2tp_1|}$, one
from the fact that $q_1$ precedes $p_1$ in the parenthesis, and one
from the fact that at the point where $q_1p_1$ is added, the broken
closed string is oriented opposite to the orientation of the knot from
Figure~\ref{fig:unknot2}.

We remark that it can be readily checked that $\delta(p_2tp_1)$ agrees
with $\dstr(p_2tp_1)$, calculated from the values for $\dstr(p_1),
\dstr(p_2), \dstr(t)$ given in Section~\ref{ssec:comb}; see
Proposition~\ref{prop:dstrdelta} below.

\end{example}

We now prove some fundamental algebraic properties of $\delta$.

\begin{definition}
Define a $\Z$-linear map $\bullet\co\cyc\to\alg$ as follows. Any word
$w$ in $\A$ corresponds to a based broken closed
string $\gamma_w\co(S^1,\tau_0,\tau_1,\dots,\tau_k)\to\Lambda$, and
$\gamma_w$ can be chosen so that whenever $\gamma_w(\tau) = \bullet$,
$\gamma_w'(\tau) \neq 0$. Then $\gamma_w$ passes through $\bullet$
some number of times $n(\gamma_w)\in\Z$, counted with sign according
to the orientation of the knot. (More precisely, since $\gamma_w$
begins and ends at $\bullet$, one should ``close up'' $\gamma_w$ and
view it as a homotopy class of unbased broken closed strings when
calculating $n(\gamma_w)$.)
Now define $\bullet\co\A\to\A$ by
\[
\bullet(w) = (-1)^{|w|} n(\gamma_w) w
\]
and $\bullet\co\cyc\to\alg$ by $\bullet([w]) =
\sum_{w'\in\CC(w)} \bullet(w')$.
\end{definition}

\begin{proposition}[Properties of $\delta$]

\label{prop:delta}
\begin{enumerate}
\item
$\delta$ gives a well-defined map from $\A$ to $\A$ and induces a
well-defined map from $\A\tocyc$ to $\A\tocyc$, as well as from $\alg$
to $\alg$ and from $\cyc$ to $\cyc$.
\label{delta:1}
\item
If $x,y\in\A$, \label{delta:2} then $\delta(xy) = (\delta x) y +
(-1)^{|x|} x (\delta y)$.
\item
If $x\in\A$, \label{delta:3} then $\delta^2(x) = 0$.
\item
If $x,y\in\A\tocyc$, \label{delta:4} then
\[
\delta\sftb{x,y} = \sftb{\delta x,y} - (-1)^{|x|} \sftb{x,\delta y}.
\]
\item
If $x\in\A\tocyc$ and $y\in\A$, \label{delta:5} then
\[
\delta\sftb{x,y} = \sftb{\delta x,y} - (-1)^{|x|} \sftb{x,\delta y} +
[\bullet(x),y].
\]
\end{enumerate}
\end{proposition}

\begin{proof}
For (\ref{delta:1}), note that any two generic based broken closed
strings with holomorphic corners that
represent the same word in $\A$ can be related by a set of local
moves, depicted in Figure~\ref{fig:bcs-homotopy}. It is easy to check
that $\delta$ is unchanged by each of these moves, and (\ref{delta:1})
follows. Items (\ref{delta:2}) and (\ref{delta:3}) are both clear mod
$2$, and are readily seen to hold over $\Z$ by the definition of the signs
$\epsilon(\gamma;\tau)$.

\begin{figure}[b]
\centerline{
\includegraphics[height=1.5in]{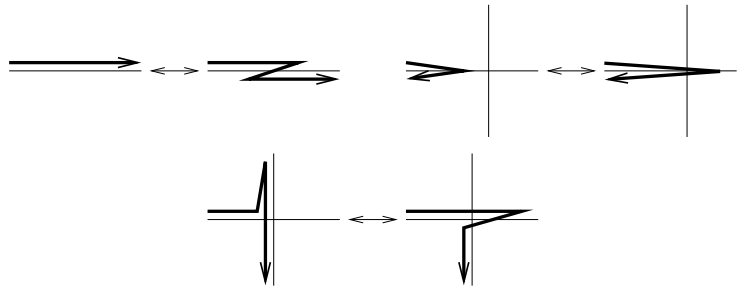}
}
\caption{
Local moves in the $xy$ projection
relating two broken closed strings that give the same word
in $\A$. The broken closed strings are drawn thickly and the
underlying $xy$ projection of $\Lambda$ is drawn thinly.
}
\label{fig:bcs-homotopy}
\end{figure}

It remains to prove (\ref{delta:4}) and (\ref{delta:5}). Assume $x$ is
a cyclic word ($x\in\A\tocyc$) and $y$ is either a cyclic word or a
word ($y\in\A\tocyc$ or $y\in\A$), and define
\[
f(x,y) = \delta\sftb{x,y} -
\sftb{\delta x,y} + (-1)^{|x|} \sftb{x,\delta y}.
\]
Most terms in
$\delta\sftb{x,y}$ have an obvious corresponding term (with the same
sign) in one of $\sftb{\delta x,y}$ or $(-1)^{|x|+1}\sftb{x,\delta y}$, and
conversely. The
exceptions are terms where the $\delta$ operation interacts with the
bracket:
\begin{itemize}
\item
every term in $\sftb{x,y}$ arises from gluing a corner $s$ in $x$ to
$s^*$ in $y$; the resulting broken closed string has a segment that
passes over the $s,s^*$ crossing, where $\delta$ can be applied;
\item
every term in $\delta x$ includes two consecutive corners not
appearing in $x$, and either can be glued to $y$;
\item
every term in $\delta y$ includes two consecutive corners not
appearing in $y$, and either can be glued to $x$.
\end{itemize}

\begin{figure}
\centerline{
\includegraphics[height=5in]{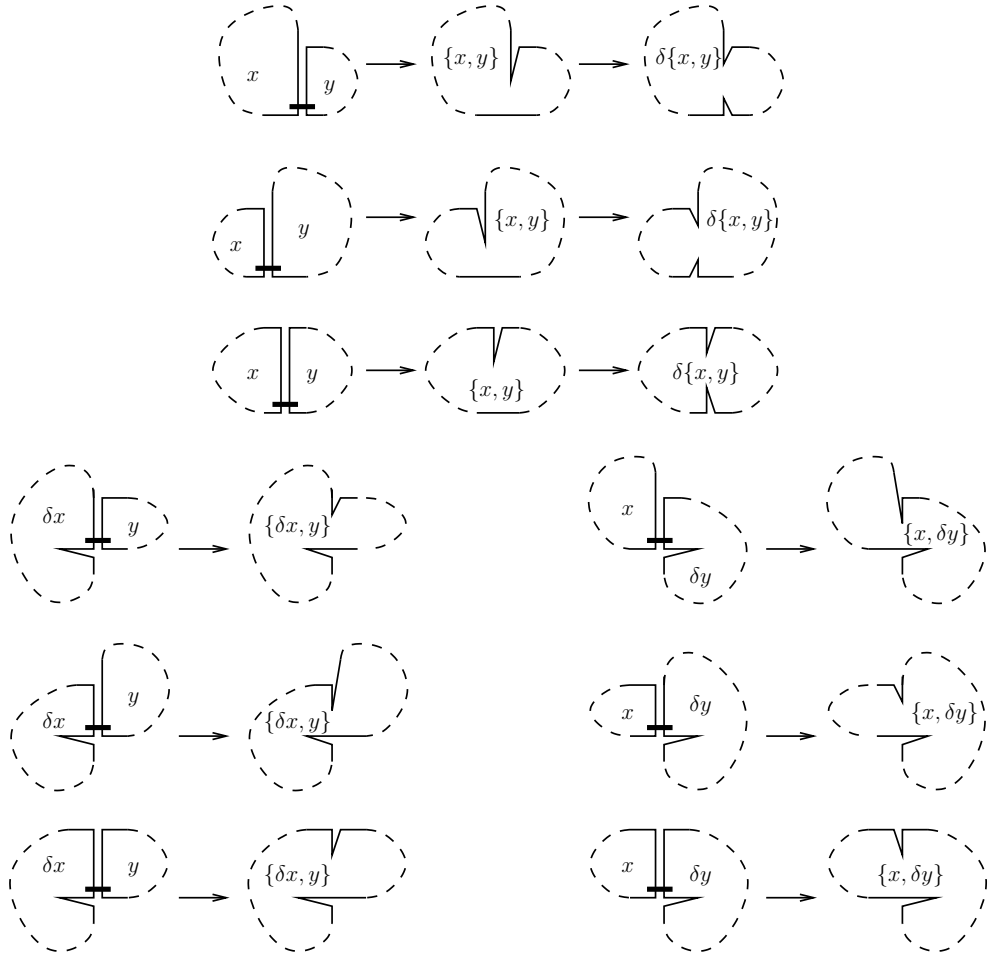}
}
\caption{
Exceptional terms in $f(x,y)$. The corners glued under the SFT bracket
are indicated by a heavy bar. This figure only shows terms
where the $x$ corner lies
counterclockwise from the $y$ corner at the gluing; there is a
corresponding set of terms where the $x$ corner lies clockwise from
the $y$ corner, obtained by reflecting each of the pictures.}
\label{fig:delta2}
\end{figure}

\begin{figure}
\centerline{
\includegraphics[height=3in]{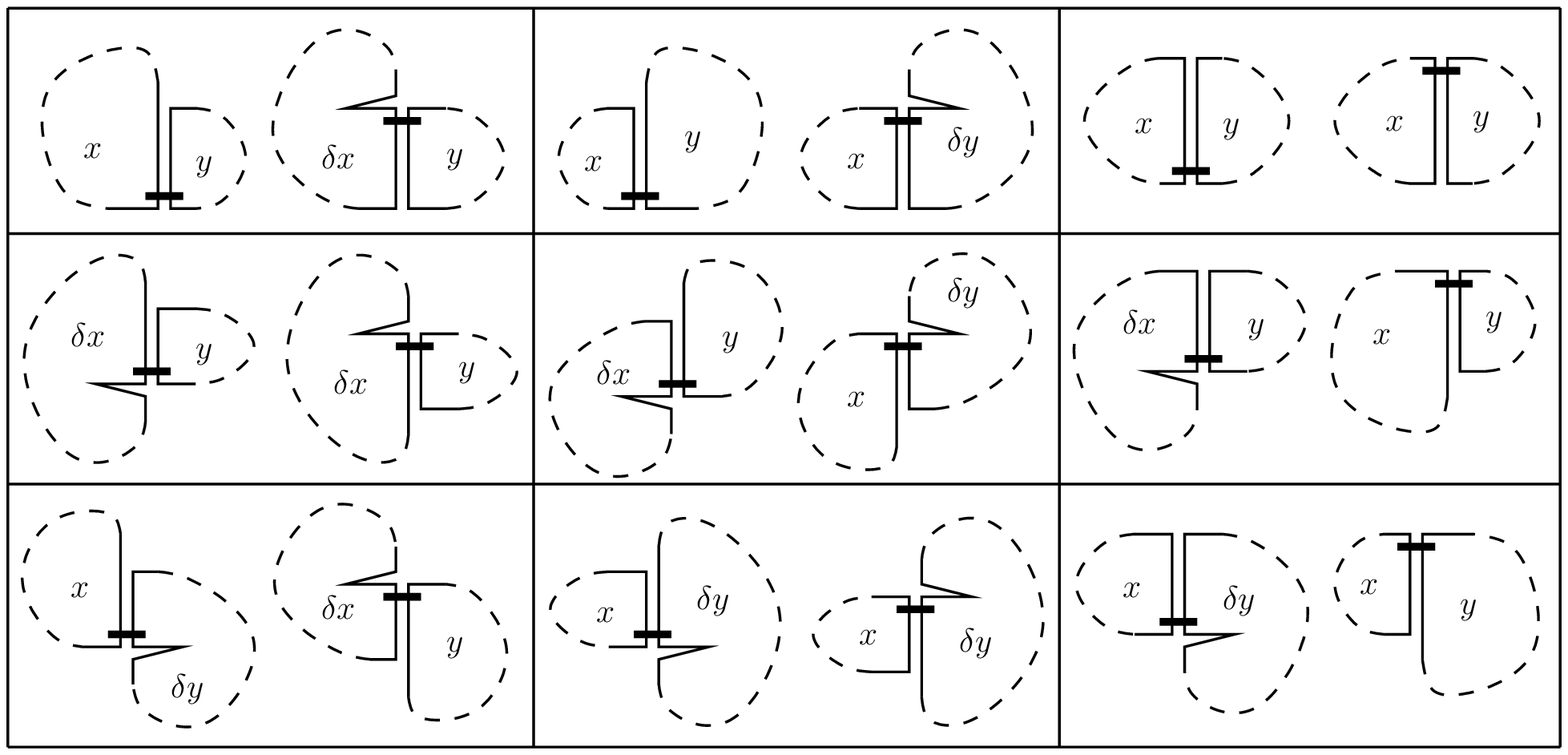}
}
\caption{Pairing the exceptional terms in $f(x,y)$. Bars indicate
  glued corners.}
\label{fig:delta3}
\end{figure}

These ``exceptional terms'' in $f(x,y)$ can be depicted as in
Figure~\ref{fig:delta2}, where in each schematic picture $x$ and $y$ are oriented counterclockwise. Note that in each case, a quadrant of $x$ or
$\delta x$ is glued to a quadrant of $y$ or $\delta
y$. Figure~\ref{fig:delta2} only shows the terms where the $x$
quadrant lies counterclockwise from the $y$ quadrant. There is an
analogous set of exceptional terms where the $x$ quadrant lies
clockwise from the $y$ quadrant.
Furthermore, there is a one-to-one correspondence between
``counterclockwise'' terms and ``clockwise'' terms; see
Figure~\ref{fig:delta3}. In $\cyc$, the terms under the one-to-one
correspondence cancel pairwise in $f(x,y)$, and (\ref{delta:4})
follows (up to sign, which will be more carefully considered below).

For (\ref{delta:5}), the same cancellation holds, but
we need to examine the position of the
base point on the based broken closed string $y$. If $x$ and $y$ do
not overlap (share a segment), then $\bullet(x) = f(x,y) = 0$ since
none of the exceptional terms exist.

\begin{figure}
\centerline{
\includegraphics[height=1in]{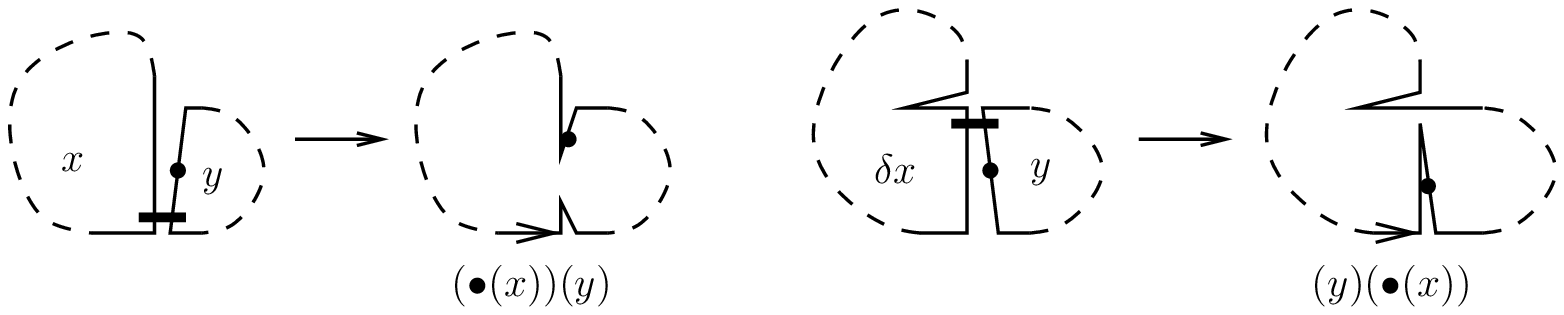}
}
\caption{
Two exceptional terms in $f(x,y)$ for $y\in\A$.
}
\label{fig:delta4}
\end{figure}

Otherwise, assume for simplicity
that $x$ and $y$ share exactly one segment, and in particular $x$
passes through the base point exactly once; the argument is similar in
the general case. It follows that there are exactly two exceptional
terms contributing to $f(x,y)$, and they are paired under the
correspondence of Figure~\ref{fig:delta3}. (The configuration of $x$
and $y$ determines which particular pair from Figure~\ref{fig:delta3}
appears.)

First work mod $2$.
If the base point does not lie on the shared segment between $x$ and
$y$, then $\bullet(x) =
f(x,y) = 0$ since the two exceptional terms cancel in $\A$. If
the base point lies on the shared segment, then one exceptional term
contributes $(\bullet(x))(y)$ and the other $(y)(\bullet(x))$, because
the base point is positioned differently on the glued broken closed
string depending on which gluing is used. See
Figure~\ref{fig:delta4}. This completes the proof of (\ref{delta:5})
mod $2$.

\begin{figure}
\centerline{
\includegraphics[height=1.25in]{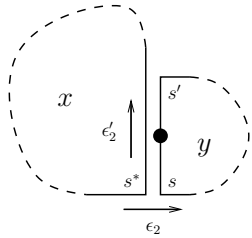}
}
\caption{
Labeling corners and orientations for $x,y$.
}
\label{fig:delta5}
\end{figure}

To establish (\ref{delta:5}) over $\Z$, we just need to check signs
for each of the nine pairs depicted in Figure~\ref{fig:delta3}. This
is completely straightforward but somewhat tedious; we do one sample
sign calculation and leave the rest to the interested reader. Suppose $x,y$
are as depicted in Figure~\ref{fig:delta4}. Label the corners of $x,y$
as shown in Figure~\ref{fig:delta5}. There are words $w_1,w_2$ such that
$x = [w_1 s^*]$, the image in $\A\tocyc$ of $w_1s^*$, and $y = s w_2
s'$. Also let $\epsilon_2,\epsilon_2'$ be the signs depicted in Figure~\ref{fig:delta5}.

The relevant term in $\sftb{x,y}$ is $\sftb{s^*,s} w_1w_2s'$,
where as usual $\sftb{s^*,s}$ is $+1$ if $s$ is a $q$, $-1$ if $s$ is
a $p$. It follows that the relevant term in $\delta\sftb{x,y}$ is
\[
\left((-1)^{|w_1|} \epsilon_2\right) w_1s^*sw_2s' = \left((-1)^{|x|}
  \epsilon_2'\right) x y = (\bullet(x))y,
\]
where the first equality follows from the fact that $\epsilon_2 =
(-1)^{|s^*|} \epsilon_2'$. On the other hand,
the relevant term in $\delta x$ is $\left((-1)^{|x|} \sftb{s'^*,s'}
\epsilon_2'\right) w_1 s^* s'^* s'$, which is equal in $\A\tocyc$ to
$\left((-1)^{|x|+|x||s'|} \sftb{s'^*,s'} \epsilon_2'\right) s' w_1 s^* s'^*$;
thus the relevant term in $\sftb{\delta x,y}$ is
\[
\left((-1)^{|x|+|x||s'|+|x||s w_2|} \epsilon_2'\right) s w_2 s' w_1 s^*
= (-1)^{|x||y|} y(\bullet(x)).
\]
Combining the terms in $\delta\sftb{x,y}$ and $-\sftb{\delta x,y}$
contributes $[\bullet(x),y]$ to $f(x,y)$, as desired.
\end{proof}

The string differential $\dstr$ from Section~\ref{ssec:comb} was
defined to satisfy the following result.

\begin{proposition}
On $\alg$, we have $\dstr = \delta$.
\label{prop:dstrdelta}
\end{proposition}

\begin{proof}
Since $\delta$ is a derivation by Proposition~\ref{prop:delta}(\ref{delta:2}), it suffices to show that $\dstr(s) = \delta(s)$ when
$s$ is a generator of $\A$.

\begin{figure}
\centerline{
\includegraphics[height=1.5in]{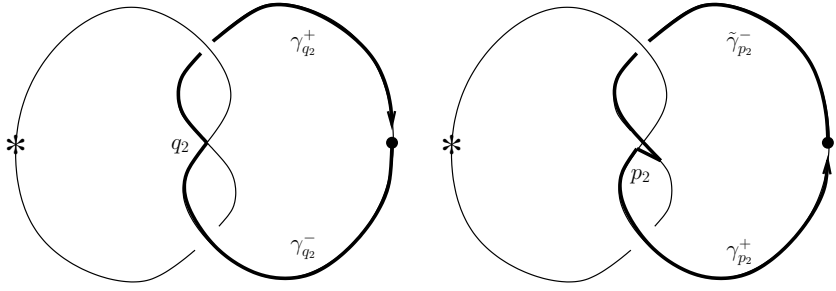}
}
\caption{
Based broken closed strings $\gamma_{q_2} = \gamma_{q_2}^- \cup
\gamma_{q_2}^+$ (left) and $\gamma_{p_2} = \tilde{\gamma}_{p_2}^- \cup
\gamma_{p_2}^+$ (right) for the knot $\Lambda_0$.
}
\label{fig:capping}
\end{figure}

If $s$ is a $p$ or $q$, we can use the paths $\gamma_s^{\pm}$ from
Section~\ref{ssec:comb} to define a based broken closed string
$\gamma_s$ with holomorphic corners such that $w(\gamma_s) = s$. More
precisely, if $s$ has holomorphic capping paths, then define $\gamma_s
= \gamma_s^- \cup \gamma_s^+$; if $s$ has antiholomorphic capping
paths, then define $\gamma_s$ to be a perturbation of $\gamma_s^- \cup
\gamma_s^+$ to have a holomorphic corner at $s$ (explicitly, let
$\gamma_s^-$ run past the crossing and return, and then join
$\gamma_s^+$ to it). See Figure~\ref{fig:capping}. If $s=t^{\pm 1}$, define paths $\gamma_t$ to run along
$\Lambda$ once, and $\gamma_{t^{-1}}$ to run along $\Lambda$ once with
the reverse orientation.

In all cases, it is now trivial to check, using $\gamma_s$, that
$\delta(s)$ from Definition~\ref{def:deltabcs} agrees with $\dstr(s)$
from Definition~\ref{def:dstr}.
\end{proof}

\subsection{The Hamiltonian and the LSFT algebra}
\label{ssec:ham}

Having introduced the SFT bracket and the $\delta$ map, we are now in
a position to redefine the LSFT algebra $(\alg,\d)$ in terms of strings.
First we introduce the Hamiltonian $h\in\cyc$ counting rigid
holomorphic disks with boundary on $\Lambda\times\R$.

Let $f$ be an
immersed disk in $\R^2$ with boundary on $\pi_{xy}(\Lambda)$ and
convex corners. More precisely, in the language of
Definition~\ref{def:disks}, $f\in\Delta(s_1,\dots,s_k)$
for some $s_1,\dots,s_k$; recall that $f$ is equally well an element of
$\Delta(s_2,\dots,s_k,s_1)$ and other cyclic permutations as well. The
boundary of $f$ is a broken closed string in $\Lambda$ with
corresponding word $s_2\cdots s_ks_1 \in\A$ (note that $s_1$ appears
\textit{last} in this word), and we can also associate a
sign $\epsilon(f;s_1)$ to $f$, as defined in
Section~\ref{ssec:comb}.

Define $\tilde{w}(f) \in \A\tocyc$ to be the
image in $\A\tocyc$ of $\epsilon(f;s_1) s_2\cdots s_ks_1$. The key point
is the following.

\begin{lemma}
\label{lem:hamsign}
The element $\tilde{w}(f)$ of $\A\tocyc$
depends only on the disk $f$ and not on which
puncture is labeled $s_1$; that is, $\tilde{w}(f)$ is independent of
whether $f$ is viewed as an element of $\Delta(s_1,\dots,s_k)$,
$\Delta(s_2,\dots,s_k,s_1),$ or any other cyclic permutation.
\end{lemma}

\begin{proof}
This is clear mod $2$. To check signs, suppose
$f\in\Delta(s_1,\dots,s_k) = \Delta(s_2,\dots,s_k,s_1)$. If we view
$f$ as an element of
$\Delta(s_1,\dots,s_k)$, then $\tilde{w}(f) = \epsilon(f;s_1)
s_2\cdots s_ks_1$, while if we view $f$ as an element of
$\Delta(s_2,\dots,s_k,s_1)$, then $\tilde{w}(f) = \epsilon(f;s_2)
s_3\cdots s_k s_1 s_2$. But we have
$\epsilon(f;s_1) \epsilon(f;s_2) = \epsilon'(f;s_1) \epsilon'(f;s_2) =
(-1)^{|s_2|}$,
where $\epsilon'(f;s_1),\epsilon'(f;s_2)$ are the signs shown in
Figure~\ref{fig:hamsign} and the second equality follows from
Lemma~\ref{lem:degsign}. Since $|f| = -2$ by Lemma~\ref{lem:hdegree}
(cf.\ Lemma~\ref{lem:ham} below),
\[
\epsilon(f;s_1) s_2\cdots s_ks_1 = \epsilon(f;s_1) (-1)^{|s_2|}
s_3\cdots s_ks_1s_2 = \epsilon(f;s_2) s_3\cdots s_ks_1s_2
\]
in $\A\tocyc$, and the lemma follows.
\end{proof}

\begin{definition}
The Hamiltonian $h\in\cyc$ is the sum of $\tilde{w}(f)$ over all
immersed disks $f$ in some $\Delta(s_1,\dots,s_k)$ for all possible
$k\geq 1$ and all possible $s_1,\dots,s_k$. (Here we mod out by cyclic permutations and count each immersed disk once; that is, we count $f\in\Delta(s_1,\dots,s_k)$ once and not $k$ times.)
\end{definition}

\begin{figure}
\centerline{
\includegraphics[width=5in]{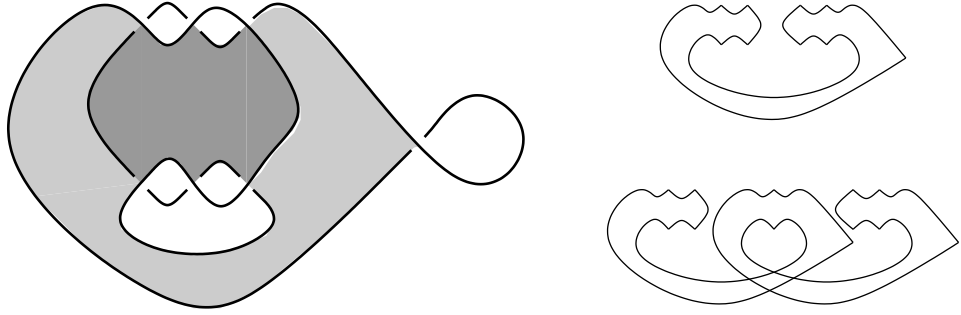}
}
\caption{
The $xy$ projection of a Legendrian knot for which $h$ is infinite.
One immersed disk is
indicated in the diagram on the left (the darker-shaded region is
covered twice), and heuristically redrawn as the
diagram on the top right (pulled apart from itself for clarity). This
can be extended by wrapping more times
around the bottom section of the projection, to produce an infinite
family of contributions to $h$, the next of which is drawn on the
bottom right.
This knot is isotopic through restricted Reidemeister II moves to the
standard Legendrian unknot.
}
\label{fig:hinf}
\end{figure}

It is entirely possible that $h$ is an infinite sum; see
Figure~\ref{fig:hinf} for an example. We do however have the following
result.

\begin{lemma}
\label{lem:ham}
The Hamiltonian $h$ has degree $-2$ and is an element of $\F^1\!\cyc$.
\end{lemma}

\begin{proof}
The fact that $h$ has degree $-2$ can be proved in
the same way as Lemma~\ref{lem:hdegree}. To show that $h\in\F^1\!\cyc$,
we claim that all terms in $h$ contain at least one $p$, and that only
finitely many terms contain at most $k$ $p$'s for any $k$.
The first part is evident from Lemma~\ref{lem:stokes}. For the second
part, Lemma~\ref{lem:stokes} implies that there are only finitely many
nonempty moduli spaces of disks $\Delta(s_1,\dots,s_l)$ for which $k$
of the $s_i$'s are $p$'s. But for fixed $l$ and $s_1,\dots,s_l$, the
set $\Delta(s_1,\dots,s_l)$ is finite by a standard argument
given in \cite{Ch02}.
\end{proof}

We now have the following result.

\begin{proposition}[Quantum master equation]
$\delta h + \frac{1}{2} \sftb{h,h} = 0$.
\label{prop:QME}
\end{proposition}

\begin{remark}
Despite the presence of $\frac{1}{2}$ in the
statement of Proposition~\ref{prop:QME}, the result can be interpreted
over $\Z$ or even $\Z/2$. To do this, rewrite $\sftb{x,y}$ as a
difference of two terms, $(x\to y) - (y\to x)$, where $(x\to y)$
counts terms in $\sftb{x,y}$ where a $p$ in $x$ is glued to the
corresponding $q$ in $y$ and similarly for $(y\to x)$. Then we can
write $(h\to h)$ instead of $\frac{1}{2} \sftb{h,h}$.
\label{rem:hh}
\end{remark}

\begin{figure}[b]
\centerline{
\includegraphics[height=1.8in]{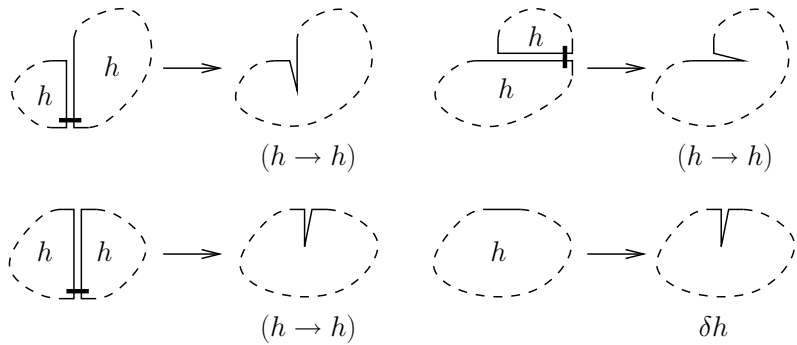}
}
\caption{
Paired contributions to $\delta h + (h \to h)$. Top row: an obtuse
disk contributes to two canceling terms in $(h \to h)$. Bottom row:
an immersed disk contributes to terms in $(h\to h)$ and $\delta h$,
which also cancel.
}
\label{fig:hh}
\end{figure}

\begin{proof}[Proof of Proposition~\ref{prop:QME}]
We wish to show that $\delta h = -(h\to h)$ in the notation of
Remark~\ref{rem:hh}. First argue mod $2$.
As in the standard proof of $d^2=0$ in Chekanov \cite{Ch02}, most of
the terms in $(h\to h)$ cancel pairwise. Terms in $(h \to h)$
correspond to gluing two immersed disks at a corner; near this corner,
the two disks overlap on an edge. If the overlapping edges are not
identical, then the result is an ``obtuse disk'' with one concave
corner, and this obtuse disk appears twice in $(h\to h)$. See the top
line of Figure~\ref{fig:hh}. If the overlapping edges are identical,
then the glued disk is also an immersed disk, and the contribution of
the glued disk to $(h\to h)$ is canceled by the contribution of the
immersed disk to $\delta h$. See the bottom line of Figure~\ref{fig:hh}.

\begin{figure}
\centerline{
\includegraphics[height=1in]{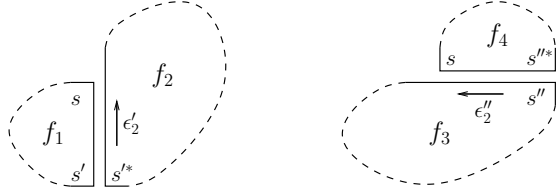}
}
\caption{
A closer look at obtuse disks.
}
\label{fig:hhsigns}
\end{figure}

As usual, to complete the proof, we need to compute signs. We
claim that the two obtuse disks (the top row of Figure~\ref{fig:hh})
give canceling contributions to $(h\to h)$; a similar calculation
shows that the bottom row of Figure~\ref{fig:hh} gives canceling
contributions to $(h\to h) + \delta h$. Consider the four disks
$f_1,f_2,f_3,f_4$ shown in Figure~\ref{fig:hhsigns}. The contribution
of, e.g., $f_1,f_2$ to $(h\to h)$ is either
$\{\tilde{w}(f_1),\tilde{w}(f_2)\}$ or
$\{\tilde{w}(f_2),\tilde{w}(f_1)\}$
depending on which of $f_1,f_2$ contains the $p$ and which the $q$,
but these two quantities are equal since $h$ has degree $-2$. It thus
suffices to show that the contributions of
$\{\tilde{w}(f_1),\tilde{w}(f_2)\}$ and
$\{\tilde{w}(f_3),\tilde{w}(f_4)\}$ to $\sftb{h,h}$ have opposite
sign.

We have
\[
\tilde{w}(f_1) = \left( \epsilon_s(f_1) \epsilon_{s'}(f_1)
  \epsilon_{w_1}(f_1) \epsilon_2' \right) s w_1 s'
\]
where $w_1$ is some word,
$\epsilon_s(f_1),\epsilon_{s'}(f_1)$ are the orientation signs
for corners $s,s'$ in $f_1$, $\epsilon_{w_1}(f_1)$ is the product of
orientation signs over all other corners of $f_1$ (i.e., the corners
corresponding to $w_1$), and $\epsilon_2'$
is the sign depicted in Figure~\ref{fig:hhsigns} (as usual, relative
to the knot orientation). Similarly, we have
\begin{align*}
\tilde{w}(f_2) &= \left( -\epsilon_{s'^*}(f_2) \epsilon_{w_2}(f_2)
  \epsilon_2' \right) s'^* w_2 \\
\tilde{w}(f_3) &= \left( \epsilon_{s''}(f_3) \epsilon_{w_3}(f_3)
  \epsilon_2'' \right) w_3 s'' \\
\tilde{w}(f_4) &= \left( -\epsilon_{s''^*}(f_4) \epsilon_s(f_4)
  \epsilon_{w_4}(f_4) \right) s''^* w_4 s.
\end{align*}
Gluing $s'$ in $f_1$ to $s'^*$ in $f_2$ yields a contribution of
\[
\left( -\epsilon_s(f_1)\epsilon_{s'}(f_1) \epsilon_{s'^*}(f_2)
  \epsilon_{w_1}(f_1) \epsilon_{w_2}(f_2) \sftb{s',s'^*} \right)
sw_1w_2 = \left(-\epsilon_s(f_1) \epsilon_{w_1}(f_1) \epsilon_{w_2}(f_2)
\right) sw_1w_2
\]
to $\sftb{\tilde{w}(f_1),\tilde{w}(f_2)}$ by Lemma~\ref{lem:orsigns};
similarly, gluing $s''$ in $f_3$ to $s''^*$ in $f_4$
yields a contribution of
$\left( -\epsilon_s(f_4) \epsilon_{w_3}(f_3)
  \epsilon_{w_4}(f_4) \right) w_3w_4s$
to $\sftb{\tilde{w}(f_3),\tilde{w}(f_4)}$. But $w_1w_2=w_3w_4$ and hence
$sw_1w_2 = w_3w_4s$ in $\A\tocyc$, while
$\epsilon_{w_1}(f_1)\epsilon_{w_2}(f_2) =
\epsilon_{w_3}(f_3)\epsilon_{w_4}(f_4)$ and hence
\[
\left(-\epsilon_s(f_1) \epsilon_{w_1}(f_1) \epsilon_{w_2}(f_2)
\right)= - \left( -\epsilon_s(f_4) \epsilon_{w_3}(f_3)
  \epsilon_{w_4}(f_4) \right)
\]
since $\epsilon_s(f_1) = -\epsilon_s(f_4)$ by Lemma~\ref{lem:orsigns}.
This shows that the obtuse disks give canceling contributions to
$(h\to h)$, as desired.
\end{proof}

The following result, the string version of the ``$\d^2=0$ result''
Theorem~\ref{thm:d2}, implies Theorem~\ref{thm:d2}.

\begin{proposition}
Define $\d\co\alg\to\alg$ by
\[
\d(x) = \sftb{h,x} + \delta x.
\]
Then $\sftb{h,x}$, $\delta x$, and $\d(x)$ coincide with $\dsft(x)$,
$\dstr(x)$, and $\d(x)$ as defined in Section~\ref{ssec:comb}, and
$(\alg,\d)$ is an LSFT algebra with $F_{\d} = \bullet(h)$.
\label{prop:d2}
\end{proposition}

\begin{proof}
We have already seen in Proposition~\ref{prop:dstrdelta} that $\delta
= \dstr$. The fact that $\sftb{h,\cdot} = \dsft$ follows from a direct
inspection of the definitions of $h$, $\dsft$, and the SFT bracket.

It remains to show that
$\d^2 x = [\bullet(h),x]$. But
Proposition~\ref{prop:sftb} implies that $\sftb{h,\sftb{h,x}} =
\sftb{\frac{1}{2}\sftb{h,h},x}$, whence by Propositions~\ref{prop:delta}
and~\ref{prop:QME},
\begin{align*}
d^2 x &= \sftb{h,\sftb{h,x}} + \sftb{h,\delta x} + \delta\sftb{h,x} +
\delta^2 x \\
&= - \sftb{\delta h,x} + \sftb{h,\delta x} + \delta\sftb{h,x} \\
&= [\bullet(h),x],
\end{align*}
as desired.
\end{proof}

\secspace

\section{Proof of Invariance}
\label{sec:invariance}

This section is devoted to the proof of the invariance result,
Theorem~\ref{thm:invariance}. The LSFT algebra structure is associated
to a Legendrian knot $\Lambda$ with two marked points
$\ast,\bullet$. Any two Legendrian-isotopic knots with marked points
can be related by a sequence of four basic moves: keeping the knot
fixed and sliding $\ast$ along it; keeping the knot fixed and sliding
$\bullet$ along it; changing the knot by a Reidemeister II move, while
keeping $\ast,\bullet$ fixed and away from the move; and changing the
knot by a Reidemeister III move, while keeping $\ast,\bullet$ fixed
and away from the move.

\begin{figure}[b]
\centerline{
\includegraphics[height=0.75in]{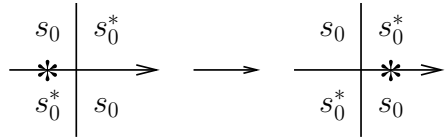}
}
\caption{
Changing the marked point $\ast$ by moving it through a crossing of
$\pi_{xy}(\Lambda)$, in the direction of the orientation of
$\Lambda$. There is an analogous diagram moving the marked point
$\bullet$ through the crossing.
}
\label{fig:markedptchange}
\end{figure}

Changing the marked points changes the LSFT algebra in a fairly
trivial way. One can readily check from the definitions that
moving $\ast$ across a crossing of $\pi_{xy}(\Lambda)$ labeled by $s_0,s_0^*$, as shown in Figure~\ref{fig:markedptchange}, has the effect
of replacing $s_0$ by $t^{-1}s_0$ and $s_0^*$ by $s_0^* t$,
and thus corresponds to a basis change. On
the other hand, moving $\bullet$ in the same way does not change $\dsft$ but does change $\dstr$ by
\[
\dstr'(s) = \dstr(s) - \sftb{s_0,s_0^*} [s_0s_0^*,s]
\]
for all generators $s$; this corresponds to a gauge change with
$z=-\sftb{s_0,s_0^*} s_0 s_0^*$, in the notation of Definition~\ref{def:gauge}.

The remaining moves, Reidemeister III and Reidemeister II, are
addressed in Sections~\ref{ssec:RIII} and \ref{ssec:RII},
respectively. These are essentially extensions of the invariance
arguments for the contact-homology differential graded algebra from \cite{Ch02}.

\subsection{Reidemeister III}
\label{ssec:RIII}

\begin{figure}
\centerline{
\includegraphics[height=1.25in]{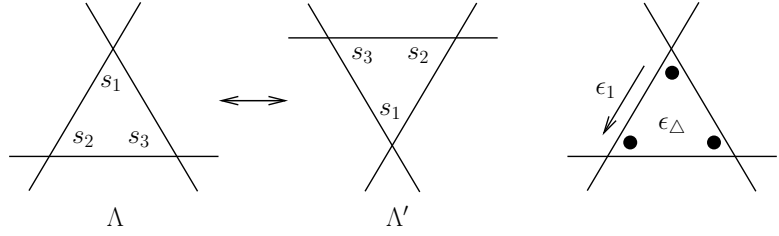}
}
\caption{
Reidemeister III move. The right-hand diagram defines signs
$\epsilon_\triangle, \epsilon_1$: $\epsilon_\triangle$
is the product of the orientation signs of the three marked corners,
while $\epsilon_1$ is the depicted orientation relative to the
orientation of $\Lambda$.}
\label{fig:ReidIII}
\end{figure}

Here we assume that $\Lambda$ and $\Lambda'$ are related by a
Reidemeister III move, as shown in Figure~\ref{fig:ReidIII}.
We may also assume without loss of generality that the points
$\ast,\bullet$ are not involved in the move and lie outside of the
local pictures. Let
$(\alg,\d)$ and $(\alg,\d')$ be the LSFT algebras associated to
$\Lambda$ and $\Lambda'$, respectively. Note that we have identified
the underlying tensor algebras by labeling the three relevant
crossings as shown; all other crossings are labeled identically for the
two pictures. Figure~\ref{fig:ReidIII} also defines two signs
$\epsilon_\triangle, \epsilon_1$.

Let $\phi$ be the basis change defined as follows:
\begin{align*}
s_1 &\mapsto s_1 + \epsilon \sftb{s_1^*,s_1} s_3^* s_2^* \\
s_2 &\mapsto s_2 + \epsilon \sftb{s_2^*,s_2} s_1^* s_3^* \\
s_3 &\mapsto s_3 + \epsilon \sftb{s_3^*,s_3} s_2^* s_1^* \\
t^{\pm 1} &\mapsto t^{\pm 1} \\
s &\mapsto s \qquad\qquad\qquad\qquad \text{ for $s \neq s_1,s_2,s_3$},
\end{align*}
where $\epsilon = \sftb{s_1^*,s_1} \sftb{s_2^*,s_2} \sftb{s_3^*,s_3}
\epsilon_\triangle$.
Note that $\phi$ preserves the filtration on $\alg$. Indeed, the only
way this would not hold would be if $s_1,s_2,s_3$ were all $p$'s, in
which case the move depicted in
Figure~\ref{fig:ReidIII} would not be a topological isotopy.

\begin{lemma}
If $x\in\A\tocyc$ and $y\in\A$, \label{lem:sftbphi}
then
\begin{equation}
\phi(\sftb{x,y}) = \sftb{\phi(x),\phi(y)}.
\label{eq:sftbphi}
\end{equation}
\end{lemma}

\noindent
We remark that Lemma~\ref{lem:sftbphi} does not hold if $\phi$ is an
arbitrary basis change.

\begin{proof}[Proof of Lemma~\ref{lem:sftbphi}]
By Proposition~\ref{prop:sftb}(\ref{it:sftb2}),
it suffices to establish (\ref{eq:sftbphi}) for $y=s$, where $s$ is
any $p$ or $q$. If $s \neq s_1,s_2,s_3$, then
(\ref{eq:sftbphi}) holds because $s$ does not interact with any
of the quadratic terms in $\phi$.
It remains to prove (\ref{eq:sftbphi}) for $y=s_1,s_2,s_3$. We may
assume without loss of generality that $y=s_1$ and further that $x$ is (the
cyclic quotient of) a word in $\A$.

Contributions to $\sftb{x,y}$ come from appearances of $s_1^*$ in $x$:
if $x = w_1 s_1^* w_2$, then $\sftb{x,y}$ contains the term
$\pm w_2 w_1$,
and $\phi(\sftb{x,y})$ contains the term $\pm \phi(w_2)\phi(w_1)$.  On the other
hand, contributions to $\sftb{\phi(x),\phi(y)} =
\sftb{\phi(x),s_1+\epsilon \sftb{s_1^*,s_1} s_3^*s_2^*}$ come from
appearances of any of
$s_1^*,s_2,s_3$ in $x$. The appearances of $s_1^*$ give a contribution
to $\sftb{\phi(x),\phi(y)}$ exactly equal to the corresponding
contribution to $\phi\sftb{x,y}$. Any appearance of $s_2$ gives two
canceling contributions to $\sftb{\phi(x),\phi(y)}$: if $x = w_1 s_2
w_2$, then $\phi(x) = \phi(w_1) s_2 \phi(w_2) + \epsilon
\sftb{s_2^*,s_2} \phi(w_1) s_1^* s_3^* \phi(w_2)$,
and the contribution to $\sftb{\phi(x),s_1+\epsilon \sftb{s_1^*,s_1}
  s_3^*s_2^*}$ is
\begin{multline*}
(-1)^{|w_1s_2||w_2|} \left( \epsilon \sftb{s_1^*,s_1} \sftb{s_2,s_2^*}
  (-1)^{(|w_1s_2w_2|+1)|s_3^*|} \right. \\
\left. + \epsilon \sftb{s_2^*,s_2}
  \sftb{s_1^*,s_1} (-1)^{|s_3^*||w_2w_1s_1^*|} \right) s_3^* \phi(w_2)
\phi(w_1) = 0,
\end{multline*}
where the cancellation occurs since $|s_2| = |s_3^*s_1^*|$.
Any appearance of $s_3$ in $x$ similarly gives two canceling
contributions to $\sftb{\phi(x),\phi(y)}$.
It follows that (\ref{eq:sftbphi})
holds for $y = s_1$, as desired.
\end{proof}

\begin{lemma}
If $\gamma$ is a path in $\Lambda$ whose endpoints do not coincide
with any of the endpoints of $s_1,s_2,s_3$, then $\phi(\delta(\gamma))
= \delta(\gamma)$.
\label{lem:RIII1}
\end{lemma}

\begin{proof}
The hypothesis of the lemma implies that $s_1,s_2,s_3$ only appear in
$\delta(\gamma)$ in pairs, namely $\sftb{s_1,s_1^*} s_1^*s_1+
\sftb{s_2,s_2^*} s_2s_2^*$ and its two cyclic permutations
(where we cyclically permute the indices $1,2,3$).
But $\phi$ preserves each of these sums.
\end{proof}

\begin{lemma}
If $s$ is a generator of $\alg$ not equal to
$s_1,s_2,s_3,s_1^*,s_2^*,s_3^*$, then $\delta'(s) = \delta(s)$. Also,
\begin{align*}
\delta'(s_1) &= \delta(s_1) - (\epsilon_1 (-1)^{|s_1|} \sftb{s_2,s_2^*})
s_1 s_2 s_2^* - (\epsilon_1 (-1)^{|s_1|} \sftb{s_3^*,s_3}) s_3^* s_3 s_1 \\
\delta'(s_1^*) &= \delta(s_1^*) + (\epsilon_1 \sftb{s_2,s_2^*}) s_2
s_2^* s_1^* - (\epsilon_1 \sftb{s_3^*,s_3}) s_1^* s_3^* s_3,
\end{align*}
with corresponding formulas for
$\delta'(s_2),\delta'(s_3),\delta'(s_2^*),\delta'(s_3^*)$ (permute the
indices cyclically).
\label{lem:RIIIdelta}
\end{lemma}

\begin{proof}
The first statement is clear by Lemma~\ref{lem:RIII1}.
The rest follows from the definition of $\delta,\delta'$ and an
examination of how capping paths change under the Reidemeister III move.
\end{proof}

For the next lemma, note that the triangle in the Reidemeister III
move contributes the term $(\epsilon_\triangle \epsilon_1) s_2s_3s_1$ to
both $h$ and $h'$ (for the latter, use Lemma~\ref{lem:orsigns}).

\begin{lemma}
\label{lem:RIII2}
Write $h_0 = h - (\epsilon_\triangle \epsilon_1) s_2s_3s_1$
and $h_0' = h'- (\epsilon_\triangle \epsilon_1) s_2s_3s_1$; then
$h_0' = \phi(h_0)$.
\end{lemma}

\begin{figure}
\centerline{
\includegraphics[height=2.5in]{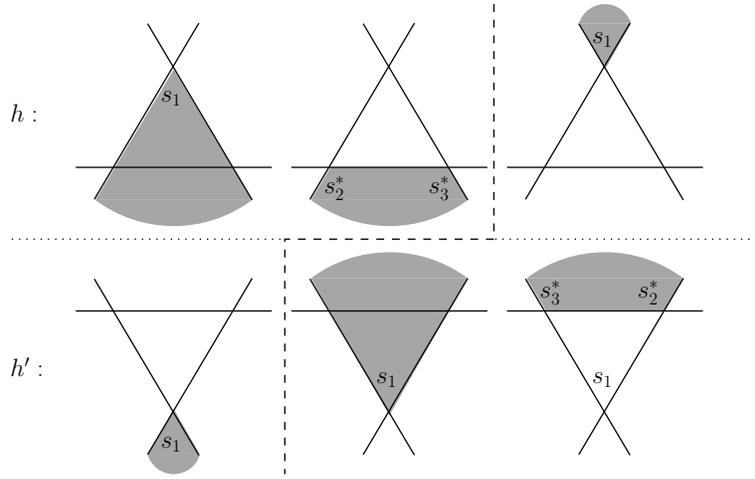}
}
\caption{
Contributions to $h$ and $h'$ for the Reidemeister III move.
}
\label{fig:RIIIdiff}
\end{figure}

\begin{proof}
This is a standard argument along the lines of Chekanov \cite{Ch02}. Let
$\triangle,\triangle'$ denote the triangles bounded by $s_1,s_2,s_3$ in
$\pi_{xy}(\Lambda),\pi_{xy}(\Lambda')$, respectively. Disks that
contribute to $h,h'$ and contain $s_1$ fall into two categories,
depending on whether they contain $\triangle,\triangle'$ or not. To a disk
in $h'$ with a corner at $s_1$ and not containing $\triangle'$ (the left
picture in the bottom row of Figure~\ref{fig:RIIIdiff}), there are
two corresponding disks in $h$, one with a corner at $s_1$ and
containing $\triangle$, the other with corners at $s_3^*$ and $s_2^*$
(the two left pictures in the top row). To a disk in $h$ with a corner
at $s_1$ and not containing $\triangle$ (right picture, top row), there
are two corresponding disks in $h'$, one with a corner at $s_1$ and
containing $\triangle'$, the other with corners at $s_3^*$ and $s_2^*$
(two right pictures, bottom row). Similar correspondences occur for
disks in $h,h'$ with corners at $s_2$ or $s_3$.
That the signs work out follows easily from the definition of
$\epsilon$ and Lemma~\ref{lem:orsigns}.
\end{proof}

\begin{proposition}
Under a Reidemeister III move, the LSFT algebra changes by a change of
basis. More precisely, $d' = \phi \circ \d \circ \phi^{-1}$.
\end{proposition}

\begin{proof}
It suffices to show that
\begin{equation}
\phi \d(s) = \d' \phi (s)
\label{eq:RIII}
\end{equation}
for $s$ any generator of $\alg$.


\noindent \textsc{Case 1:} $s \neq
s_1,s_2,s_3,s_1^*,s_2^*,s_3^*$.
By Lemma~\ref{lem:RIIIdelta},
$\phi(\delta(s)) = \delta(s) = \delta'(s)$,
while by Lemmas~\ref{lem:sftbphi}
and~\ref{lem:RIII2}, $\phi\sftb{h_0,s} = \sftb{h_0',s}$. It follows that
\[
\phi d(s) = \phi\sftb{h_0+(\epsilon_\triangle \epsilon_1)
  s_2s_3s_1,s}+
\phi\delta(s)
= \phi\sftb{h_0,s} + \delta'(s) = d'\phi(s),
\]
and (\ref{eq:RIII}) holds.

\noindent \textsc{Case 2:} $s=s_1^*,s_2^*,s_3^*$. By symmetry, we
may assume that $s = s_1^*$. Then
\begin{align*}
\phi \d s_1^* &= \phi\sftb{h_0,s_1^*} + \epsilon_\triangle \epsilon_1 \phi\sftb{s_2s_3s_1,s_1^*} +
\phi\delta(s_1^*) \\
&= \sftb{h_0',s_1^*} + \epsilon_\triangle \epsilon_1 \sftb{s_1,s_1^*} (s_2+ \epsilon \sftb{s_2^*,s_2} s_1^*s_3^*)(s_3+\epsilon \sftb{s_3^*,s_3} s_2^*s_1^*) + \phi\delta(s_1^*)
\end{align*}
and
\begin{align*}
\d'\phi s_1^* = \d' s_1^* &= \sftb{h_0',s_1^*} +
\epsilon_\triangle \epsilon_1 \sftb{s_2s_3s_1,s_1^*} + \delta'(s_1^*) \\
&= \sftb{h_0',s_1^*} + \epsilon_\triangle \epsilon_1 \sftb{s_1,s_1^*}
s_2s_3 + \delta(s_1^*) + \epsilon_1 \sftb{s_2,s_2^*} s_2 s_2^* s_1^* -
\epsilon_1 \sftb{s_3^*,s_3} s_1^* s_3^* s_3;
\end{align*}
thus to establish (\ref{eq:RIII}) for $s=s_1^*$, it suffices to
show that
\[
\phi\delta(s_1^*) - \delta(s_1^*) = \epsilon_1 \epsilon s_1^* s_3^* s_2^* s_1^*.
\]
Because of the form of $\phi$, one can disregard all terms in
$\delta(s_1^*)$ except those of the
form $s_1^* s_1 s_1^*$, $s_1^* s_3^* s_3$, $s_2 s_2^* s_1^*$. It is
straightforward to check that the total contribution of these terms to
$\phi\delta(s_1^*) - \delta(s_1^*)$ is precisely $\epsilon_1 \epsilon
s_1^* s_3^* s_2^* s_1^*$ for any of the four possible configurations of
$\gamma_{s_1^*}$ near $s_1^*$.

\noindent \textsc{Case 3:} $s=s_1,s_2,s_3$. We may assume
that $s=s_1$. Now
\[
\phi \d(s_1) = \phi\sftb{h,s_1} + \phi\delta(s_1) =
\phi\sftb{h_0,s_1} + \phi\delta(s_1)
\]
and
\begin{align*}
\d' \phi(s_1) &=
\sftb{h_0',\phi(s_1)} + \epsilon_\triangle \epsilon_1 \sftb{s_2s_3s_1,s_1+\epsilon \sftb{s_1^*,s_1} s_3^* s_2^*} + \delta'(s_1+\epsilon \sftb{s_1^*,s_1} s_3^* s_2^*) \\
&= \phi\sftb{h_0,s_1} +\delta'(s_1+\epsilon \sftb{s_1^*,s_1} s_3^* s_2^*) \\
& \quad + \epsilon_1 \sftb{s_2^*,s_2}\sftb{s_3^*,s_3} \left( (-1)^{|s_1|} \sftb{s_3,s_3^*} s_1s_2s_2^* + (-1)^{|s_2|+|s_3^*|} \sftb{s_2,s_2^*} s_3^*s_3s_1 \right) \\
&= \phi\sftb{h_0,s_1} + \delta(s_1) + \epsilon \sftb{s_1^*,s_1} \delta'(s_3^*s_2^*)
%
\end{align*}
by Lemmas~\ref{lem:sftbphi}, \ref{lem:RIIIdelta}, and \ref{lem:RIII2};
thus to establish (\ref{eq:RIII}) for $s=s_1$, it suffices to show
that
\begin{equation}
\phi\delta(s_1) - \delta(s_1) = \epsilon \sftb{s_1^*,s_1} \delta'(s_3^* s_2^*). \label{eq:s1}
\end{equation}
By Lemma~\ref{lem:RIII1}, (\ref{eq:s1}) simply states that replacing each appearance of $s_1$ in $\delta(s_1)$ by $s_3^*s_2^*$ results in $\delta'(s_3^* s_2^*)$. But given a based broken closed string in $\Lambda$ with a single holomorphic corner at $s_1$, a small perturbation yields a based broken closed string in $\Lambda'$ whose word is $s_3^*s_2^*$; the correspondence between these strings yields (\ref{eq:s1}).
\end{proof}

\subsection{Reidemeister II}
\label{ssec:RII}

\begin{figure}[b]
\centerline{
\includegraphics[height=0.9in]{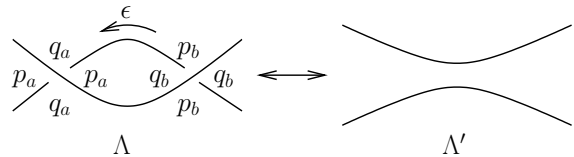}
}
\caption{
Reidemeister II move. The two crossings in $\pi_{xy}(\Lambda)$ are
given labels $p_a,q_a,p_b,q_b$ in their quadrants, as shown. A sign
$\epsilon$ is also shown.
}
\label{fig:ReidII}
\end{figure}

Here we assume that $\Lambda$ and $\Lambda'$ are related by a
Reidemeister II move, as shown in Figure~\ref{fig:ReidII}.
At some point it will become important that the move is a restricted
Reidemeister II move; we will indicate where we use this fact in the proof.

As in the
Reidemeister III case, assume that the points $\ast,\bullet$ are not
involved in the move and lie outside of the local pictures. Let
$(\alg,\d)$ and $(\alg',\d')$ be the LSFT algebras associated to $\Lambda$
and $\Lambda'$, respectively. View the algebra $\alg$ as a stabilization
of $\alg'$ by adding four generators $p_a,q_a,p_b,q_b$ corresponding to
the two new crossings in $\Lambda$. Then we can extend $\d'$ to $\alg$
by setting
\begin{align*}
\d'(q_a) &= q_b \\
\d'(q_b) &= [F_{\d'},q_a] \\
\d'(p_b) &= p_a \\
\d'(p_a) &= [F_{\d'},p_b].
\end{align*}
Note that this makes $(\alg,\d')$ an LSFT algebra and a stabilization of
$(\alg',\d')$.

We claim that $(\alg,\d)$ and $(\alg,\d')$ are related by a basis change;
this will prove invariance under restricted Reidemeister II.

The bigon in Figure~\ref{fig:ReidII} contributes the term $\epsilon
p_aq_b$ to the Hamiltonian for $\Lambda$, where $\epsilon$ is the sign
depicted in Figure~\ref{fig:ReidII}; this then contributes $-\epsilon
p_a,\epsilon (-1)^{|p_a|}q_b$ to
$\d(p_b),\d(q_a)$ respectively.
Write
\begin{align*}
\d(p_b) &= -\epsilon p_a + u \\
\d(q_a) &= \epsilon (-1)^{|p_a|} q_b + v
\end{align*}
for some $u,v\in\alg$.
Let $\phi_0$ be the algebra map on $\alg$ defined by
\begin{align*}
\phi_0(q_b) &= \epsilon (-1)^{|p_a|} q_b + v \\
\phi_0(p_a) &= -\epsilon p_a + u \\
\phi_0(s) &= s \qquad\qquad\qquad \text{for $s\neq q_b,p_a$}.
\end{align*}

\begin{lemma}
$\phi_0$ is a basis change on $\alg$.
\end{lemma}

\begin{proof}
We have $\phi_0 = \phi_b \circ \phi_a$, where $\phi_b$ is the
elementary automorphism supported at $q_b$ that
sends $q_b$ to $\epsilon (-1)^{|p_a|} q_b +v$, and
$\phi_a$ is the elementary automorphism supported at $p_a$ that sends
$p_a$ to $-\epsilon p_a + \phi_b^{-1}(u)$. Now by Lemma~\ref{lem:stokes}, any term
in $v$ either involves a $p$ or only involves $t^{\pm 1}$ and $q$'s of
smaller height than $q_b$, and so $\phi_b$ is an elementary
automorphism of $\alg$. Also by Lemma~\ref{lem:stokes}, any term in $u$
either involves two $p$'s or only involves $t^{\pm 1}$, $q$'s, and
$p$'s of greater height than $p_a$; it follows that the only terms in
$\phi_b^{-1}(u)$ not containing at least two $p$'s do not involve
$p_a$, and so $\phi_a$ is an elementary automorphism as well.
\end{proof}


So far, the argument given here is a straightforward extension of Chekanov's proof of Reidemeister-II invariance in \cite[section~8.4]{Ch02}, which then hinges on the following two points \cite[Lemma~8.2]{Ch02}:
\begin{enumerate}
\item
$d'$ and $\phi_0^{-1} \circ d \circ \phi_0$ agree on $q$ generators whose height is at most the height of $q_b$;
\label{item:ch1}
\item
on all $q$ generators, $d'$ and $\phi_0^{-1} \circ d \circ \phi_0$ agree modulo terms involving either $q_a$ or $q_b$.
\label{item:ch2}
\end{enumerate}
Chekanov uses these two properties to bootstrap $\phi_0$ up to an automorphism that intertwines $d'$ and $d$ for all $q$ generators.

The proof in the current circumstance is complicated by the fact that (\ref{item:ch1}) no longer holds, due to the possible presence of $p$'s in the differentials of any $q$ generator. Nevertheless, an analogue of (\ref{item:ch2}) still holds and is presented as Lemma~\ref{lem:ReidII} below. We will use this, along with a property of the differential that we call ``ordered''-ness (Definition~\ref{def:ordered}) that loosely generalizes property (\ref{item:ch1}), to perform a bootstrapping argument similar to Chekanov's.

In order to prove the analogue of (\ref{item:ch2}), we need to establish a lemma that extends the central argument in Chekanov's proof of (\ref{item:ch2}) \cite[section~8.5]{Ch02}.
Define a graded algebra map
$\psi\co\alg\to\alg$ by $\psi(p_a) = \epsilon u$, $\psi(q_b) =
-\epsilon (-1)^{|p_a|} v$, $\psi(q_a) = \psi(p_b) = 0$, and $\psi(s) =
s$ for all other generators of $\alg$.

\begin{lemma}
Let $s$ be a generator of $\alg$ besides $p_a,q_a,p_b,q_b,t^{\pm
  1}$, and let $\psi^n$ denote the $n$-th iterate of $\psi$. In $\alg$, the limit $\lim_{n\to\infty} \psi^n \dsft(s)$ exists
and equals $\dsft'(s)$.
\label{lem:limit}
\end{lemma}

\begin{proof}
We first show that the limit exists. We can assume that $h(p_a) -
h(p_b) > 0$ is arbitrarily small, since by Lemma~\ref{lem:stokes} this
is the area of the bigon determined by $p_a$ and $q_b$. By
Lemma~\ref{lem:stokes}
again, this implies that any term in $v$ or $u$ involving $q_b$ must
be $O(p)$, while any term in $u$ involving $p_a$ must be $O(p^2)$. It
follows that for any $m$, the portion of $\psi^n \dsft(s)$ involving
at most $m$ $p$'s stabilizes as $n\to\infty$, and thus the limit
exists.

It remains to show that
\begin{equation}
\lim_{n\to\infty} \psi^n \dsft(s) = \dsft'(s).
\label{eq:limit}
\end{equation}
First assume for simplicity that $u,v\in\alg'$, i.e., that $u,v$
do not involve $p_a,q_a,p_b,q_b$. In this case, the left hand side of (\ref{eq:limit}) is $\psi \dsft(s)$.

\begin{figure}
\centerline{
\includegraphics[height=1.6in]{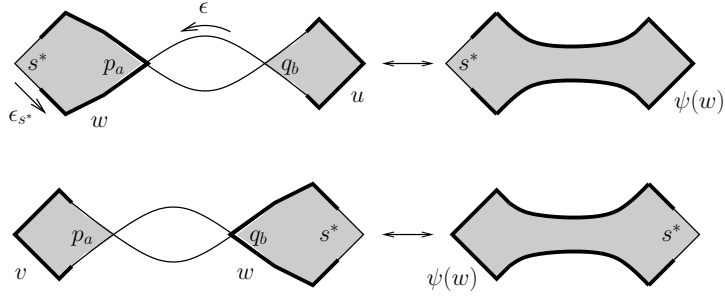}
}
\caption{
Top: gluing a term $w$ in $\dsft(s)$ to $u$ to obtain a term (or sum of terms)
$\psi(w)$ in $\dsft'(s)$; bottom: gluing a term $w$ in $\dsft(s)$ to $v$ to
obtain a term $\psi(w)$ in $\dsft'(s)$. Two signs
$\epsilon,\epsilon_{s^*}$ are also shown.
}
\label{fig:ReidIIglue}
\end{figure}

Any term in $\dsft(s)$ not involving any of $p_a,q_a,p_b,q_b$
corresponds to a disk preserved by the Reidemeister II move and thus
appears in $\dsft'(s)$ as well. Any term in $\dsft(s)$ involving
either $q_a$ or $p_b$ is killed by
$\psi$. The remaining terms in $\dsft(s)$ involve $p_a$ or $q_b$ but
not $q_a$ or $p_b$; call such a term $w$. Then $\psi(w)$ appears in
$\dsft'(s)$: at each $p_a$ or $q_b$ corner of $w$, glue $u$ or $v$,
respectively; this gives the disks in $\dsft'(s)$ passing through the
neck in $\pi_{xy}(\Lambda')$. See Figure~\ref{fig:ReidIIglue}. This
proves (\ref{eq:limit}) mod $2$ in this case.

In fact, the signs in the definition of $\psi$ work out so that
$\dsft'(s) = \psi \dsft(s)$ over $\Z$. Consider for instance gluing
$u$ to $p_a$, where we abuse notation and
use $u$ to denote a particular term in $u$. Let
$\epsilon_w,\epsilon_u$ be the product of the orientation signs over
the corners in $w,u$, and let $\epsilon_{s^*}$ be the sign shown in
Figure~\ref{fig:ReidIIglue}. Then $w$ appears in $\dsft(s)$ with sign
$\epsilon_{s^*} \epsilon_w$, while $u$ has sign $\epsilon
\epsilon_u$. On the other hand, $\psi(w)$ appears in $\dsft'(s)$ with
sign $\epsilon_{s^*} \epsilon_w \epsilon_u$ since the orientation
signs for the corners of $w,u$ at $p_a,q_b$, respectively, are
equal. This agrees with the fact that $\psi$ sends $p_a$ to $\epsilon
u$. We can similarly check the sign in $\psi(q_b) = -\epsilon
(-1)^{|p_a|}v$. This completes the proof of (\ref{eq:limit}) when
$u,v\in\alg'$.

In general, even if $u,v$ involve $p_a,q_a,p_b,q_b$, the same argument
shows (\ref{eq:limit}).
To get $\dsft'(s)$, one starts with $\dsft(s)$ and
keeps replacing any appearances of $p_a$, $q_b$, $q_a$, $p_b$ by
$\epsilon u$, $-\epsilon (-1)^{|p_a|} v$, $0$, $0$, respectively; but
this is precisely what $\lim_{n\to\infty} \psi^n \dsft(s)$ gives.
\end{proof}





\begin{definition}
If $x,y\in\alg$, write $x\equiv y \pmod{a,b}$ if $x-y$ only includes
terms involving at least one of $p_a,q_a,p_b,q_b$ (in the notation
from Section~\ref{ssec:alg}, $\pi(x-y) = 0$). If $f,g$ are two maps
from $\alg$ to $\alg$, write $f\equiv g \pmod{a,b}$ if $f(x) \equiv g(x)
\pmod{a,b}$ for all $x\in\alg$.
\end{definition}

\begin{lemma}
On $\alg$, $\d' \equiv \phi_0^{-1} \circ \d \circ \phi_0 \pmod{a,b}$.
\label{lem:ReidII}
\end{lemma}

\begin{proof}
It suffices to show that $\d'(s) \equiv \phi_0^{-1} \circ \d \circ
\phi_0(s) \pmod{a,b}$ for all generators $s$ of $\alg$.
We have
\[
\phi_0^{-1} \d \phi_0(p_a) 
= \phi_0^{-1} \d^2 (p_b) = \phi_0^{-1} [F_{\d},p_b]
\equiv 0 \equiv [F_{\d'},p_b] = \d'(p_a)
\pmod{a,b}
\]
and similarly $\phi_0^{-1}\d\phi_0(s) \equiv \d'(s) \equiv 0
\pmod{a,b}$ for $s=q_a,p_b,q_b$.
Also, $\d'(t)-\d(t)$ consists of four terms, one each involving
$p_aq_a$, $q_ap_a$, $p_bq_b$, $q_bp_b$, and it follows easily that
$\phi_0^{-1}\d\phi_0(t^{\pm 1}) \equiv \d'(t^{\pm 1}) \pmod{a,b}$.

Now assume that $s\neq p_a,q_a,p_b,q_b,t^{\pm 1}$; we want to show
that $\d's \equiv \phi_0^{-1} \d s \pmod{a,b}$. Since $\dstr's \equiv
\phi_0^{-1} \dstr s \pmod{a,b}$ as before, it suffices to show that
$\dsft'(s) \equiv \phi_0^{-1} \dsft(s) \pmod{a,b}.$ By
Lemma~\ref{lem:limit}, this follows from establishing that
$\lim_{n\to\infty} \psi^n \equiv \phi_0^{-1} \pmod{a,b}$.

In fact, we
claim that on $\alg$, $\lim_{n\to\infty} \psi^n = \pi \circ \phi_0^{-1}$, or
equivalently $\lim_{n\to\infty} \psi^n \circ \phi_0 = \pi$, where $\pi$ is the projection map from $\alg$ to $\alg'$ as
usual. Indeed, since both sides are algebra maps, it suffices to check
that $\lim_{n\to\infty} \psi^n \circ \phi_0(s) = \pi(s)$ for all
generators $s$ of $\alg$. This holds trivially unless $s = q_b$ or
$p_a$; it also holds for $s=q_b$ since
\[
\epsilon(-1)^{|p_a|} \lim_{n\to\infty} \psi^n \phi_0(q_b) =
\lim_{n\to\infty} \psi^n (q_b+\epsilon (-1)^{|p_a|}v) =
\lim_{n\to\infty} (\psi^n q_b - \psi^{n+1} q_b) = 0,
\]
and similarly for $s=p_a$.
\end{proof}

Write $\d_0 = \phi_0^{-1}\circ \d\circ \phi_0$ on $\alg$. We claim that
$\d_0$ is related to $d'$ by a basis change; this will imply that $\d$
is related to $\d'$ by a basis change, which will complete the
invariance proof for the restricted Reidemeister II move.

\begin{definition}
We say that a derivation $\d$ on $\alg$ is \textit{ordered} if
\begin{align*}
\d(q_j) &= (\text{function of $t^{\pm 1},q_a,q_b,q_1,\dots,q_{j-1}$}) + O(p) \\
\d(p_j) &= (\text{function of $t^{\pm
    1},p_a,q_a,p_b,q_b,q_1,\dots,q_n,p_{j+1},\dots,p_n$}) + O(p^2) \\
\d(t) &= O(p).
\end{align*}
\label{def:ordered}
\end{definition}

\noindent
Order the crossings of $\Lambda$ (or $\Lambda'$) in such a way that $h(q_1) \leq
h(q_2) \leq \dots \leq h(q_n)$; then by Lemma~\ref{lem:stokes}, $\d$ is
automatically ordered.

\begin{lemma}
$d_0 = \phi_0^{-1} \circ \d \circ \phi_0$ is ordered.
\label{lem:d0ordered}
\end{lemma}

\begin{proof}
Since $\phi_0$ preserves the $p$ filtration, it is clear that $\d_0(t)
= O(p)$. Next consider $\d_0(q_j) = \phi_0^{-1}(\d(q_j))$, and note
that $\phi_0^{-1}$ fixes all generators of $\alg$ except for $q_b$ and
$p_a$. We wish to show that the order $p^0$ terms (that is, the terms
that are not $O(p^1)$) in $\d_0(q_j)$ do
not involve $q_j,\dots,q_n$.
Since $\d$ is ordered and $\phi_0$ fixes all words that do not
involve $q_b$ or $p_a$, it suffices to show that if $w$ is a word in
$\d(q_j)$ involving $q_b$, then the order $p^0$ terms in
$\phi_0^{-1}(w)$ do not involve $q_j,\dots,q_n$.

We may assume that $|h(q_j)| > |h(q_b)|$, since otherwise any term in
$\d(q_j)$ involving $q_b$ must be $O(p)$ by
Lemma~\ref{lem:stokes}. As in the proof of Lemma~\ref{lem:ReidII}, we
may also assume that $|h(q_a)| - |h(q_b)| > 0$ is arbitrarily small
(more precisely, that no $|h(q_j)|$ lies in the interval $[|h(q_b)|,|h(q_a)|]$).
Then since $\d(q_a) = q_b + v$, any term in $v$ involving
$q_j,\dots,q_n$ must be $O(p)$; otherwise, by Lemma~\ref{lem:ReidII},
$|h(q_a)| > |h(q_i)|$ for some $i \geq j$, and $|h(q_i)| \geq
|h(q_j)|$. But then to order $p^0$, $v$ and hence $\phi_0^{-1}(w)$
does not involve $q_j,\dots,q_n$.

Finally, we need to prove that the order $p^1$ terms (that is, the
terms that are not $O(p^2)$) in
$\d_0(p_j) = \phi_0^{-1}(\d(p_j))$ do not involve
$p_1,\dots,p_{j-1}$. As before, since $\d$ is ordered, the only place
$p_1,\dots,p_{j-1}$ can appear in the order $p^1$ terms in $\d_0(p_j)$
is in contributions from $\phi_0^{-1}(q_b)$ or $\phi_0^{-1}(p_a)$. Any
contribution from $\phi_0^{-1}(q_b)$ is $O(p^2)$, since it contains
one of $p_1,\dots,p_{j-1}$ along with some
other $p$ (from the fact that $\d(p_j) = O(p)$). Now if $h(p_j) >
h(p_a)$, then all terms in $\d(p_j)$ involving $p_a$ are $O(p^2)$,
while if $h(p_j) < h(p_a)$, then any term in $u$ involving
$p_1,\dots,p_{j-1}$ is $O(p^2)$ by Lemma~\ref{lem:stokes} again. In
either case, we conclude that the order $p^1$ terms in
$\phi_0^{-1}(\d(p_j))$ cannot involve $p_1,\dots,p_{j-1}$.
\end{proof}

\begin{figure}[b]
\centerline{
\includegraphics[height=2in]{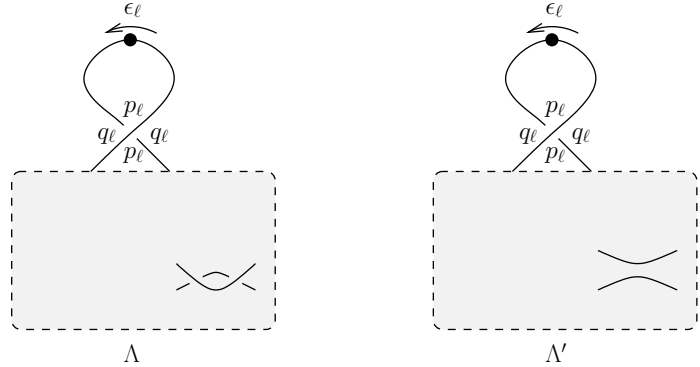}
}
\caption{
Labels for a restricted Reidemeister II move.
}
\label{fig:ReidIIreslabel}
\end{figure}

To find the basis change relating $d_0$ to $d'$, we need to use the fact that the Reidemeister II move is restricted. Let
$p_\ell,q_\ell$ label the crossing whose loop lies outside the rest of the
diagram and does not interact with the move, and choose $\bullet$ to
lie on this loop and $\ast$ not to lie on this loop.

\begin{lemma}
For a restricted Reidemeister II move, we have: \label{lem:restricted}
\begin{itemize}
\item
$\d(p_\ell) = \d'(p_\ell) = 0$;
\item
$F_{\d} = F_{\d'} = \epsilon_\ell p_\ell$,
where $\epsilon_\ell$ is the sign shown in Figure~\ref{fig:ReidIIreslabel};
\item
$d'(s) = d_0(s)$ if $s = p_a,q_a,p_b,$ or $q_b$.
\end{itemize}
\end{lemma}

\begin{proof}
The only term in the Hamiltonians for $\Lambda$ or $\Lambda'$
involving $q_\ell$ is $p_\ell$. The first two properties then follow from
the definitions of $\d,\d'$ and Proposition~\ref{prop:d2}. The third
property can be proven by trivially modifying the first paragraph of the
proof of Lemma~\ref{lem:ReidII}, where we now use the fact that
$\phi_0^{-1}(F_\d) = F_\d = F_{\d'}$.
\end{proof}

For ease of notation, define $q_{-n},\dots,q_0$ by $q_{-j} = p_j$ for $1 \leq j\leq n$, $q_0
= t$. Then the condition for a derivation $\d$ to be ordered is that
for all $j$ with $-n\leq j\leq n$, any term in $\d(q_j)$ involving one of
$q_j,q_{j+1},\ldots,q_n$ also involves another $p$. Here a word $w$ involving $q_i$ involves another $p$ if
$i\geq 0$ and $w=O(p)$, or if $i<0$ and $w=O(p^2)$.

Starting with $\d_{0,1} := \d_0$, we will inductively
define a sequence of differentials
$\d_{j,1},\phi_{j,1}$ for $1\leq j\leq n$ and
$\d_{j,k},\phi_{j,k}$ for $-n\leq j\leq n$ and $k\geq 2$.

\begin{claim}
We can inductively construct $\d_{j,k},\phi_{j,k}$ to satisfy the following properties:

\renewcommand{\theenumi}{\roman{enumi}}
\begin{enumerate}
\item
we have
\begin{align*}
\d_{j,1} &= \phi_{j,1}^{-1} \circ \d_{{j-1},1} \circ
\phi_{j,1} && 1\leq j\leq n \\
\d_{j,k} &= \phi_{j,k}^{-1} \circ \d_{{j-1},k} \circ
\phi_{j,k} && 1-n \leq j\leq n,~~~k\geq 2 \\
\d_{{-n},k} &= \phi_{{-n},k}^{-1} \circ \d_{n,k-1} \circ
\phi_{{-n},k} && k\geq 2;
\end{align*}
\item \label{it:5}
$\d_{j,k}$ is ordered;
\item
$\d_{j,k}(s) = \d'(s)$ for $s=p_a,q_a,p_b,q_b$, $\d_{j,k}(p_\ell) = 0$, and $F_{d_{j,k}} = \epsilon_\ell p_\ell$; \label{it:p}
\item
$\phi_{j,k} \equiv \Id \pmod{p^{k-1}}$ for all $j,k$;
\label{it:2}
\item
$\d_{j,k} \equiv \d' \pmod{a,b}$ and $\d_{j,k} \equiv \d'
\pmod{p^{k-1}}$ for all $j,k$;
\label{it:3}
\item
$\d_{j,k}(q_i) \equiv \d'(q_i) \pmod{p^k}$ for $i\leq j$, for all $j,k$.
\label{it:4}
\end{enumerate}
\end{claim}

Note that $\d_{0,1}$ satisfies (\ref{it:5}) by
Lemma~\ref{lem:d0ordered}, (\ref{it:p}) by Lemma~\ref{lem:restricted},
(\ref{it:3}) by Lemma~\ref{lem:ReidII}, and
(\ref{it:4}) because $\d_{0,1}(q_i) \equiv \d'(q_i) \equiv 0
\pmod{p^1}$ for $i\leq 0$.

The following diagram summarizes the inductive order of the construction.
Given $\d_{j-1,k}$ for $j\geq 1-n$, we construct
$\phi_{j,k},\d_{j,k}$; given $\d_{n,k-1}$, we construct $\phi_{-n,k},\d_{-n,k}$.
Each differential agrees with $\d'$ to the specified order in $p$ when
applied to a generator corresponding to its column or any column to
its left, and to order one less in $p$ when applied to any generator
corresponding to a column to its right.
\[
\xymatrix@C=8pt@R=0pt@W=0pt@H=0pt{
&& *=0{} \ar@{-}[dddddd] \\
&&&  & *{p_n} & p_{n-1} & \cdots & p_1 & t & q_1 & q_2 &
\dots & q_n \\
*=0{} \ar@{-}[rrrrrrrrrrrrr] &&&&&&&&&&&&& \\
& \rule[-8pt]{0pt}{18pt} O(p^1)
&&&&&&&
d_0=d_{0,1} \ar[r] & d_{1,1} \ar[r] & d_{2,1} \ar[r]
& \cdots \ar[r] & d_{n,1} \ar `/6pt[d] `[l] `[dlllllllll] `[dllllllll]
 [dllllllll] \\
& \rule[-8pt]{0pt}{18pt} O(p^2)  &&& *{d_{-n,2}} \ar[r] & d_{1-n,2}
\ar[r] & \cdots \ar[r] &
d_{-1,2} \ar[r] &
d_{0,2} \ar[r] & d_{1,2} \ar[r] & d_{2,2} \ar[r]
& \cdots \ar[r] & d_{n,2} \ar `/6pt[d] `[l] `[dlllllllll] `[dllllllll]
 [dllllllll] \\
& O(p^3) \rule[-8pt]{0pt}{18pt} &&& *{d_{-n,3}} \ar[r] & \cdots &&&&&&& \\
&&
}
\]


\begin{proof}[Proof of Claim]
Suppose that $\d_{j-1,k}$ satisfies (\ref{it:5}) through (\ref{it:4})
for some $j\leq n$. Define the elementary automorphism $\phi_{j,k}$ on
$\alg$, supported at $q_j$, by
\[
\phi_{j,k}(q_j) = q_j - H \d_{j-1,k}(q_j),
\]
where $H\co\alg\to\alg$ is the operator defined in
Section~\ref{ssec:alg}. (Note that the notation has changed slightly
from Section~\ref{ssec:alg}: $\alg$ there is $\alg'$ here, and
$S_i\alg$ there is $\alg$ here.) Observe that $\phi_{j,k}$ is
elementary because $\d_{j-1,k}$ is ordered by assumption.

Define $\d_{j,k} = \phi_{j,k}^{-1} \circ \d_{j-1,k} \circ
\phi_{j,k}$. We wish to show that $\d_{j,k},\phi_{j,k}$ satisfy
(\ref{it:5}) through (\ref{it:4}).
(A corresponding construction produces $\d_{-n,k}$ from $\d_{n,k-1}$; here
$\phi_{-n,k}$ is supported at $q_{-n}=p_n$ and $\phi_{-n,k}(q_{-n}) =
q_{-n} - H \d_{n,k-1}(q_{-n})$. The proof that $\d_{-n,k},\phi_{-n,k}$
satisfy (\ref{it:5}) through (\ref{it:4}) is entirely
similar and will be omitted here.)

We first check (\ref{it:5}): for all $i$, any term in $\d_{j,k}(q_i)$
involving $q_i,\dots,q_n$ must include another $p$ as well.
If $i < j$, then since $\d_{j-1,k}$ is
ordered, any term in $\d_{j-1,k}(q_i)$ involving $q_j$ must include
another $p$, and the condition holds. If $i\geq j$, then since
any term in $H\d_{j-1,k}(q_j)$ involving $q_j,\dots,q_n$ must include
another $p$, it follows that any term in $\phi_{j,k}^{-1}(q_j)$
involving $q_j,\dots,q_n$ must also include another $p$, and the
condition holds here as well. This demonstrates (\ref{it:5}) for $\d_{j,k}$.

As for the other conditions, note that (\ref{it:p}) holds for
$\d_{j,k}$ since $\phi_{j,k}$ preserves $p_a,q_a,p_b,q_b,p_\ell$ by
construction (for $p_\ell$, use
the induction hypothesis $\d_{j-1,k}(p_\ell) = 0$).
Since $H$ preserves the $p$ filtration and
$\d_{j-1,k}(q_j) \equiv \d'(q_j) \pmod{p^{k-1}}$, we have
\[
H \d_{j-1,k}(q_j) \equiv H \d'(q_j) = 0 \pmod{p^{k-1}}
\]
and thus $\phi_{j,k} \equiv \Id \pmod{p^{k-1}}$, as required by
(\ref{it:2}). It follows that (\ref{it:3}) holds for $\d_{j,k}$ since
it holds for $\d_{j-1,k}$. As for (\ref{it:4}), if $i<j$, then
\[
\d_{j,k}(q_i) = \phi_{j,k}^{-1}\d_{j-1,k}\phi_{j,k}(q_i)
= \phi_{j,k}^{-1}\d_{j-1,k}(q_i) \equiv
\phi_{j,k}^{-1} \d'(q_i) \equiv \d'(q_i) \pmod{p^k},
\]
where the second-to-last equality holds since $\d_{j-1,k}(q_i) \equiv
\d'(q_i) \pmod{p^k}$ by induction assumption, and the final equality
holds because the only terms in $\d'(q_i)$ that can involve $q_j$ must
also involve another $p$.

Thus to complete the induction step,
we need to establish that $\d_{j,k}(q_j) \equiv \d'(q_j)
\pmod{p^k}$. Since the only terms $\d'(q_j)$
involving $q_j$ itself must also involve another $p$, we have
$\phi_{j,k}^{-1} \d'(q_j) \equiv \d'(q_j) \pmod{p^k}$. By the construction
of $\d_{j,k}$ in terms of $\d_{j,k-1}$, it now suffices to show that
\begin{equation}
\d_{j-1,k} \phi_{j,k}(q_j) \equiv \d'(q_j) \pmod{p^k}.
\label{eq:RIIinduct}
\end{equation}

For ease of notation, write $\tilde{\d} = \d_{j-1,k}$.
Recall the map $H\co \alg\to\alg$ from
Lemma~\ref{lem:htpy}; by Lemma~\ref{lem:htpy}, we have $H\d'+\d'H = \Id
- \iota\circ\pi$, where $\iota\circ\pi$ is the map on $\alg$ that
projects away any term involving $p_a,q_a,p_b,q_b$. It follows that
\[
H\d'\tilde{\d}q_j + \d'H\tilde{\d}q_j = \tilde{\d}q_j -
\iota\pi\tilde{\d}q_j = \tilde{\d}q_j - \d'q_j
\]
where the last equality holds by property (\ref{it:3}). Thus
\begin{align*}
\tilde{\d}\phi_{j,k}q_j &= \tilde{\d}(q_j-H\tilde{\d}q_j) \\
&= \d'q_j - (\tilde{\d}-\d')H\tilde{\d}q_j - H(\tilde{\d}-\d')\tilde{\d}q_j +
H[F_{\tilde{\d}},q_j].
\end{align*}

Now
from
properties (\ref{it:5}), (\ref{it:3}), and (\ref{it:4}) for
$\tilde{\d}$, it follows that
\[
(\tilde{\d}-\d')\tilde{\d}q_j \equiv (\tilde{\d}-\d')H\tilde{\d}q_j
\equiv 0 \pmod{p^k}.
\]
Also $F_{\tilde{\d}} = F_{\d'} \in\alg'$ so $H[F_{\tilde{\d}},q_j] =
0$. The desired equation (\ref{eq:RIIinduct}) now follows.

This completes the induction step and the proof of the claim.
\end{proof}

To finish the proof of invariance
under restricted Reidemeister II, we note that since
$\d_0 = \phi_0^{-1} \circ \d
\circ \phi_0$, we can write $\d' = \phi^{-1} \circ \d \circ \phi$,
where
\[
\phi = \phi_0 \phi_{1,1} \phi_{2,1} \cdots \phi_{n,1} \phi_{-n,2}
\cdots \phi_{n,2} \phi_{-n,3} \cdots.
\]
This is an infinite composition, but for any $k$, all but finitely many
terms in this composition are congruent to the identity
$\pmod{p^k}$. For $\phi$ to be a change of basis, we need rewrite this
as a composition of finitely many elementary automorphisms. The
following result thus completes the invariance proof.

\begin{lemma}
If, for $-n\leq j\leq n$ and $k\geq 2$, $\phi_{j,k}$ are elementary
automorphisms with $\phi_{j,k}$ supported at $q_j$ and $\phi_{j,k}
\equiv \Id \pmod{p^{k-1}}$ for all $j,k$, then we can write the
infinite composition
\[
\phi_{-n,2} \cdots \phi_{n,2} \phi_{-n,3} \cdots
\]
as a finite composition $\phi_{(-n)} \cdots \phi_{(n)}$, where
$\phi_{(j)}$ is an elementary automorphism supported at $q_j$.
\end{lemma}

\begin{proof}
Consider two elementary automorphisms $\phi_1,\phi_2$ on
$\alg$ supported at two different generators $s_1,s_2$ with
$\phi_1(s_1) = s_1+v_1$, $\phi_2(s_2)=s_2+v_2$, and $v_1,v_2=O(p^2)$.
Then the composition $\phi_2\circ\phi_1$ can also be
written as $\phi_1'\circ\phi_2'$, where $\phi_1',\phi_2'$ are
elementary automorphisms supported at $s_1,s_2$ with
\begin{align*}
\phi_1'(s_1) &= s_1 + \phi_2(v_1) \\
\phi_2'(s_2) &= s_2 + (\phi_1')^{-1}(v_2);
\end{align*}
that is, we can rewrite a composition of elementary automorphisms
supported at $s_1$ and $s_2$ as a similar composition with the roles
of $s_1$ and $s_2$ reversed. In addition, if $\phi_1 \equiv \Id
\pmod{p^k}$ for some $k$, then $\phi_1' \equiv \Id
\pmod{p^k}$ and $\phi_2' \equiv \phi_2 \pmod{p^k}$.

Through this
trick, we can rewrite a ``partial convergent''
\[
(\phi_{-n,2} \cdots \phi_{n,2}) \cdots (\phi_{-n,k} \cdots \phi_{n,k})
\]
of the infinite composition as $\phi_{(-n,k)} \cdots \phi_{(n,k)}$, where
$\phi_{(j,k)}$ is an elementary automorphism supported at $q_j$ with
$\phi_{(j,k)}(q_j) = q_j + v_{j,k}$ for some $v_{j,k} = O(p^2)$. (To
this end, note that the composition of two elementary automorphisms
supported at $q_j$ is another.) It is easy to see that $v_{j,k} \equiv
v_{j,k+1} \pmod{p^k}$ for all $k$, and thus that $v_{j,k}$ has a limit
in $\alg$ as $k\to\infty$. Defining $\phi_{(j)}$ for $-n\leq j\leq n$
to be the elementary
automorphism supported at $q_j$ with $\phi_{(j)}(q_j) = q_j +
\lim_{k\to\infty} v_{j,k}$ then completes the proof of the lemma.
\end{proof}

\secspace

\section*{Appendix A: Orientation Signs}

\setcounter{figure}{0}
\renewcommand{\thefigure}{A.\arabic{figure}}
\setcounter{theorem}{0}
\renewcommand{\thetheorem}{A.\arabic{theorem}}

To define the Hamiltonian $h$ and the SFT differential $\dsft$ over
$\Z$, we chose particular orientation signs as shown in
Figure~\ref{fig:pqsigns}; see also Remark~\ref{rem:LCH}.
These are not the only possible orientation
signs leading to a viable LSFT algebra. Here we find all possible
combinatorial choices for orientation signs and show that they are all
equivalent
under basis change. As a corollary, we obtain a refinement of a result
in \cite{EES05b}. There two sign choices for Legendrian
contact homology in $\R^3$ are given, one of which recovers the signs
from \cite{ENS} and one of which appears to be different; we show that
the two choices are in fact equivalent.

\begin{figure}[b]
\centerline{
\includegraphics[height=0.9in]{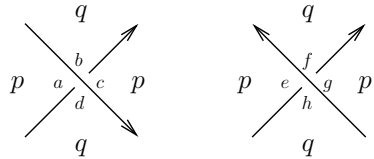}
}
\caption{
Possible orientation signs for corners. Each of $a,b,c,d,e,f,g,h$ is
$\pm 1$. The left figure is a positive crossing in the usual
knot-theoretic sense, with $q$ even and $p$ odd; the right figure is a
negative crossing, with $q$ odd and $p$ even.
}
\label{fig:pqsigns2}
\end{figure}

The most general set of orientation signs has eight degrees of freedom
$a,b,c,d,e,f,g,h\in\{\pm 1\}$,
one for each quadrant of a positive and a negative crossing; see
Figure~\ref{fig:pqsigns2}. Note that in the formulation from
Section~\ref{sec:string}, the orientation signs only figure in the
definition of $h$.

\begin{figure}[b]
\centerline{
\includegraphics[height=0.9in]{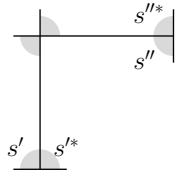}
}
\caption{
Signs from an obtuse disk (cf.\ Figure~\ref{fig:hhsigns}).}
\label{fig:hhsigns2}
\end{figure}

To give rise to an LSFT algebra structure,
signs must be chosen so that an identity like $\frac{1}{2} \sftb{h,h} +
\delta h = 0$ (Proposition~\ref{prop:QME}) holds. In particular, the two
terms in $\sftb{h,h}$ arising from an obtuse disk must cancel. From
the proof of Proposition~\ref{prop:QME}, we see that we must have
$\epsilon \sftb{s',s'^*} \sftb{s'',s''^*} = 1$ for all configurations
of the form depicted in Figure~\ref{fig:hhsigns2}, where $\epsilon$ is
the product of the orientation signs over the six shaded corners. One
readily deduces that we must have $ab=gh=-ad=-eh$ and $cd=ef=-bc=-fg$,
whence $d=-b$, $g=-e$, $f = -bc/e$, and $h = ab/e$. This reduces us to
four degrees of freedom $a,b,c,e$.

We can get rid of three further degrees of freedom as
follows. Replacing $(a,b,c,e)$ by $(-a,-b,-c,e)$ has the effect in the
LSFT algebra of replacing all $p_j,q_j$ corresponding to positive
crossings by $-p_j,-q_j$, thus only modifying the LSFT algebra by a
basis change. Similarly, replacing $(a,b,c,e)$ by $(a,b,c,-e)$ just
replaces all $p_j,q_j$ corresponding to negative crossings by
$-p_j,-q_j$. Furthermore, replacing $(a,b,c,e)$ by $(-a,b,-c,-e)$ has
no effect on $h$ or the LSFT algebra, since this simply changes each
term $w$ in $h$ by $(-1)^{o(w)}$, where $o(w)$ is the number of
odd-degree generators in $w$ and is always even since $|h|=-2$.

\begin{figure}
\centerline{
\includegraphics[height=1.1in]{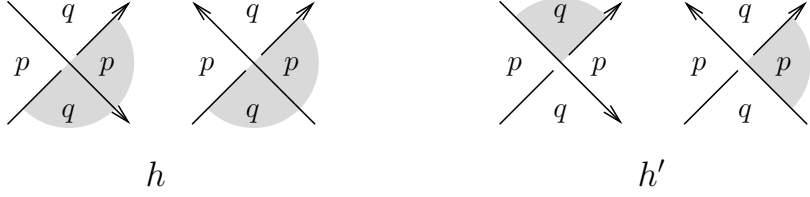}
}
\caption{
The orientation signs used in this paper (left) and another choice
agreeing with \cite{ENS} (right), resulting in Hamiltonians $h$ and $h'$.}
\label{fig:pqsigns3}
\end{figure}

Eliminating these three degrees of freedom, we are left with two
possibly different equivalence classes of orientation signs,
represented by $(a,b,c,e) = (1,1,-1,1)$ and $(1,-1,1,1)$ and depicted in
Figure~\ref{fig:pqsigns3}. These orientation signs yield
two Hamiltonians $h,h'\in\alg$ that agree mod $2$. The
first is the Hamiltonian used in this paper, satisfying the quantum
master equation $\frac{1}{2} \sftb{h,h} + \delta h = 0$. It can
readily be shown that the second satisfies the equation $\frac{1}{2}
\sftb{h',h'} - \delta h' = 0$. We can then define a derivation $\d' =
\sftb{h,\cdot} - \delta x$, and $(\alg,\d')$ is an LSFT algebra.

Each of $(\alg,\d)$ and $(\alg,\d')$ induces a choice of signs for
the differential on Legendrian contact homology
$\F^0\!\alg/\F^1\!\alg$. In \cite[Theorem~4.32]{EES05b}, two sign rules
for Legendrian contact homology are given, essentially corresponding
to the two orientations on $\C$; these correspond to our
$(a,b,c,e) = (1,-1,1,1)$ and $(-1,-1,1,1)$ and hence to $(\alg,\d')$
and $(\alg,\d)$, respectively. The first sign rule in \cite{EES05b}
also agrees with the signs given in \cite{ENS}.\footnote{To translate
  between our
signs and the signs for contact homology in \cite{EES05b,ENS},
we must incorporate the sign
$\epsilon(f;s_1)$ (cf.\ Section~\ref{ssec:comb}) measuring the
orientation of the disk after the $p$ puncture. This has the effect of
negating the signs for the corners marked $a$ and $g$ in
Figure~\ref{fig:pqsigns2}.}

At the time of the writing of \cite{EES05b}, it was not known whether
the two sign rules led to different contact homology differential
graded algebras. In fact, we
shall see that they are equivalent. This follows from the
corresponding result for the LSFT algebras.

\begin{proposition}
\label{prop:twosigns}
Let $\Lambda$ be a Legendrian knot in standard contact $\R^3$.
The LSFT algebra $(\alg,\d)$ for $\Lambda$ is related by a basis
change to the LSFT
algebra obtained from $(\alg,\d')$ by conjugation with the
involution $t^{\pm 1} \mapsto (-1)^{r(\Lambda)}t^{\pm 1}$, where
$r(\Lambda)$ is the
rotation number of $\Lambda$.
\end{proposition}

\begin{corollary}
\label{cor:twosigns}
The DGAs given by the two sign rules in
\cite[Theorem~4.32]{EES05b}
are tamely isomorphic if we first replace $t$ in one of the DGAs by
$(-1)^{r(\Lambda)}t$. Here the tame isomorphism can be chosen to
extend the identity map on the base ring $\Z[t,t^{-1}]$.
\end{corollary}

\begin{proof}[Proof of Proposition~\ref{prop:twosigns}]
Define an involution $\phi_1$ on $\alg$ that negates $p_j,q_j$ for all
$j$ such that $p_j,q_j$ corresponds to a positive crossing, and let
$h''=\phi_1(h')$. The orientation signs defining $h''$ are the same as
those defining $h'$, except that all signs for positive crossings are
reversed. An examination of Figure~\ref{fig:pqsigns3} then shows that
the orientation signs between $h$ and $h''$ only differ at corners
where the knot is oriented into the crossing on both sides of the
corner.

Define another involution $\phi_2$ on $\alg$ as follows:
\[
\phi_2(s) = \begin{cases}
s & \text{if $|s| \equiv 0,1 \pmod 4$} \\
-s & \text{if $|s| \equiv 2,3 \pmod 4$} \\
(-1)^{r(\Lambda)} t^{\pm 1} & \text{if $s=t^{\pm 1}$}.
\end{cases}
\]
(Note that the third line is superfluous but has been included for
clarity.)
We claim that $h'' = -\phi_2(h)$. Indeed, the difference in signs
between the appearances of a word $w$ in $h$ and in $h''$ is
$(-1)^{o(w)/2}$, where $o(w)$ is the number of odd $s$'s appearing in
$w$, so that $o(w)/2$ is the number of corners where the sign changes between $h$ and $h''$, cf.\ Lemma~\ref{lem:degsign}. Suppose that the word $w$ contains $m_j$ generators (counting multiplicity) whose degree is $j \pmod 4$ for $j=0,1,2,3$. Since $|w| = -2$, we have $m_1 + 2m_2 + 3m_3 \equiv 2 \pmod 4$ and hence
\[
\frac{o(w)}{2} = \frac{m_1+m_3}{2} \equiv m_2 + m_3 + 1 \pmod 2.
\]
The claim follows.

We conclude that $h' = - \phi(h)$ where $\phi = \phi_1 \circ
\phi_2$. Note that $\phi$ is a basis change composed with the map
$t^{\pm 1} \mapsto (-1)^{r(\Lambda)} t^{\pm 1}$.
By construction, $\phi$ negates exactly one of each $p_j,q_j$ pair; it
follows that $\phi\sftb{x,y} = -\sftb{\phi(x),\phi(y)}$ and
$\phi(\delta(x)) = -\delta(\phi(x))$ for all $x,y$. Hence
\[
\phi\d\phi^{-1}(x) = \phi\sftb{h,\phi^{-1}(x)} +
\phi\delta\phi^{-1}(x)
= -\sftb{\phi(h),x} - \delta(x) = \d'(x)
\]
and this establishes the proposition.
\end{proof}

\secspace

\section*{Appendix B: Stabilized Knots}

\setcounter{figure}{0}
\renewcommand{\thefigure}{B.\arabic{figure}}
\setcounter{theorem}{0}
\renewcommand{\thetheorem}{B.\arabic{theorem}}

In this section, we show that the LSFT algebra of any stabilized
Legendrian knot $\Lambda$ is equivalent via a basis change to an LSFT algebra
that depends only on $tb(\Lambda)$ and $r(\Lambda)$, the classical
Legendrian invariants associated to $\Lambda$. This implies that the
LSFT algebra of a stabilized knot contains no interesting information
about the knot.

\begin{figure}
\centerline{
\includegraphics[height=1.5in]{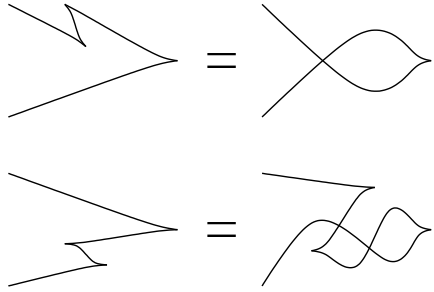}
}
\caption{
Obtaining a bubble from a stabilization, in the front
projection. Which of the two pictures (top or bottom) applies depends
on the sign of the stabilization and the orientation of the rightmost
cusp of the front.
}
\label{fig:stab1}
\end{figure}

\begin{figure}
\centerline{
\includegraphics[height=0.6in]{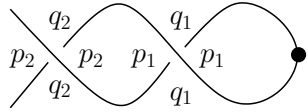}
}
\caption{
The $xy$ projection of a stabilized knot $\Lambda$.
}
\label{fig:stab2}
\end{figure}

Let $\Lambda$ be a stabilized knot, i.e., a knot Legendrian isotopic
to one whose front diagram contains a zigzag. Up to isotopy, we
can assume
that the front of $\Lambda$ has a zigzag next to its rightmost
cusp. We can further isotop the front to obtain a ``bubble'' at the
rightmost cusp; see Figure~\ref{fig:stab1}. The resolution of a bubble
is shown in Figure~\ref{fig:stab2}. It follows that, up to equivalence
of LSFT algebras, we can assume that the $xy$ diagram for $\Lambda$,
given by resolving its front, contains the piece shown in
Figure~\ref{fig:stab2}, and no part of the diagram lies further to the
right than the depicted part.

With $q_1,p_1,q_2,p_2$ as labeled and the base point $\bullet$ as
shown in Figure~\ref{fig:stab2}, the LSFT algebra for $\Lambda$
satisfies
\begin{align*}
F_{\d} &= p_1 \\
\d(p_1) &= 0 \\
\d(q_1) &= 1 - p_2 + q_1 p_1 q_1.
\end{align*}
(These signs are correct if the rightmost loop in
Figure~\ref{fig:stab2} is oriented counterclockwise; if it is oriented
clockwise, then $\d(q_1) = -1 + p_2 - q_1p_1q_1$, but this is
equivalent to the given signs after we replace $q_1$ by $-q_1$.)
By further conjugating by the basis change given by the elementary
automorphism $p_2 \mapsto -p_2 + q_1p_1q_1$, we can write $\d(q_1) =
1+p_2$ instead.

\begin{definition}
The LSFT algebra $(\alg,\d)$ generated by
$q_1,\dots,q_n,p_1,\dots,p_n,t^{\pm 1}$ is of \textit{ordered
stabilized type} if it has the form:
\begin{align*}
F_{\d} &= p_1 \\
\d(p_1) &= 0 \\
\d(q_1) &= 1 + p_2 \\
\d(p_2) &= [p_1,q_1] \\
\d(q_j) &= (\text{function of } t^{\pm 1},q_1,\dots,q_{j-1}) + O(p) \\
\d(p_j) &= (\text{function of } t^{\pm
  1},q_1,\dots,q_n,p_1,p_2,p_{j+1},\dots,p_n) + O(p^2) \\
\d(t) &= O(p)
\end{align*}
where $j=1,\dots,n$. \label{def:ordstab}
\end{definition}

\begin{lemma}
Up to equivalence, the LSFT algebra for a stabilized knot is of
ordered stabilized type.
\end{lemma}

\begin{proof}
From the discussion prior to Definition~\ref{def:ordstab}, we can
write the LSFT algebra of any stabilized knot to satisfy $F_{\d} =
p_1$, $\d(p_1) = 0$, $\d(q_1) = 1+p_2$, and hence $\d(p_2) =
[p_1,q_1]$ as well. Now order the crossings of the $xy$ projection to
satisfy $h(q_3) \leq h(q_4) \leq \cdots \leq h(q_n)$, where $h$ is the
height of the corresponding Reeb chord. The remainder of the
conditions in Definition~\ref{def:ordstab} then follow automatically
from Lemma~\ref{lem:stokes}.
\end{proof}

A \textit{linear map} $f\co\alg\to\alg$ is a map that sends $y\in\alg$
to a sum of terms, each of which includes $y$ once. Given a linear map
$f$ and $y\in\alg$, we can define another linear map
$D_yf\co\alg\to\alg$
as the ``derivative'' of $f$:
$(D_yf)(z) = \d(f(y))|_{\d(y) = z}$. Thus if $f(y) = x_1 y x_2$ for
$x_1,x_2\in\alg$, then
\[
(D_yf)(z) = \d(x_1) y x_2 + (-1)^{|x_1|} x_1 z x_2 + (-1)^{|x_1y|} x_1
y \d(x_2).
\]

Let $(\alg,\d)$ be of ordered stabilized type.
For $k \geq -1$, define a linear map $f_k \co \alg \to \alg$
inductively as follows: $f_{-1}(y) = f_0(y) = 0$ and
\[
f_k(y) = f_{k-1}(y) + q_1 [p_1,y] - q_1 (D_yf_{k-1})(f_{k-1}(y)) \pmod
{p^{k+1}},
\]
where $\pmod{p^{k+1}}$ indicates that we drop all terms of order $k+1$
or higher in the $p$'s (i.e., in $\F^{k+1}\!\alg$ for all $y$).
The first few $f_k$'s are given as follows:
\begin{align*}
f_0(y) &= 0 \\
f_1(y) &= q_1 [p_1,y] \\
f_2(y) &= q_1 [p_1,y] + q_1^2[p_1,q_1[p_1,y]] - q_1p_2[p_1,y] \\
f_3(y) &= q_1 [p_1,y] + q_1^2[p_1,q_1[p_1,y]] - q_1p_2[p_1,y] +
q_1^2 [p_1,q_1^2[p_1,q_1[p_1,y]]] \\
&\quad -q_1^3 [p_1,p_2[p_1,y]]
- q_1^2[p_1,q_1p_2[p_1,y]] -
q_1p_2q_1[p_1,q_1[p_1,y]] + q_1p_2^2[p_1,y].
\end{align*}

\begin{lemma}
For all $k$, $f_k(y)$ is a linear function in $y$ with coefficients
involving only $p_1,q_1,p_2$, and $f_k(y) \in \F^1\!\alg$ (i.e., every
term in $f_k(y)$ involves some $p$). Furthermore, $f_k(y)$ does not
depend on $(\alg,\d)$.
\label{lem:fk}
\end{lemma}

\begin{proof}
Clear by induction.
\end{proof}

Now suppose that $(\alg,\d)$ is of ordered stabilized type.
Let $s$ be a generator of $\alg$ besides $p_1,q_1,p_2$. We
inductively define a sequence of differentials $d_k$, elements
$x_k\in\alg$, and elementary
automorphisms $\phi_k$ supported at $s$, as follows:
\begin{itemize}
\item
$d_{-1} = d$;
\item
$x_k = d_{k-1}(s) - f_{k-1}(s) \pmod{p^{k+1}}$;
\item
$\phi_k(s) = s - q_1 x_k$;
\item
$d_k = \phi_k^{-1} \circ d_{k-1} \circ \phi_k$.
\end{itemize}
Note that it is not clear a priori that $\phi_k$ is an elementary
automorphism (in particular, invertible).

\begin{lemma}
For all $k \geq 0$, we have: \label{lem:dk}
\begin{itemize}
\item
$\d_k(x) = \d(x)$ for $x=p_1,q_1,p_2$ and $F_{d_k} = p_1$;
\item
$x_k = O(p^k)$;
\item
$\phi_k$ is an elementary automorphism;
\item
$d_k(s) = f_k(s) + O(p^{k+1})$.
\end{itemize}
\end{lemma}

\begin{proof}
We prove the lemma by induction. For $k=0$, since $d$ is ordered,
$x_0 = d(s) \pmod{p}$ does not involve $s$, so $\phi_0$ is
elementary. We then have $d(x_0) = d^2(s) + O(p) = [p_1,s] + O(p) =
O(p)$ and thus
\[
d_0(s) = \phi_0^{-1}(d(s-q_1x_0)) =
\phi_0^{-1}(x_0-x_0+q_1d(x_0)+O(p)) = \phi_0^{-1}(O(p)) = O(p).
\]

Now assume the lemma holds for $k-1 \geq 0$. Since $\phi_k$ is
supported at $s$ and $\d_{k-1}(x)$ does not involve $s$ for
$x=p_1,q_1,p_2$, it follows that $\d_k(x) = \d_{k-1}(x)$ for these
values of $x$, while $F_{\d_k} = \phi_k^{-1}(F_{\d_{k-1}}) = p_1$.

Next, since $d_{k-1}(s) =
f_{k-1}(s) +
O(p^k)$, $x_k = O(p^k)$ by definition. It follows that $\phi_k$ is
elementary. (If $s$ is a $p$ and $k=1$, then $x_0=0$ and $x_1 = d(s)
\pmod{p^2}$; since $\d$ is ordered, $x_1$ does not involve $s$ and so
$\phi_k$ is elementary.)

Finally we check that $d_k(s) = f_k(s) + O(p^{k+1})$. By Lemma~\ref{lem:fk},
$f_{k-1}(s)$ is linear in $s$ with coefficients involving only
$p_1,q_1,p_2$; since $d_{k-1}=d$ on $p_1,q_1,p_2$, we have
$d_{k-1}(f_{k-1}(s)) = (D_sf_{k-1})(d_{k-1}(s))$. By the induction
assumption, $d_{k-1}(s) \equiv f_{k-1}(s) \pmod{p^k}$; since every
term in $f_{k-1}$ involves a $p$ by Lemma~\ref{lem:fk} again, we
conclude that
\[
d_{k-1}(f_{k-1}(s)) = (D_sf_{k-1})(d_{k-1}(s)) \equiv
(D_sf_{k-1})(f_{k-1}(s)) \pmod{p^{k+1}}.
\]
Thus
\begin{align*}
d_{k-1}\phi_k(s) = d_{k-1}(s-q_1x_k) &=
d_{k-1}(s)-(1+p_2)x_k+q_1d_{k-1}(x_k) \\
&= f_{k-1}(s) + q_1(d_{k-1}^2(s)-d_{k-1}f_{k-1}(s)) + O(p^{k+1}) \\
&= f_{k-1}(s) + q_1[p_1,s] - q_1(D_sf_{k-1})(f_{k-1}(s)) + O(p^{k+1}) \\
&= f_k(s) + O(p^{k+1}).
\end{align*}
Now $\phi_k$ is the identity mod
$p^k$ and $f_k(s) \in \F^1\!\alg$,
whence $d_k(s) = \phi_k^{-1}d_{k-1}\phi_k(s) = \phi_k^{-1}(f_k(s)) +
O(p^{k+1}) = f_k(s) + O(p^{k+1})$.
\end{proof}

\begin{lemma}
For all $k$, we have $f_k(s) = f_{k-1}(s) + O(p^k)$. \label{lem:fkstab}
\end{lemma}

\begin{proof}
This follows directly from Lemma~\ref{lem:dk}. Since $\phi_k$ is the
identity mod $p^k$, $d_k(s) = d_{k-1}(s) + O(p^k)$. From
Lemma~\ref{lem:dk} again, $f_k(s) = d_k(s) + O(p^{k+1})$ and
$f_{k-1}(s) = d_{k-1}(s) + O(p^k)$. The lemma follows.
\end{proof}

Because of Lemma~\ref{lem:fkstab}, we can define an element $f(s) \in
\alg$ to be the limit $\lim_{k\to\infty} f_k(s)$.

\begin{lemma}
If $(\alg,\d)$ is of ordered stabilized type and $s$ is one of
$q_2,\dots,q_n,p_3,\dots,p_n,t$, then there is an
elementary automorphism $\phi$ supported at $s$ such that if $\d' =
\phi^{-1} \circ d \circ \phi$, then $(\alg,\d')$ is of ordered
stabilized type and $\d'(s) = f(s)$.
\label{lem:ordstab}
\end{lemma}

\begin{proof}
Using Lemma~\ref{lem:dk},
set $\phi = \lim_{k\to\infty} (\phi_0 \circ \phi_1 \circ \cdots \circ
\phi_k)$, a well-defined limit since $\phi_k$ is the identity mod
$p^k$ for all $k$.
\end{proof}

\begin{proposition}
An LSFT algebra of ordered stabilized type is equivalent under basis
change to the LSFT algebra whose differential is given by
\label{prop:stab}
\begin{align*}
F_{\d} &= p_1 \\
\d(p_1) &= 0 \\
\d(q_1) &= 1 + p_2 \\
\d(p_2) &= [p_1,q_1] \\
\d(q_j) &= f(q_j) \text{ for $2\leq j\leq n$} \\
\d(p_j) &= f(p_j) \text{ for $3 \leq j \leq n$} \\
\d(t) &= f(t).
\end{align*}
\end{proposition}

\begin{proof}
Successively apply Lemma~\ref{lem:ordstab} with
$s=q_2,\dots,q_n,p_3,\dots,p_n,t$.
\end{proof}

For any Legendrian knot $\Lambda$,
the LSFT algebra of $\Lambda$ encodes both
the Thurston--Bennequin and rotation numbers of $\Lambda$:
$r(\Lambda)$ is $-1/2$ times the degree of $t$, while $tb(\Lambda)$ is
the difference between the number of $q$ generators of even degree and
the number of odd degree. However, if $\Lambda$ is
stabilized, Proposition~\ref{prop:stab} implies that the LSFT algebra
encodes nothing else.

\begin{corollary}
The LSFT algebra of a stabilized knot $\Lambda$
is equivalent to an LSFT algebra
depending only on $tb(\Lambda)$ and $r(\Lambda)$.
\end{corollary}

\secspace

\bibliographystyle{alpha}
\bibliography{biblio}

\end{document}